\documentclass{article}
\usepackage{xcolor}
\usepackage{amsmath}
\usepackage{amssymb}
\usepackage{amsthm}
\usepackage{hyperref}
\hypersetup{
	colorlinks=true,
	linkcolor=cyan,
	filecolor=blue,
	urlcolor=red,
	citecolor=green,
}

\newtheorem{proposition}{Proposition}[section]
\newtheorem{theorem}{Theorem}[section]
\newtheorem{lemma}{Lemma}[section]

\title{Optimal edge fault-tolerant-prescribed hamiltonian laceability of balanced hypercubes}
\author{Ningning Song and Yuxing Yang
\footnote{They are in School of Mathematics and Information Science,
Henan Normal University,
Xinxiang, Henan 453007, China.
E-mails:\ yxyangcn@163.com,\ yyx@htu.edu.cn}}
\date{}

\begin{document}
\maketitle

\begin{abstract}
{\bf Aims:}
Try to prove the $n$-dimensional balanced hypercube $BH_n$ is $(2n-2)$-fault-tolerant-prescribed hamiltonian laceability.
{\bf Methods:}
Prove it by induction on $n$. It is known that the assertation holds for $n\in\{1,2\}$. Assume it holds for $n-1$ and prove it holds for $n$, where $n\geq 3$. If there are $2n-3$ faulty links and they are all incident with a common node, then we choose some dimension such that there is one or two faulty links and no prescribed link in this dimension; Otherwise, we choose some dimension such that the total number of faulty links and prescribed links does not exceed $1$. No matter which case, partition $BH_n$ into $4$ disjoint copies of $BH_{n-1}$ along the above chosen dimension.
{\bf Results:}
On the basis of the above partition of $BH_n$, in this manuscript, we complete the proof for the case that there is at most one faulty link in the above chosen dimension. 
\end{abstract}

\section{Introduction}
\label{section1}

The $n$-dimensional balanced hypercube $BH_n$ was proposed by Huang and Wu \cite{Huang1995} as a variant of the well-known hypercube, and it has most of the good properties of the hypercube, such as bipartite structure, recursiveness, regularity, vertex-symmetry \cite{WuHuang1997} and edge-symmetry \cite{Zhou2015}. Particularly, each vertex of the balanced hypercube has a backup vertex that has the same neighborhood as the original one. Let $(X,Y)$ be a bipartition of $BH_{n}$.

Cheng et al. \cite{Cheng2014} investigated the disjoint paths cover problem of balanced hypercubes, and they proved the following.

\begin{theorem} (see \cite{Cheng2014})\label{th-cheng2014}
Let $u, x\in X$ and $v, y\in Y$ be pairwise distinct. Then there exist two vertex-disjoint paths $P[u,v]$ and $P[x,y]$ in $BH_{n}$ such that each vertex of $BH_n$ lies on one of the two paths.
\end{theorem}

The hamiltonian property is a major requirement in designing network topologies since a topology structure containing hamiltonian paths or cycles can efficiently simulate many algorithms designed on linear arrays or rings (see for example, \cite{Chen2007,Xu2007} and references therein). A bipartite graph is \emph{hamiltonian laceable} if there is a hamiltonian path between any two vertices in different bipartite sets. Xu et al. \cite{Xu2007} investigated the hamiltonian laceability of  balanced hypercubes and they obtained the following.

\begin{theorem}(see \cite{Xu2007})\label{th-xu2007}
$BH_n$ is hamiltonian laceable.
\end{theorem}

In parallel computer systems, failures of processors and/or physical links are inevitable. Thus, the problem of fault-tolerant embedding of hamiltonian paths and cycles has become an important issue and has been studied in depth (see, for example \cite{WuHuang1997,LvZhang2014, Hao2014,Zhou2015,Cheng2015,Li2019}). For any set $F$ of at most $k$ edges of a bipartite graph $G$, if $G-F$ is hamiltonian laceable, then $G$ is said to be \emph{$k$-fault-tolerant hamiltonian laceable}. Zhou at el. \cite{Zhou2015} investigated fault-tolerant hamiltonian laceability of balanced hypercubes. One of their main results can be restated as follows:

\begin{theorem} (see \cite{Zhou2015})\label{th-zhou2015}
$BH_n$ is $(2n-2)$-fault-tolerant hamiltonian laceable.
\end{theorem}

In \cite{Li2019}, Li et al. investigated the problem of embedding hamiltonian cycle into balanced hypercubes with conditional faulty edges and they obtained the following.

\begin{theorem} (see \cite{Li2019})\label{th-li2019}
Let $F\subset E(BH_{n})$ with $|F|\leq 4n-5$ such that the minimum degree of $BH_n-F$ is at least $2$. Then each edge in $BH_{n}-F$ lies on a hamiltonian cycle of $BH_n-F$.
\end{theorem}

In \cite{LvZhang2014}, L\"{u} and Zhang proved the following result on the problem of embedding hamiltonian paths into $BH_n$ with a faulty vertex.

\begin{theorem} (see \cite{LvZhang2014})\label{th-lv2014}
Let $u\in X$ be a vertex of $BH_{n}$, and let $x,y\in Y$. Then there is hamiltonian path of $BH_n-u$ connecting $x$ and $y$.
\end{theorem}

As a complementary to fault-tolerant embedding problem, Dvo\v{r}\'{a}k \cite{Dvorak2005} proposed the prescribed embedding problem which requires that the embedded paths and cycles pass through a given number of prescribed edges. Following Dvo\v{r}\'{a}k's work, prescribed embedding problems were studied in literatures (see, for example, \cite{Caha2006,Chen2017,Cheng2019,Lv2019} and references therein). A set $\{u,v\}$ of two vertices in a graph $G$ \emph{is compatible to} a given linear forest $L$ of $G$ if none of the paths in $L$ has $u$ or $v$ as internal vertices or both of them as end vertices. A bipartite graph $G$ is \emph{$k$-prescribed hamiltonian laceable} if $G$ admits a hamiltonian path between $u$ and $v$ passing through any prescribed linear forest $L$ with at most $k$ edges provided that $\{u,v\}$ is compatible to $L$.
Cheng \cite{Cheng2019} investigated prescribed hamiltonian laceability of balanced hypercubes and she obtained the following.

\begin{theorem} (see \cite{Cheng2019})\label{th-cheng2019}
$BH_n$ is $(2n-2)$-prescribed hamiltonian laceable.
\end{theorem}

In faulty interconnection networks, the embedded  fault-free paths and/or cycles may be required to pass through a prescribed linear forest.
A bipartite graph $G$ is \emph{$k$-fault-tolerant-prescribed hamiltonian laceable} if $G-F$ is $(k-|F|)$-prescribed hamiltonian laceable for any set $F$ with at most $k$ edges in $G$.
In \cite{Yang2019}, Yang and Zhang investigated fault-tolerate-prescribed hamiltonian laceability of balanced hypercubes and they proved the following.

\begin{theorem} (see \cite{Yang2019})\label{th-yang2019}
$BH_n$ is $(n-1)$-fault-tolerant-prescribed hamiltonian laceable.
\end{theorem}

Inspired by the above works, in this manuscript, we try to prove the following and we complete the proof of the major case.

\begin{theorem}\label{th-main}
$BH_n$ is $(n-1)$-fault-tolerant-prescribed hamiltonian laceable. 
\end{theorem}

\section{Preliminaries}
\label{section2}

The \emph{neighborhood} $N_G(v)$ of a vertex $v$ in a graph $G$ is the set of neighbors of $v$ in $G$. Denote by $P[u,v]$ a path between $u$ and $v$, and abbreviate the terms ``hamiltonian path" and ``hamiltonian cycle" as ``H-path" and ``H-cycle", respectively. \emph{A maximal path} is one that can not be extended to a longer path from either end. For notations and operations used without defining here we follow \cite{Bondy2008}. Denote by $N_k$ the set of non-negative integers less than $k$ for any positive integer $k$. In the rest of the paper, all the additions and subtractions on the superscript and subscript of a symbol are modulo $4$.

The \emph{$n$-dimensional balanced hypercube} $BH_{n}$ is a simple graph that consists of $4^{n}$ vertices, and each of which is labelled by $x=x_{0}x_{1}\cdots x_{n-1}$, where $x_{i}\in N_{4}$ for any $i\in N_n$. A vertex $\alpha=\alpha_{0}\alpha_{1}\cdot\cdot\cdot \alpha_{n-1}\in BH_{n}$ has $2n$ neighbors $\alpha^{0\pm}, \alpha^{j\pm}$ of $\alpha$ in $BH_{n}$, where $\alpha^{0\pm}=(\alpha_{0}\pm 1\mod 4)\alpha_{1}\alpha_{2}\cdot\cdot\cdot \alpha_{n-1}$, and $\alpha^{j\pm}=(\alpha_{0}\pm 1\mod 4)\alpha_{1}\cdot\cdot\cdot (\alpha_{j}+(-1)^{\alpha_{0}}\mod 4)\alpha_{j+1}\cdot\cdot\cdot \alpha_{n-1}$ for $j\in N_{n}\setminus \{0\}$. We call $(\alpha,\alpha^{i\pm})\in E(BH_{n})$ \emph{$i$-dimensional edges}, for $i\in N_{n}$. The \emph{shadow vertex} $\alpha^{s}=(\alpha^{i+})^{i+}=(\alpha^{i-})^{i-}=(\alpha_{0}+2)\alpha_{1}\cdots \alpha_{n-1}$ of $\alpha$ is unique, and $\alpha$ and $\alpha^{s}$ have the same neighbor-set. Clearly, $(\alpha^{i+})^s=\alpha^{i-}$.
$BH_{n}$ has a recursive structure, more precisely, for $n\geq 2$, $BH_{n}$ can be partitioned into $4$ disjoint copies of $BH_{n-1}$ along some dimension $d^*\in N_{n}$ by deleting all the $d^*$ dimensional edges of $BH_{n}$ \cite{Lv2012}.

\begin{lemma} (see Section 3 in \cite{Yang-Song})\label{le-yang-song}
Let $f$ and $e$ be any two different edges in $BH_2$. For any two vertices $u$ and $v$ in different parts of $BH_2$, $BH_2-f$ admits a hamiltonian path passing through $e$.
\end{lemma}

\begin{lemma}\label{le-base}
{$BH_2$ is $2$-fault-tolerant-prescribed hamiltonian laceable.}
\end{lemma}

\begin{proof}
By Theorems \ref{th-zhou2015}, \ref{th-cheng2019} and Lemma \ref{le-yang-song}, the lemma is immediate.
\end{proof}

In the rest of the paper, we try to prove Theorem \ref{th-main}. Let $F\subset E(BH_{n})$ be a set of faulty edges and $L$ be linear forest in $BH_{n}-F$ such that $|E(L)\cup F|\leq 2n-2$. Let $u,v$ be two vertices in opposite partite set of $BH_n$ such that $\{u,v\}$ is compatible to $L$.
It suffices to prove that $BH_n-F$ admits a hamiltonian path between $u$ and $v$ passing through $L$, and it is enough to consider the case that the total number of the edges in $F$ and $L$ is up to $2n-2$. Theorems \ref{th-zhou2015} and \ref{th-cheng2019} imply the result holds for $E(L)=\emptyset$ and $F=\emptyset$, respectively. In the following, we consider the case that $E(L)\neq \emptyset$ and $F\neq \emptyset$. Prove the result by induction on $n$. The result holds trivially for $BH_1$. Lemma \ref{le-base} implies that $BH_2$ is $2$-fault-tolerant-prescribed hamiltonian laceable. In the remainder, we will assume that the result holds for $BH_{n-1}$ and prove it also holds for $BH_n$ for $n\geq 3$. We partition $BH_n$ into $4$ copies of $BH_{n-1}$ along some dimension according to the following rules.

If $|F|=2n-3$ and all of the faulty edges are incident to a common vertex, then there is exactly one edge $e$ in $L$. Assume that $e$ is an $i$-dimensional edge for some $i\in N_{n}$. Then by the Pigeonhole Principle, there exists some $j\in N_n\setminus \{i\}$ such that $F$ has at least one edge in dimension $j$. Clearly, $F$ has at most $2$ edges in dimension $j$, and $L$ has no edge in dimension $j$.

If $|F|\leq 2n-4$, or $|F|=2n-3$ and not all of the faulty edges are incident to a common vertex, then there exists some $j\in N_n$ such that there is at most one edge of $F\cup E(L)$ in this dimension.

No matter which case above, we can assume that $j=n-1$ and partition $BH_{n}$ into $4$ disjoint copies, $B^{0},B^{1},B^{2},B^{3}$, of $BH_{n-1}$ along $n-1$ dimension, where the rightmost digit of any vertex in $B^{i}$ is $i$ for $i\in N_{4}$. For simplification, abbreviate $V(B^{i})$ as $V_{i}$. Denote by $L_{i}$ and $F_{i}$ the restriction of $L$ and $F$ in $B^{i}$, respectively. Without loss of generality, assume that $|E(L_{0})\cup F_{0}|=\max\{|E(L_{i})\cup F_{i}|: i\in N_{4}\}$.
Denote by $E_{i,j}$ the set of edges between $B^{i}$ and $B^{j}$, and denote $L^c=E(L)\cap E_{i,j}, F^{c}=F\cap E_{i,j}$, where $i,j\in N_{4}$ and $i\neq j$. For an arbitrary vertex $x\in V_{i}$, abbreviate the neighbors $x^{(n-1)\pm}$ of $x$ as $x^{\pm}$. The fact that for any two distinct vertices $a,b\in V_{i}$, then $a^{+}\neq b^{+}$, and $a^{-}\neq b^{-}$ will be used often in the remainder, where $i\in N_{4}$. A vertex in $BH_n$ is an \emph{even vertex} (resp. \emph{odd vertex}) if the leftmost digit of which is even (resp. odd). Let $X$ and $Y$ are the sets of even vertices and odd vertices in $BH_n$, respectively. Then $(X,Y)$ is a bipartition of $BH_{n}$. Without loss of generality, assume that $u\in X$ and $v\in Y$.

On the basis of the above way that we partition $BH_n$, there are four cases to consider, i.e., the cases $F_c=E(L_c)=\emptyset$, $F_c=\emptyset$ and $|E(L_c)|=1$, $|F_c|=1$ and $E(L_c)=\emptyset$, and $|F_c|=2$ and $E(L_c)=\emptyset$. Sections \ref{section3}, \ref{section4} and \ref{section5} will deal with the former cases.

The following lemmas will be used in the proof of our main result in \ref{section3}, \ref{section4}, \ref{section5}.

\begin{lemma}\label{le-1}
If $|E(L_{0})\cup F_{0}|=2n-3$ and $E(L_{0})\neq \emptyset$,
then $B^{0}-F_{0}$ contains a H-cycle passing through $L_{0}$.
\end{lemma}

\begin{proof}
Let $(x,y)\in E(L_{0})$. Then $\{x,y\}$ is compatible to $L_{0}-(x,y)$. Since $B^{0}\cong BH_{n-1}$, by the induction hypothesis, $B^{0}-F_{0}$ has a H-path $P[x,y]$ passing through $L_{0}-(x,y)$. Hence, $P[x,y]+(x,y)$ is a H-cycle passing through $L_{0}$ in $B^{0}-F_{0}$.
\end{proof}

\begin{lemma}\label{le-2}
Given a $s\in V_i\cap X$ (resp. $s\in V_i\cap Y$). Let $i\in N_4$. If $|E(L_{0})\cup F_{0}|\leq 2n-4$, then there is a vertex $x\in V_{i}\cap X$ (resp. $x\in V_{i}\cap Y$) such that

\noindent ($i$). $x$ is incident with none of $E(L_{i})$; and

\noindent ($ii$). none of $x^{\pm}$ is incident with an edge of $E(L_{i+1})\cup F_{i+1}$ (resp. $E(L_{i-1})\cup F_{i-1}$); and

\noindent ($iii$). $x\neq s$.

\end{lemma}

\begin{proof}
The proofs of the cases $x\in X$ and $x\in Y$ are analogous. We here only consider the case $x\in X$. A vertex $x\in V_{i}\cap X\setminus \{s\}$ fails the lemma only if

\noindent ($a$). $x$ is incident with an edge of $L_{i}$; or

\noindent ($b$). $x^{+}$ or $x^{-}$ is incident with an edge of $E(L_{i+1})\cup F_{i+1}$; or


There are $|V_{i}\cap X\setminus \{s\}|=4^{n-1}/{2}-1$ vertex candidates.
Since there are at most $|E(L_i)|$ even vertices in $L_i$, the number of such $x$ that supports $(a)$ does not exceed $|E(L_i)|$.
Since there are at most $|E(L_{i+1})|+|F_{i+1}|$ odd vertices incident with an edge of $E(L_{i+1})\cup F_{i+1}$, each of which makes at most two vertex candidates support $(b)$, the number of such $x$ that supports $(b)$ does not exceed $2(|E(L_{i+1})|+|F_{i+1}|)$. 
Thus, the total number of vertex candidates that fail the lemma does not exceed $|E(L_i)|+2|E(L_{i+1})\cup F_{i+1}|+1\leq |E(L)\cup F|+|E(L_{i+1})\cup F_{i+1}|\leq (2n-2)+(2n-4)\leq 4n-6$. Since $|V_{i}\cap X\setminus \{s\}|-(4n-5)=(4^{n-1}/{2}-1)-(4n-6)>0$ for $n\geq 3$, there is an $x\in V_{i}\cap X$ supporting the lemma.
\end{proof}

\begin{lemma}\label{le-3}
Given a $y\in V_i$. Let $i\in N_4$ and let $P[z,w]$ be a H-path of $B^{i}-F_{i}$ passing through $L_{i}$. If $|E(L_{i})\cup F_{i}|\leq 2n-3$. Then there is an edge $(s,t)\in E(P[z,w])\setminus E(L_{i})$ for some $s\in X$ and $t\in Y$ such that


\noindent $(i)$. $\{s,t\}\cap \{z,w\}=\emptyset$; and

\noindent $(ii)$. if $|E(L_{i})\cup F_{i}|\leq 2n-4$, then $s^{+}$ or $s^{-}$ (resp. $t^{+}$ or $t^{-}$) is incident with none of $E(L_{i+1})\cup F_{i+1}$ (resp. $E(L_{i-1})\cup F_{i-1}$); and

\noindent $(iii)$. if $|E(L_{i})\cup F_{i}|=2n-3$, then $s^{\pm}$
(resp. $t^{\pm}$) are incident with none of $E(L_{i+1})\cup F_{i+1}$ (resp. $E(L_{i-1})\cup F_{i-1}$); and

\noindent $(\romannumeral4)$. $y\notin \{s,t\}$.

\end{lemma}

\begin{proof}

An edge $(s,t)\in E(P[z,w])\setminus E(L_{i})$ fails the lemma only if


\noindent ($a$). $\{s,t\}\cap \{z,w\}\neq \emptyset$; or

\noindent ($b$). if $|E(L_{i})\cup F_{i}|\leq 2n-4$, then both $s^{+}$ and $s^{-}$ (resp. both $t^{+}$ and $t^{-}$) are incident with an edge of $E(L_{i+1})\cup F_{i+1}$ (resp. $E(L_{i-1})\cup F_{i-1}$); or

\noindent ($c$). if $|E(L_{i})\cup F_{i}|=2n-3$, then $s^{+}$ or $s^{-}$, and $t^{+}$ or $t^{-}$ are incident with an edge of $E(L_{i+1})\cup F_{i+1}$ and $E(L_{i-1})\cup F_{i-1}$, respectively.

\noindent ($d$). $y\in \{s,t\}$.

There are $|E(P[z,w])|-|E(L_{i})|$ edge candidates.
Clearly, the number of such $(s,t)$ that supports $(a)$ and $(d)$ does not exceed $2+2=4$.

Suppose first that $|E(L_{i})\cup F_{i}|\leq 2n-4$.
Since there are at most $|E(L_{i+1})|+|F_{i+1}|$ (resp. $|E(L_{i-1})|+|F_{i-1}|$) odd (resp. even) vertices incident with an edge of $E(L_{i+1})\cup F_{i+1}$ (resp. $E(L_{i-1})\cup F_{i-1}$), each of which makes at most two edge candidates support $(b)$, the number of such $(s,t)$ that supports $(b)$ does not exceed $2(|E(L_{i+1})|+|F_{i+1}|)+2(|E(L_{i-1})|+|F_{i-1}|)$. Thus, the total number of edge candidates that fail the lemma does not exceed $2|E(L_{i+1})\cup F_{i+1}|+2|E(L_{i-1})\cup F_{i-1}|+4$. Since $|E(P[z,w])|-|E(L_{i})|-(2|E(L_{i+1})\cup F_{i+1}|+2|E(L_{i-1})\cup F_{i-1}|+4)\geq |E(P[z,w])|-(2|E(L)\cup F|+4)\geq 4^{n-1}-1-4n>0$, there is an edge $(s,t)\in E(P[z,w])$ supporting the lemma.

Suppose now that $|E(L_{i})\cup F_{i}|=2n-3$. The number of such $(s,t)$ that supports $(c)$ does not exceed $4(|E(L_{i+1})|+|F_{i+1}|)+4(|E(L_{i-1})|+|F_{i-1}|)$. Note that $(|E(L_{i+1})|+|F_{i+1}|)+(|E(L_{i-1})|+|F_{i-1}|)\leq |E(L)\cup F|-|E(L_{i})\cup F_{i}|\leq 1$.
Thus, the total number of edge candidates that fail the lemma does not exceed $4|E(L_{i+1})\cup F_{i+1}|+4|E(L_{i-1})\cup F_{i-1}|+4\leq 4+4=8$. Since $|E(P[z,w])|-|E(L_{i})|-8\geq 4^{n-1}-1-(2n-3)-8>0$, there is an edge $(s,t)\in E(P[z,w])$ supporting the lemma.
\end{proof}


\begin{lemma}\label{le-4}
Let $i\in N_4$. If $|E(L_{i})\cup F_{i}|\leq 2n-4$ and {$|E(L_{i+1})\cup F_{i+1}|\leq 2n-6$ (resp. $|E(L_{i-1})\cup F_{i-1}|\leq 2n-6$)}, then there is an even (resp. odd) vertex $s\in V_{i}$ such that

\noindent $(i)$. $s$ is incident with none of $E(L_{i})$; and

\noindent $(ii)$. neither $s^{+}$ nor $s^{-}$ is incident with an edge of $E(L_{i+1})\cup F_{i+1}$ (resp. $E(L_{i-1})\cup F_{i-1}$); and

\noindent $(iii)$. $u$ (resp. $v$) is not adjacent to $s^{\pm}$ in $B^{i+1}$ (resp. $B^{i-1}$); and

{\noindent $(\romannumeral4)$. furthermore, for $n\geq 4$, if $L^c\cup F^c=\{(x,y)\}$ for some $x\in X$ and $y\in Y$, then $s\notin \{x,y\}$ and $x$ (resp. $y$) is not adjacent to $s^{\pm}$ in $B^{i+1}$ (resp. $B^{i-1}$).}
\end{lemma}

\begin{proof}
The proofs of the cases $s\in X$ and $s\in Y$ are analogous. We here only consider the case $s\in X$. A vertex $s\in V_{i}\cap X$ fails the lemma only if

\noindent ($a$). $s$ is incident with an edge of $E(L_{i})$; or

\noindent ($b$). $s^{+}$ or $s^{-}$ is incident with an edge of $E(L_{i+1})\cup F_{i+1}$; or

\noindent ($c$). $(s^{+},u)\in E(B^{i+1})$ or $(s^{-},u)\in E(B^{i+1})$; or 

\noindent ($d$). $(s^{+},x)\in E(B^{i+1})$ or $(s^{-},x)\in E(B^{i+1})$.

There are $|V_{i}\cap X|=4^{n-1}/{2}$ vertex candidates. Since there are at most $|E(L_{i})\cup F_i|$ even vertices in $L_i$ and $F_i$, the number of such $s$ that supports $(a)$ does not exceed $|E(L_{i})\cup F_i|$. Since there are at most $|E(L_{i+1})|+|F_{i+1}|$ odd vertices incident with at least one edge of $E(L_{i+1})\cup F_{i+1}$, each of which makes at most two vertex candidates support $(b)$, the number of such $s$ that supports $(b)$ does not exceed $2(|E(L_{i+1})|+|F_{i+1}|)$. Clearly, the number of such $s$ that supports $(c)$ does not exceed $2|N_{B^{i+1}}(u)|/{2}$.

Suppose first that the condition of \noindent $(\romannumeral4)$ holds (i.e., $L^c\cup F^c=\{(x,y)\}$). The number of such $s$ that supports $(d)$ does not exceed $2|N_{B^{i+1}}(x)|/{2}$. Thus, the total number of vertex candidates that fail the lemma does not exceed $|E(L_{i})|+2|E(L_{i+1})\cup F_{i+1}|+2|N_{B^{i+1}}(u)|/{2}+2|N_{B^{i+1}}(x)|/{2}\leq |E(L_{i})\cup F_{i}|+2|E(L_{i+1})\cup F_{i+1}|+|N_{B^{i+1}}(u)|+|N_{B^{i+1}}(x)|\leq 3(2n-4)+2(2n-2)<10n-16$. Since $|V_{i}\cap X|-(10n-16)=4^{n-1}/{2}-(10n-16)>0$ for $n\geq 4$, there is a vertex $s$ supporting the lemma.

{Suppose now that the condition of \noindent $(\romannumeral4)$ does not hold.

If $n=3$, $|E(L_{i+1})\cup F_{i+1}|\leq 2n-6=0$, the total number of vertex candidates that fail the lemma does not exceed $|E(L_{i})|+2|E(L_{i+1})\cup F_{i+1}|+2|N_{B^{i+1}}(u)|/{2}\leq |E(L_i)|+0+|N_{B^{i+1}}(u)|\leq (2n-4)+(2n-2)=4n-6$. Since $|V_{i}\cap X|-(4n-6)=4^{n-1}/{2}-(4n-6)>0$, there is a vertex $s$ supporting the lemma.

If $n\geq 4$, the total number of vertex candidates that fail the lemma does not exceed $|E(L_{i})|+2|E(L_{i+1})\cup F_{i+1}|+2|N_{B^{i+1}}(u)|/{2}\leq |E(L)\cup F|+|E(L_{i+1})\cup F_{i+1}|+|N_{B^{i+1}}(u)|\leq (2n-2)+(2n-6)+(2n-2)=6n-10$. Since $|V_{i}\cap X|-(6n-10)=4^{n-1}/{2}-(6n-8)>0$, there is a vertex $s$ supporting the lemma.}
\end{proof}

\begin{lemma}\label{le-9}

Let $i\in N_4$ and let $r\in V_i\cap X$ (resp. $r\in V_i\cap Y$) such that

\noindent $(1)$. $r$ is incident with none of $E(L_{i})\cup F_{i}$; and

\noindent $(2)$. $v$ (resp. $u$) is not adjacent to $r$ in $B^i$; and

\noindent $(3)$. if $F^c=\emptyset$ and $L^c=\{(x,y)\}$ for some $x\in X$ and $y\in Y$, $y$ (resp. $x$) is not adjacent to $r$ in $B^i$.

\noindent If $|E(L_{i})\cup F_{i}|\leq 2n-5$, then $r$ has two neighbors $s$ and $t$ in $B^{i}$ such that

\noindent $(i)$. $L_{i}+\{(r,s),(r,t)\}$ is a linear forest; and

\noindent $(ii)$. $s^{+}$ or $s^{-}$ is incident with none of $E(L_{i-1})$ (resp. $E(L_{i+1})$); and

\noindent $(iii)$. $t^{+}$ or $t^{-}$ is incident with none of $E(L_{i-1})$ (resp. $E(L_{i+1})$).

\end{lemma}

\begin{proof}
The proofs of the cases $r\in V_i\cap X$ and $r\in V_i\cap Y$ are analogous. We here only consider the case $r\in V_i\cap X$.

A vertex $s\in N_{B^{i}}(r)$ fails the lemma only if

\noindent ($a$). $s$ is incident with an edge of $E(L_{i})$; or

\noindent ($b$). $s^{\pm}$ are incident with an edge of $L_{i-1}$.

There are $|N_{B^{i}}(r)|=2n-2$ vertex candidates.
Since there are at most $|E(L_i)|$ odd vertices in $L_i$, the number of such $s$ that supports $(a)$ does not exceed $|E(L_i)|$. Let $H$ be the set of even vertices which are not singletons in $L_{i-1}$. Then $|H|\leq |E(L_{i-1})|$. For two distinct $z,w\in H$, if $z$ is the shadow vertex of $w$, then the two vertices $z^+$ (i.e., $w^-$) and $z^-$ (i.e., $w^+$) support $(b)$. Therefore, the $|H|$ vertices in $H$ will make at most $|H|$ vertices of $N_{B^{i}}(r)$ support $(b)$.
Thus, the total number of such $s\in N_{B^{i}}(r)$ failing the lemma does not exceed $|E(L_i)|+|H|\leq |E(L_i)|+ |E(L_{i-1})|\leq |E(L)|+(|F|-1)\leq 2n-3$. Since $|N_{B^{i}}(r)|-(2n-3)=(2n-2)-(2n-3)>0$, there is a vertex $s\in N_{B^{i}}(r)$ supporting the lemma.

A vertex $t\in N_{B^{i}}(r)\setminus \{s\}$ fails the lemma only if

\noindent ($c$). $t$ is an internal vertex of $L_{i}$; or

\noindent ($d$). $t^{\pm}$ are incident with an edge of $L_{i-1}$.

Since there are at most $\lceil|E(L_i)|-1\rceil/{2}$ odd internal vertices in $L_i$, the number of such $t$ that supports $(c)$ does not exceed $\lceil|E(L_i)|-1\rceil/{2}$.
Similarly to the computation of such $s$ that supports $(b)$, we can obtain that the number of such $t$ that supports $(d)$ does not exceed $|E(L_{i-1})|$.
Thus, the total number of such $t\in N_{B^{i}}(r)\setminus \{s\}$ failing the lemma does not exceed $\lceil|E(L_i)|-1\rceil{2}+|E(L_{i-1})|\leq |E(L_i)|/{2}+|E(L_{i-1})|\leq (|E(L)\cup F|-1)/{2}+|E(L_{i-1})\cup F_{i-1}|/{2}\leq (2n-2)-1/{2}+2n-5/{2}\leq 2n-4$. Since $|N_{B^{i}}(r)\setminus \{s\}|-(2n-4)=(2n-3)-(2n-4)>0$, then there is a $t\in N_{B^{i}}(r)\setminus \{s\}$ supporting the lemma.
\end{proof}

\section{$L^c=F^c=\emptyset$.}
\label{section3}

\begin{proposition}\label{pr-1}
If $|E(L_{0})\cup F_{0}|=2n-2$, then $B^{0}-F_{0}$ contains a H-path $P[a,b]$ passing through $L_{0}$ for some $a\in X$ and $b\in Y$.
\end{proposition}

\begin{proof}
Since $F_{0}\neq \emptyset$, there is an edge $f\in F_{0}$. By Lemma \ref{le-1}, $B^{0}-F_{0}\setminus \{f\}$ has a H-cycle $C_{0}$ passing through $L_{0}$. Let $(a,b)=f$ if $f$ lies on $C_{0}$, and let $(a,b)$ be an arbitrary edge in $C_{0}\setminus E(L_{0})$ otherwise. Then $C_0-(a,b)$ is a desired path.
\end{proof}

\begin{lemma}\label{main-1}

If $u,v\in V_i$ for some $i\in N_{4}$, then $BH_{n}-F$ has a H-path $P[u,v]$ passing through $L$.
\end{lemma}

\begin{proof}
According to the total number of edges in $L_0$ and $F_0$, we consider the following three cases.

{\it Case 1.}  $|E(L_{0})\cup F_{0}|\leq 2n-4$.

Since $|E(L_{0})\cup F_{0}|=\max\{|E(L_{k})\cup F_{k}|: k\in N_{4}\}$, then $|E(L_{k})\cup F_{k}|\leq 2n-4$ for $k\in N_{4}$. In this scenario, the proofs of the cases $i=0$, $i=1$, $i=2$ and $i=3$ are almost the same. We here only consider the case $i=0$.

Since $|E(L_{0})\cup F_{0}|\leq 2n-4$, by the induction hypothesis, there is a H-path $P[u,v]$ passing through $L_{0}$ in $B^{0}-F_{0}$. Lemma \ref{le-3} implies that there is an edge $(a,b)\in P[u,v]\setminus E(L_{0})$ for some $a\in X$ and $b\in Y$ such that $a^+$ or $a^-$ (resp. $b^+$ or $b^-$), say $a^+$ (resp. $b^+$), is not incident with an edge of $L_{1}$ (resp. $L_3$). By Lemma \ref{le-2}, there is an $x\in V_{1}\cap X$ such that $x$ (resp. $x^+$) is not incident with an edge of $L_1$ (resp. $L_2$). Again by Lemma \ref{le-2}, there is a $y\in V_{2}\cap X$ such that $y$ (resp. $y^+$) is not incident with an edge of $L_2$ (resp. $L_3$). Thus, $\{a^{+},x\}$ is compatible to $L_{1}$, $\{x^{+},y\}$ is compatible to $L_{2}$, and $\{y^{+},b^+\}$ is compatible to $L_{3}$. Combining these with $|E(L_{k})\cup F_{k}|\leq 2n-4$ for $k\in N_{4}$, there are H-paths $P[a^{+},x]$ passing through $L_{1}$ in $B^{1}-F_{1}$, $P[x^{+},y]$ passing through $L_{2}$ in $B^{2}-F_{2}$, and $P[y^{+},b^{+}]$ passing through $L_{3}$ in $B^{3}-F_{3}$. Hence $P[u,v]\cup P[a^{+},x]\cup P[x^{+},y]\cup P[y^{+},b^{+}]+\{(a,a^{+}),(b^{+},b),(x,x^{+}),(y,y^{+})\}-(a,b)$ is a H-path of $BH_{n}-F$ passing through $L$.

{\it Case 2.} $|E(L_{0})\cup F_{0}|=2n-3$.

By Lemma \ref{le-1}, there is a H-cycle $C_{0}$ passing through $L_{0}$ in $B^{0}-F_{0}$. In this case, $|E(L_{j})\cup F_{j}|\leq 1$ for any $j\in N_{4}\setminus\{0\}$.

{\it Case 2.1.}  $i=0$.

{\it Case 2.1.1.}  $u$ is adjacent to $v$ on $C_0$.

In this case, $P[u,v]=C_0-(u,v)$ is a H-path passing through $L_{0}$ of $B^0-F_0$. Similarly to Case 1, it is easy to construct a H-path of $BH_n-F$ passing through $L$.

{\it Case 2.1.2.}  $u$ is not adjacent to $v$ on $C_0$.

Since $\{u,v\}$ is compatible to $L$, there are two edges $(u,a),(v,b)\in E(C_{0})\setminus E(L_{0})$. Since $|E(L_{j})|\leq 1$ for $j\in N_{4}\setminus\{0\}$, $a^+$ or $a^-$ (resp. $b^+$ or $b^-$) is not incident with an edge of $L_3$ (resp. $L_1$). Without loss of generality, assume $a^+$ (resp. $b^+$) is not incident with an edge of $L_3$ (resp. $L_1$). By Lemma \ref{le-2}, there is an $x\in V_{1}\cap X$ such that $x$ (resp. $x^+$) is not incident with an edge of $L_1$ (resp. $L_2$). Again by Lemma \ref{le-2}, there is a $y\in V_{3}\cap Y$ such that $y$ (resp. $y^+$) is not incident with an edge of $L_3$ (resp. $L_2$).

If each of the two paths between $u$ and $v$ on $C_0$ contains exactly one of $\{a,b\}$,
combining these with the fact that $|E(L_{j})\cup F_{j}|\leq 1$ for $j\in N_{4}\setminus \{0\}$, Theorem \ref{th-yang2019} implies that there are H-paths $P[b^{+},x]$ passing through $L_{1}$ in $B^{1}-F_{1}$, $P[x^{+},y^{+}]$ passing through $L_{2}$ in $B^{2}-F_{2}$, and $P[a^{+},y]$ passing through $L_{3}$ in $B^{3}-F_{3}$. Thus, $C_{0}\cup P[b^{+},x]\cup P[x^{+},y^{+}]\cup P[a^{+},y]+\{(a,a^{+}),(b,b^{+}),(x,x^{+}),(y,y^{+})\}-\{(u,a),(v,b)\}$ is a desired H-path of $BH_n-F$.

In the following, we consider the case that there is a path between $u$ and $v$ on $C_0$ containing both $a$ and $b$. Denote by $P[u,v]$ the other path between $u$ and $v$ on $C_0$. Since $\{u,v\}$ is compatible to $L$, there is an edge, say $(s,t)$, on $E(P[u,v])\setminus E(L_{0})$ for some $s\in X$ and $t\in Y$.

Suppose first that $E(L_{2})\cup F_2=\emptyset$. Note that $|E(L_{k})\cup F_k|\leq 1$ for $k\in \{1,3\}$.
By Theorem \ref{th-yang2019}, $B^{1}-F_{1}$ has a H-path $P[s^{+},x]$ passing through $L_{1}$, $B^{3}-F_{3}$ has a H-path $P[t^{+},y]$ passing through $L_{3}$. Let $c$ be the neighbor of $b^+$ on the segment of $P[s^+,x]$ between $s^+$ and $b^+$, and $d$ be the neighbor of $a^+$ on the segment of $P[t^+,y]$ between $a^+$ and $y$. By Theorem \ref{th-cheng2014}, $B^{2}$ has two vertex-disjoint paths $P[c^{+},d^+]$ and $P[x^{+},y^+]$ such that each vertex of $B^2$ lies on one of the two paths. Thus, $C_{0}\cup P[s^{+},x]\cup P[c^+,d^+]\cup P[x^{+},y^{+}]\cup P[t^{+},y]+\{(a,a^{+}),(b,b^{+}),(c,c^+),(d,d^+),(s,s^+),(t,t^+),(x,x^{+}),(y,y^{+})\}-\{(u,a),(v,b),$ $(s,t),(b^+,c),(a^+,d)\}$ is a desired H-path of $BH_n-F$.

Suppose now that $E(L_{2})\cup F_2\neq \emptyset$. Then $|E(L_{2})\cup F_2|=1$ and $E(L_{k})\cup F_k= \emptyset$ for $k\in \{1,3\}$. Let $c\in V(B^1)\cap X$ such that $c\neq x$. By Theorem \ref{th-yang2019}, $B^{2}-F_{2}$ has a H-path $P[c^{+},y^+]$ passing through $L_{2}$. Let $d$ be the neighbor of $x^+$ on the segment of $P[c^+,y^+]$ between $c^+$ and $x^+$. Note that at least one of the two neighbors of $d$ in $B^3$, say $d^+$, is not $y$. By Theorem \ref{th-cheng2014}, $B^{1}$ has two vertex-disjoint paths $P[b^{+},c]$ and $P[s^{+},x]$ such that each vertex of $B^1$ lies on one of the two paths, and $B^{3}$ has two vertex-disjoint paths $P[a^{+},y]$ and $P[t^{+},d^+]$ such that each vertex of $B^3$ lies on one of the two paths. Thus, $C_{0}\cup P[s^{+},x]\cup P[b^+,c]\cup P[c^+,y^+]\cup P[a^{+},y]\cup P[t^{+},d^+]+\{(a,a^{+}),(b,b^{+}),(c,c^+),(d,d^+),(s,s^+),$ $(t,t^+),(x,x^{+}),(y,y^{+})\}-\{(u,a),(v,b),(s,t),(x^+,d)\}$ is a desired H-path of $BH_n-F$.

{\it Case 2.2.}  $i\in \{1,3\}$.

By symmetry, it suffices to consider that $i=1$. Note that $|E(L_{j})\cup F_{j}|\leq 1$ for $j\in N_{4}\setminus \{0\}$. There is an edge $(a,b)\in E(C_0)\setminus E(L_{0})$ for some $a\in X$ and $b\in Y$ such that $a^+$ or $a^-$ (resp. $b^+$ or $b^-$), say $a^+$ (resp. $b^+$), is not incident with an edge of $L_{1}$ (resp. $L_3$). Theorem \ref{th-yang2019} implies that $B^{1}-F_{1}$ has a H-path $P[u,v]$ passing through $L_{1}$. Let $(a^{+},x)\in E(P[u,v])$. Then $(a^{+},x)\notin E(L_1)$. There is a neighbor of $x$ in $B^2$, say $x^+$, incident with none of $E(L_{2})$. Let $y\in V_{2}\cap X$. By Theorem \ref{th-yang2019}, there are H-paths $P[x^{+},y]$ passing through $L_{2}$ in $B^{2}-F_{2}$ and $P[y^{+},b^{+}]$ passing through $L_{3}$ in $B^{3}-F_{3}$. Thus, $C_{0}\cup P[u,v]\cup P[x^{+},y]\cup P[y^{+},b^{+}]+\{(a,a^{+}),(b,b^{+}),(x,x^{+}),(y,y^{+})\}-\{(a,b),(a^{+},x)\}$ is a H-path of $BH_{n}-F$ passing through $L$.

{\it Case 2.3.}  $i=2$.

Combining these with $|E(L_{j})\cup F_{j}|\leq 1$ for $j\in N_{4}\setminus \{0\}$. By Theorem \ref{th-yang2019}, $B^{2}-F_{2}$ contains a H-path $P[u,v]$ passing through $L_{2}$. By Lemma \ref{le-3}, there is an edge $(a,b)\in P[u,v]\setminus E(L_{2})$ for some $a\in X$ and $b\in Y$ such that $a^+$ or $a^-$ (resp. $b^+$ or $b^-$), say $a^+$ (resp. $b^+$), is incident with none of $E(L_{3})$ (resp. $E(L_1)$). By Lemma \ref{le-2}, there is a vertex $c\in V_{1}\cap Y$ such that $c^{+}$ is not incident with an edge of $L_{0}$. Let $(c^{+},d)\in E(C_{0})$. Then $(c^{+},d)\notin E(L_{0})$. Theorem \ref{th-yang2019} implies that $B^{1}-F_{1}$ has a H-path $P[b^{+},c]$ passing through $L_{1}$ and $B^{3}-F_{3}$ has a H-path $P[a^{+},d^{+}]$ passing through $L_{3}$. Hence, $C_{0}\cup P[b^{+},c]\cup P[u,v]\cup P[a^{+},d^{+}]+\{(a,a^{+}),(b,b^{+}),(c,c^{+}),(d,d^{+})\}-\{(a,b),(c^{+},d)\}$ is a H-path of $BH_{n}-F$ passing through $L$.



{\it Case 3.} $|E(L_{0})\cup F_{0}|=2n-2$.

In this case, $E(L_0)=E(L)\neq \emptyset$, $F_0=F\neq \emptyset$ and $E(L_{j})=F_{j}=\emptyset$ for $j\in N_{4}\setminus \{0\}$.

{\it Case 3.1.}  $i=0$.

Since $\{u,v\}$ is compatible to $L_0$ and $E(L_0)\neq \emptyset$, there is a path in $L_0$ such that at least one of the two end vertices, say $x$, is not in $\{u,v\}$. Without loss of generality, assume that $x\in X$. Let $(x,y)\in E(L_{0})$ and $f\in F_{0}$. By the induction hypothesis, $B^{0}-F_{0}\setminus \{f\}$ has a H-path $P[u,v]$ passing through $L_{0}-(x,y)$. Let $c, z\in V_{1}\cap X$ and $d, w\in V_{2}\cap X$ be pair-wise distinct.

{\it Case 3.1.1.}  $(x,y)\in E(P[u,v])$.

If $f\in E(P[u,v])$, let $(a,b)=f$; otherwise, let $(a,b)$ be an arbitrary edge in $P[u,v]\setminus E(L_{0})$ for some $a\in X$ and $b\in Y$. By Theorem \ref{th-xu2007}, $B^{1}$ has a H-path $P[a^{+},c]$, $B^{2}$ has a H-path $P[c^{+},d]$, $B^{3}$ has a H-path $P[b^{+},d^{+}]$. Thus, $P[u,v]\cup P[a^{+},c]\cup P[c^{+},d]\cup P[b^{+},d^{+}]+\{(a,a^{+}),(b,b^{+}),(c,c^{+}),(d,d^{+})\}-(a,b)$ is a desired H-path of $BH_n-F$.

{\it Case 3.1.2}  $(x,y)\notin E(P[u,v])$.

No matter $y$ is $v$ or not, there is a neighbor $s$ of $y$ on $P[u,v]$ such that $(y,s)\notin E(L_{0})$. Let $(x,t)\in E(P[u,v])$ such that exactly one of $\{s,t\}$ lies on the segment of $P[u,v]$ between $x$ and $y$.

Suppose first that $f\notin E(P[u,v])$ or $f\in \{(x,t),(y,s)\}$. By Theorem \ref{th-xu2007}, $B^{1}$ has a H-path $P[s^{+},z]$, $B^{2}$ has a H-path $P[z^{+},w]$, $B^{3}$ has a H-path $P[t^{+},w^{+}]$. Thus, $P[u,v]\cup P[s^{+},z]\cup P[z^{+},w]\cup P[t^{+},w^{+}]+\{(x,y),(s,s^{+}),(t,t^{+}),$ $(w,w^{+}),(z,z^{+})\}-\{(x,t),(y,s)\}$ is a desired H-path of $BH_n-F$.

Suppose now that $f\in E(P[u,v])$ and $f\notin \{(x,t),(y,s)\}$. Then let $(a,b)=f$ for some $a\in X$ and $b\in Y$. Let $g=b^-$ if $b=t$ and $g=b^+$ otherwise. Let $h=a^-$ if $a=s$ and $h=a^+$ otherwise. By Theorem \ref{th-cheng2014}, $B^{1}$ has two vertex-disjoint paths $P[s^{+},z]$ and $P[h,c]$ such that each vertex of $B^1$ lies on one of the two paths, $B^{2}$ has two vertex-disjoint paths $P[z^{+},w]$ and $P[c^{+},d]$ such that each vertex of $B^2$ lies on one of the two paths, and $B^{3}$ has two vertex-disjoint paths $P[w^{+},t^{+}]$ and $P[d^{+},g]$ such that each vertex of $B^3$ lies on one of the two paths. Thus, $P[u,v]\cup P[s^{+},z]\cup P[h,c]\cup P[z^{+},w]\cup P[c^{+},d]\cup P[w^{+},t^{+}]\cup P[d^{+},g]+\{(x,y),(a,h),(b,g),(c,c^{+}),(d,d^{+}),(s,s^{+}),(t,t^{+}),(z,z^{+}),$ $(w,w^{+})\}-\{(x,t),(y,s),(a,b)\}$ is a desired H-path of $BH_n-F$.

{\it Case 3.2.}  $i\in \{1,3\}$.

By symmetry, it suffices to consider that $i=1$. By Proposition \ref{pr-1}, $B^{0}-F_{0}$ has a H-path $P[a,b]$ passing through $L_{0}$ for some $a\in X$ and $b\in Y$. By Theorem \ref{th-xu2007}, $B^{1}$ has a H-path $P[u,v]$. Let $(a^{+},c)\in E(P[u,v])$ and $d\in V_{2}\cap X$. Then there are H-paths $P[c^{+},d]$ in $B^{2}$ and $P[d^{+},b^{+}]$ in $B^{3}$. Thus, $P[a,b]\cup P[u,v]\cup P[c^{+},d]\cup P[d^{+},b^{+}]+\{(a,a^{+}),(b,b^{+}),(c,c^{+}),(d,d^{+})\}-(a^{+},c)$ is a H-path of $BH_{n}-F$ passing through $L$.

{\it Case 3.3.}  $i=2$.

By Proposition \ref{pr-1}, $B^{0}-F_{0}$ has a H-path $P[a,b]$ passing through $L_{0}$ for some $a\in X$ and $b\in Y$. Let $c\in V_{1}\cap X$. Theorem \ref{th-xu2007} implies that $B^{1}$ and $B^{2}$ have H-paths $P[a^{+},c]$ and $P[u,v]$, respectively. Let $(c^{+},d)\in E(P[u,v])$. Then $B^{3}$ has a H-path $P[d^{+},b^{+}]$. Thus, $P[a,b]\cup P[a^{+},c]\cup P[u,v]\cup P[d^{+},b^{+}]+\{(a,a^{+}),(b,b^{+}),(c,c^{+}),(d,d^{+})\}-(c^{+},d)$ is a H-path of $BH_{n}-F$ passing through $L$.
\end{proof}



\begin{lemma}\label{main-2}
If $|E(L_{0})\cup F_{0}|\leq 2n-4$, $u\in V_{i}$, $v\in V_{j}$ for $i,j\in N_{4}$, and $i\neq j$, then $BH_{n}-F$ has a H-path $P[u,v]$ passing through $L$.
\end{lemma}

\begin{proof}
In this case, $|E(L_{k})\cup F_{k}|\leq 2n-4$ for $k\in N_{4}$. By symmetry, it suffices to consider the following two cases.

{\it Case 1.}   $i=0$.

By Lemma \ref{le-2}, there is an $x\in V_{0}\cap Y$ such that $x$ and $x^{\pm}$ are not incident with an edge of $L_0$ and $E(L_{3})\cup F_3$ respectively. By the induction hypothesis, $B^{0}-F_{0}$ has a H-path $P[u,x]$ passing through $L_{0}$.


{\it Case 1.1.}  $j=1$.

By Lemma \ref{le-2}, there is a $y\in V_{1}\cap X$ such that $y$ (resp. $y^{+})$ is not incident with an edge of $L_1$ (resp. $L_{2}$), and a $z\in V_2\cap X$ such that $z$ (resp. $z^{+}$) is not incident with an edge of $L_2$ (resp. $L_3$). By the induction hypothesis, there are H-paths $P[v,y]$ passing through $L_{1}$ in $B^{1}-F_{1}$, $P[y^{+},z]$ passing through $L_{2}$ in $B^{2}-F_{2}$ and $P[z^{+},x^{+}]$ passing through $L_{3}$ in $B^{3}-F_{3}$. Thus, $P[u,x]\cup P[v,y]\cup P[y^{+},z]\cup P[z^{+},x^{+}]+\{(x,x^{+}),(y,y^{+}),(z,z^{+})\}$ is a H-path of $BH_{n}-F$ passing through $L$.

{\it Case 1.2.}  $j=2$.


{\it Case 1.2.1.}  $|E(L_{2})\cup F_{2}|\geq 2n-5$.

In this scenario, $|E(L_{1})\cup F_{1}|\leq \min\{|E(L)\cup F|-\sum_{k\in N_{4}\setminus \{1\}}|E(L_{k})\cup F_{k}|,\ |E(L_{0})\cup F_{0}|\}\leq 1$. By Lemma \ref{le-2}, there is a $z\in V_{3}\cap Y$ such that $z$ and $z^{\pm}$ are incident with none of $E(L_3)$ and $E(L_{2})$, respectively, and an $s\in V_3\cap X$ such that $s$ and $s^{\pm}$ are incident with none of $E(L_3)$ and $E(L_0)$, respectively. By the induction hypothesis, $B^{3}-F_{3}$ has a H-path $P[s,z]$ passing through $L_{3}$. By Lemma \ref{le-3}, there is an edge $(a,b)\in E(P[s,z])\setminus E(L_3)$ for some $a\in X$ and $b\in Y$ such that $a^+$ or $a^-$ (resp. $b^+$ or $b^-$), say $a^+$ (resp. $b^+$), is not incident with an edge of $L_{0}$ (resp. $L_2$) and $\{a,b\}\cap \{s,z\}=\emptyset$. By the induction hypothesis, $B^0-F_0$ has a H-path $P[u,s^+]$ passing through $L_0$, $B^2-F_2$ has a H-path $P[z^+,v]$ passing through $L_2$. Let $c$ be the neighbor of $a^+$ on the segment of $P[u,s^+]$ between $s^+$ and $a^+$ and let $d$ be the neighbor of $b^+$ on the segment of $P[z^+,v]$ between $z^+$ and $b^+$. Then $(a^+,c)\notin E(L_0)$ and $(b^+,d)\notin E(L_2)$. Recalling that $|E(L_{1})\cup F_{1}|\leq 1$, there is a neighbor of $c$ in $B^1$, say $c^+$, incident with none of $E(L_{1})$. By Theorem \ref{th-yang2019}, $B^{1}-F_{1}$ has a H-path $P[c^{+},d^+]$ passing through $L_{1}$. Thus, $P[u,s^+]\cup P[c^{+},d^{+}]\cup P[z^{+},v]\cup P[s,z]+\{(a,a^{+}),(b,b^{+}),(c,c^{+}),(d,d^{+}),(s,s^{+}),(z,z^{+})\}-\{(a,b),(a^{+},c),(b^{+},d)\}$ is a desired H-path of $BH_n-F$.

{\it Case 1.2.2.}  $|E(L_{2})\cup F_{2}|\leq 2n-6$.

By Lemma \ref{le-3}, there is an edge $(a,b)\in E(P[u,x])\setminus E(L_{0})$ for some $a\in X$ and $b\in Y$ such that $a^+$ or $a^-$ (resp. $b^+$ or $b^-$), say $a^+$ (resp. $b^+$), is incident with none edges of $E(L_{1})$ (resp. $E(L_3)\cup F_3$) and $\{a,b\}\cap \{u,x\}=\emptyset$.

Suppose first that $|E(L_{k})\cup F_{k}|\leq 2n-6$ for any $k\in \{1,3\}$. Lemma \ref{le-9} implies that there are two neighbors $d$ and $t$ of $b^{+}$ in $B^{3}$ such that $d^+$ or $d^-$, and $t^+$ or $t^-$, say $d^+$ and $t^+$, are incident with none of $E(L_2)\cup F_2$ and $L_{3}+\{(d,b^{+}),(b^{+},t)\}$ is a linear forest. By Lemma \ref{le-4}, there is a $z\in V_{3}\cap Y$ such that $z$ and $z^{\pm}$ are not incident with an edge of $L_3$ and $E(L_{2})\cup F_2$, respectively, and $z^{\pm}$ is not adjacent to $v$. Note that $\{x^+,z\}$ is compatible to $L_{3}+\{(d,b^{+}),(b^{+},t)\}$, and $|E(L_{3}+\{(d,b^{+}),(b^{+},t)\})\cup F_{3}|\leq 2n-4$. By the induction hypothesis, $B^{3}-F_{3}$ has a H-path $P[x^{+},z]$ passing through $L_{3}+\{(d,b^{+}),(b^{+},t)\}$. Exactly one of $d$ and $t$, say $d$, lies on the segment of $P[x^+,z]$ between $x^+$ and $b^+$. Recall that $d^+$ or $d^-$, say $d^+$, is incident with none of $E(L_2)\cup F_2$. By Lemma \ref{le-9}, $z^{+}$ has two neighbors $c$ and $g$ in $B^{2}$ such that $c^+$ or $c^-$ (resp. $g^+$ or $g^-$) is incident with none of $E(L_1)$, and $L_2+\{(z^+,c),(z^+,g)\}$ is a linear forest. Note that $\{d^+,v\}$ is compatible to $L_2+\{(z^+,c),(z^+,g)\}$, and $|E(L_{2}+\{(z^+,c),(z^{+},g)\})\cup F_{2}|\leq 2n-4$. By the induction hypothesis, $B^{2}-F_{2}$ has a H-path $P[d^{+},v]$ passing through $L_{2}+\{(z^{+},c),(z^{+},g)\}$. Exactly one of $c$ and $g$ lies on the segment of $P[z^+,v]$ between $z^+$ and $d^+$. $c$ lies on the segment of $P[z^+,v]$ between $z^+$ and $d^+$ if $a$ lies on the segment of $P[u,x]$ between $u$ and $b$, and $g$ lies on the segment of $P[z^+,v]$ between $z^+$ and $d^+$ otherwise. Recall that $c^+$ or $c^-$, say $c^+$, is incident with none of $E(L_1)$. By the induction hypothesis, $B^{1}-F_{1}$ has a H-path $P[a^{+},c^{+}]$ passing through $L_{1}$. Thus, $P[u,x]\cup P[a^{+},c^{+}]\cup P[z^{+},v]\cup P[x^{+},z]+\{(a,a^{+}),(b,b^{+}),(c,c^{+}),(d,d^{+}),(x,x^{+}),(z,z^{+})\}-\{(a,b),(d^{+},c),(b^{+},d)\}$ is a desired H-path of $BH_n-F$.

Suppose now that $|E(L_{k})\cup F_{k}|\geq 2n-5$ for some $k\in \{1,3\}$. If $n=3$, then $|E(L_{2})\cup F_{2}|\leq 2n-6\leq 0$. If $n\geq 4$, then $|E(L_{2})\cup F_{2}|\leq |E(L)\cup F|-|E(L_{0})\cup F_{0}|-|E(L_{k})\cup F_{k}|\leq (2n-2)-2(2n-5)\leq 0$. Therefore, $E(L_{2})\cup F_{2}=\emptyset$ for $n\geq 3$. By Lemma \ref{le-2}, there is a $c\in V_{1}\cap X$ such that $c$ is not incident with an edge of $L_1$, and a $z\in V_3\cap Y$ such that $z$ is not incident with $L_3$. By the induction hypothesis, $B^{1}-F_{1}$ has a H-path $P[a^{+},c]$ passing through $L_{1}$, $B^{3}-F_{3}$ has a H-path $P[x^+,z]$ passing through $L_{3}$. Let $d$ be the neighbor of $b^+$ on the segment of $P[x^+,z]$ between $x^+$ and $b^+$. By Theorem \ref{th-cheng2014}, there exist two vertex-disjoint paths $P[z^{+},c^+]$ and $P[d^+,v]$ in $B^2$ such that each vertex of $B^2$ lies on one of the two paths. Thus, $P[u,a]\cup P[a^+,c]\cup P[z^+,c^+]\cup P[d^+,v]\cup P[x^+,z]+\{(a,a^{+}),(b,b^{+}),(c,c^{+}),(d,d^{+}),(x,x^{+}),(z,z^{+})\}-\{(a,b),(b^+,d)\}$ is a desired H-path of $BH_n-F$.

{\it Case 1.3.}  $j=3$.

Lemma \ref{le-3} implies that there is an edge $(a,b)\in E(P[u,x])\setminus E(L_{0})$ for some $a\in X$ and $b\in Y$ such that $a^+$ or $a^-$ (resp. $b^+$ or $b^-$), say $a^+$ (resp. $b^+$), is incident with none of $E(L_{1})$ (resp. $E(L_3)$) and $\{a,b\}\cap \{u,x\}=\emptyset$.

{\it Case 1.3.1.}  $|E(L_{3})\cup F_{3}|\geq 2n-5$.

In this case, $|E(L_{k})\cup F_{k}|\leq 1$ for $k\in \{1,2\}$. By the induction hypothesis, $B^{3}-F_{3}$ has a H-path $P[b^{+},v]$ passing through $L_{3}$. Let $d$ be the neighbor of $x^+$ on the segment of $P[b^+,v]$ between $b^+$ and $x^+$. There is a neighbor of $d$ in $B^2$, say $d^+$, incident with none of $E(L_{2})$. Let $c\in V_{2}\cap Y$. By Theorem \ref{th-yang2019}, $B^{2}-F_{2}$ has a H-path $P[d^{+},c]$ passing through $L_{2}$, $B^{1}-F_{1}$ has a H-path $P[a^{+},c^{+}]$ passing through $L_{1}$. Thus, $P[u,x]\cup P[a^{+},c^{+}]\cup P[d^{+},c]\cup P[b^{+},v]+\{(a,a^{+}),(b,b^{+}),(c,c^{+}),(d,d^{+}),(x,x^{+})\}-\{(a,b),(x^{+},d)\}$ is a desired H-path of $BH_n-F$.

{\it Case 1.3.2.}  $|E(L_{3})\cup F_{3}|\leq 2n-6$.

{Suppose first that $|E(L_{k})\cup F_{k}|\leq 2n-6$ for any $k\in \{1,2\}$.} By Lemma \ref{le-4}, there is a $y\in V_0\cap Y$ such that $y$ and $y^{\pm}$ are incident with none of $E(L_0)$ and $E(L_3)$, respectively, and $v$ is not adjacent to $y^{\pm}$. By the induction hypothesis, $B^0-F_0$ has a H-path $P[u,y]$ passing through $L_0$. By Lemma \ref{le-3}, there is an edge $(s,t)\in E(P[u,y])\setminus E(L_0)$ for some $s\in X$ and $t\in Y$ such that $s^+$ or $s^-$ (resp. $t^+$ or $t^-$), say $s^+$ (resp. $t^+$), is incident with none of $E(L_{1})$ (resp. $E(L_3)$) and $\{s,t\}\cap \{u,y\}=\emptyset$. Lemma \ref{le-9} implies that there are two neighbors $d$ and $h$ of $y^{+}$ in $B^{3}$ such that $d^+$ or $d^-$, and $h^+$ or $h^-$ are incident with none of $E(L_2)$ and $L_{3}+\{(y^{+},h),(y^{+},h)\}$ is a linear forest. Note that $\{t^+,v\}$ is compatible to $L_{3}+\{(y^{+},h),(y^{+},h)\}$, and $|E(L_{3}+\{(y^{+},h),(y^{+},h)\})\cup F_3|\leq 2n-4$. By the induction hypothesis, $B^{3}-F_{3}$ has a H-path $P[t^{+},v]$ passing through $L_{3}+\{(y^{+},h),(y^{+},h)\}$. Exactly one of $d$ and $h$, say $d$, lies on the segment of $P[t^+,v]$ between $t^+$ and $y^+$. Recall that $d^+$ or $d^-$, say $d^+$, is incident with none of $E(L_2)$. By Lemma \ref{le-2}, there is a $c\in V_2\cap Y$ such that $c$ (resp. $c^+$) is not incident with an edge of $L_2$ (resp. $L_1$). By the induction hypothesis, $B^{2}-F_{2}$ has a H-path $P[d^{+},c]$ passing through $L_{2}$, $B^{1}-F_{1}$ has a H-path $P[s^{+},c^{+}]$ passing through $L_{1}$. Thus, $P[u,y]\cup P[s^{+},c^{+}]\cup P[d^{+},c]\cup P[t^{+},v]+\{(c,c^{+}),(d,d^{+}),(s,s^{+}),(t,t^{+}),(y,y^{+})\}-\{(s,t),(y^{+},d)\}$ is a desired H-path of $BH_n-F$.

Suppose now that $|E(L_{k})\cup F_{k}|\geq 2n-5$ for some $k\in \{1,2\}$. In this case, $E(L_{3})\cup F_{3}=\emptyset$ for $n\geq 3$. By Lemma \ref{le-2}, there is a $c\in V_1\cap X$ such that $c$ (resp. $c^+$) is not incident with an edge of $L_1$ (resp. $L_2$) and a $d\in V_2\cap X$ such that $d$ is not incident with $L_2$. There is a neighbor of $d$ in $B^3$, say $d^+$, being not $v$. By the induction hypothesis, $B^{1}-F_{1}$ has a H-path $P[a^{+},c]$ passing through $L_{1}$, $B^{2}-F_{2}$ has a H-path $P[c^+,d]$ passing through $L_{2}$. By Theorem \ref{th-cheng2014}, there exist two vertex-disjoint paths $P[b^{+},d^+]$ and $P[x^+,v]$ in $B^3$ such that each vertex of $B^3$ lies on one of the two paths. Thus, $P[u,x]\cup P[a^+,c]\cup P[c^+,d]\cup P[x^+,v]\cup P[b^+,d^+]+\{(a,a^{+}),(b,b^{+}),(c,c^{+}),(d,d^{+}),(x,x^{+})\}-(a,b)$
is a desired H-path of $BH_n-F$.

{\it Case 2.} $i\neq 0$.

Without loss of generality, assume that $j>i$. By Lemma \ref{le-2}, there are $x\in X\cap V_0$ and $y\in Y\cap V_0$ such that $x$ and $y$ are incident with none of $E(L_0)$ and $x^{\pm}$ (resp. $y^{\pm}$) are incident with none of $E(L_1)\cup F_1$ (resp. $E(L_3)\cup F_3$). Since $|E(L_{0})\cup F_{0}|\leq 2n-4$, by the induction hypothesis, $B^{0}-F_{0}$ has a H-path $P[x,y]$ passing through $L_{0}$.

{\it Case 2.1.} $i=1,j=2$.

By Lemma \ref{le-2}, there is a $z\in V_2\cap X$ such that $z$ (resp. $z^+$) is not incident with an edge of $L_{2}$ (resp. $L_3$). Since $|E(L_{k})\cup F_{k}|\leq 2n-4$, for $k\in N_{4}\setminus \{0\}$, by the induction hypothesis, $B^{1}-F_{1}$ has a H-path $P[u,x^{+}]$ passing through $L_{1}$, $B^{2}-F_{2}$ has a H-path $P[v,z]$ passing through $L_{2}$ and $B^{3}-F_{3}$ has a H-path $P[z^{+},y^{+}]$ passing through $L_{3}$. Thus, $P[x,y]\cup P[u,x^{+}]\cup P[v,z]\cup P[z^{+},y^{+}]+\{(x,x^{+}),(y,y^{+}),(z,z^{+})\}$ is a H-path of $BH_{n}-F$ passing through $L$.

{\it Case 2.2.} $i=1,j=3$.

By Lemma \ref{le-3}, there is an edge $(a,b)\in E(P[x,y])\setminus E(L_{0})$ for some $a\in X$ and $b\in Y$ such that $a^+$ or $a^-$ (resp. $b^+$ or $b^-$), say $a^+$ (resp. $b^+$), is incident with none of $E(L_{1})$ (resp. $E(L_3)$) and $\{a,b\}\cap \{x,y\}=\emptyset$.

{\it Case 2.2.1}  $|E(L_{3})\cup F_{3}|\geq 2n-5$.

In this case, $|E(L_{k})\cup F_{k}|\leq 1$ for any $k\in \{1,2\}$. By the induction hypothesis, $B^{3}-F_{3}$ has a H-path $P[b^{+},v]$ passing through $L_{3}$. Let $d$ be the neighbor of $y^+$ on the segment of $P[b^+,v]$ between $b^+$ and $y^+$. There is a neighbor of $d$ in $B^2$, say $d^+$, incident with none of $E(L_2)$. Theorem \ref{th-yang2019} implies that $B^{1}-F_{1}$ has a H-path $P[a^{+},u]$ passing through $L_{1}$. Let $c$ be the neighbor of $x^+$ on the segment of $P[a^+,u]$ between $x^+$ and $a^+$. Thus, $B^{2}-F_{2}$ has a H-path $P[d^{+},c^{+}]$ passing through $L_{2}$. Hence, $P[x,y]\cup P[u,a^+]\cup P[d^+,c^+]\cup P[b^+,v]+\{(a,a^{+}),(b,b^{+}),(c,c^{+}),(d,d^{+}),(x,x^{+}),(y,y^{+})\}-\{(a,b),(x^+,c),(y^+,d)\}$ is a H-path passing through $L$ in $BH_n-F$.

{\it Case 2.2.2}  $|E(L_{k})\cup F_{k}|\leq 2n-6$ for any $k\in N_4\setminus \{0\}$.

{By Lemma \ref{le-4},} there are vertices $z\in X\cap V_0$ and $w\in Y\cap V_0$ such that $z$ and $w$ are incident with none of $E(L_0)$, $z^{\pm}$ (resp. $w^{\pm}$) are incident with none of $E(L_1)\cup F_1$ (resp. $E(L_3)\cup F_3$) and $u$ (resp. $v$) is not adjacent to $z^{\pm}$ (resp. $w^{\pm}$). By the induction hypothesis, $B^{0}-F_{0}$ has a H-path $P[z,w]$ passing through $L_{0}$. By Lemma \ref{le-3}, there is an edge $(s,t)\in E(P[z,w])\setminus E(L_{0})$ for some $s\in X$ and $t\in Y$ such that $s^+$ or $s^-$ (resp. $t^+$ or $t^-$), say $s^+$ (resp. $t^+$), is incident with none of $E(L_{1})$ (resp. $E(L_3)$) and $\{s,t\}\cap \{z,w\}=\emptyset$.
Lemma \ref{le-9} implies that there are two neighbors $d$ and $h$ of $w^{+}$ in $B^{3}$ such that $d^+$ or $d^-$, and $h^+$ or $h^-$ are incident with none of $E(L_2)$ and $L_{3}+\{(w^{+},d),(w^{+},h)\}$ is a linear forest. Note that $\{t^+,v\}$ is compatible to $L_{3}+\{(w^{+},d),(w^{+},h)\}$, and $|E(L_{3}+\{(w^{+},d),(w^{+},h)\})\cup F_3|\leq 2n-4$. By the induction hypothesis, $B^{3}-F_{3}$ has a H-path $P[t^{+},v]$ passing through $L_{3}+\{(w^{+},d),(w^{+},h)\}$. Exactly one of $d$ and $h$, say $d$, lies on the segment of $P[t^+,v]$ between $w^+$ and $t^+$. By Lemma \ref{le-9}, $z^{+}$ has two neighbors $c$ and $g$ in $B^{1}$ such that $c^+$ or $c^-$, and $g^+$ or $g^-$ are incident with none of $E(L_2)$. Note that $\{s^+,u\}$ is compatible to $L_1+\{(z^+,c),(z^{+},g)\}$, $|E(L_{1}+\{(z^+,c),(z^{+},g)\})\cup F_{1}|\leq 2n-4$. By the induction hypothesis, $B^{1}-F_{1}$ has a H-path $P[u,s^+]$ passing through $L_{1}+\{(z^{+},c),(z^{+},g)\}$. Exactly one of $c$ and $g$, say $c$, lies on the segment of $P[u,s^+]$ between $z^+$ and $s^+$. Since $d^+$ or $d^-$ (resp. $c^+$ or $c^-$), say $d^+$ (resp. $c^+$), is incident with none of $E(L_2)$, we have that $\{d^+,c^+\}$ is compatible to $L_2$. By the induction hypothesis, $B^{2}-F_{2}$ has a H-path $P[d^{+},c^{+}]$ passing through $L_{2}$. Thus, $P[z,w]\cup P[u,s^+]\cup P[d^+,c^+]\cup P[t^+,v]+\{(c,c^{+}),(d,d^{+}),(s,s^{+}),(t,t^{+}),(w,w^{+}),(z,z^{+})\}-\{(s,t),(z^+,c),(w^+,d)\}$ is a H-path passing through $L$ in $BH_n-F$.

{\it Case 2.2.3}  $|E(L_{3})\cup F_{3}|\leq 2n-6$ and $|E(L_{1})\cup F_{1}|\geq 2n-5$.

In this case, $E(L_{3})\cup F_{3}=\emptyset$ and $|E(L_{2})\cup F_{2}|\leq 1$ for $n\geq 3$. By the induction hypothesis, $B^{1}-F_{1}$ has a H-path $P[u,a^+]$ passing through $L_{1}$. Let $c$ be the neighbor of $x^+$ on the segment of $P[u,a^+]$ between $x^+$ and $a^+$. Note that $|E(L_2)|\leq 1$. There is a $d\in V_2\cap X$ such that $d$ is incident with none of $E(L_2)$. By Theorem \ref{th-yang2019}, $B^2-F_2$ has a H-path $P[c^+,d]$ passing through $L_2$. There is a neighbor of $d$ in $B^3$, say $d^+$, is not $v$. By Theorem \ref{th-cheng2014}, there exist two vertex-disjoint paths $P[b^{+},d^+]$ and $P[y^+,v]$ in $B^3$ such that each vertex of $B^3$ lies on one of the two paths. Thus, $P[x,y]\cup P[u,a^+]\cup P[c^+,d]\cup P[b^+,d^+]\cup P[y^+,v]+\{(a,a^{+}),(b,b^{+}),(c,c^{+}),(d,d^{+}),(x,x^{+}),(y,y^{+})\}-\{(a,b),(x^+,c)\}$ is a H-path passing through $L$ in $BH_n-F$.

{\it Case 2.2.4}  $|E(L_{3})\cup F_{3}|\leq 2n-6$ and $|E(L_{2})\cup F_{2}|\geq 2n-5$.

In this case, $E(L_{3})\cup F_{3}=\emptyset$ and $|E(L_{1})\cup F_{1}|\leq 1$ for $n\geq 3$. By Lemma \ref{le-2}, there is a $z\in V_1\cap Y$ such that $z$ (resp. $z^{\pm}$) is not incident with an edge of $L_1$ (resp. $L_0$). By the induction hypothesis, $B^{1}-F_{1}$ has a H-path $P[u,z]$ passing through $L_{1}$. By Lemma \ref{le-3}, there is an edge $(s,t)\in E(P[u,z])\setminus E(L_{1})$ for some $s\in Y$ and $t\in X$ such that $s^+$ or $s^-$ (resp. $t^+$ or $t^-$), say $s^+$ (resp. $t^+$), is incident with none of $E(L_{0})$ (resp. $E(L_2)$) and $\{s,t\}\cap \{u,z\}=\emptyset$. By Lemma \ref{le-2}, there is a $w\in V_0\cap Y$ such that $w$ is not incident with an edge of $L_0$. By the induction hypothesis, $B^{0}-F_{0}$ has a H-path $P[z^+,w]$ passing through $L_{0}$. Let $c$ be the neighbor of $s^+$ on the segment of $P[z^+,w]$ between $z^+$ and $s^+$. By Lemma \ref{le-2}, there is a $d\in V_2\cap X$ such that $d$ is incident with none of $E(L_2)$. By Theorem \ref{th-yang2019}, $B^2-F_2$ has a H-path $P[t^+,d]$ passing through $L_2$. There is a neighbor of $d$ in $B^3$, say $d^+$, is not $v$. By Theorem \ref{th-cheng2014}, there exist two vertex-disjoint paths $P[w^{+},d^+]$ and $P[c^+,v]$ in $B^3$ such that each vertex of $B^3$ lies on one of the two paths. Thus, $P[z^+,w]\cup P[u,z]\cup P[t^+,d]\cup P[w^+,d^+]\cup P[c^+,v]+\{(s,s^{+}),(t,t^{+}),(c,c^{+}),(d,d^{+}),(z,z^{+}),(w,w^{+})\}-\{(s^+,c),(s,t)\}$
is a H-path passing through $L$ in $BH_n-F$.

{\it Case 2.3.} $i=2,j=3$.

By Lemma \ref{le-2}, there is a $z\in V_{1}\cap X$ such that $z$ (resp. $z^+$) is not incident with an edge of $L_{1}$ (resp. $L_2$). Recall that $|E(L_{k})\cup F_{k}|\leq 2n-4$ for $k\in N_{4}\setminus\{0\}$. By the induction hypothesis, $B^{1}-F_{1}$ has a H-path $P[x^{+},z]$ passing through $L_{1}$, $B^{2}-F_{2}$ has a H-path $P[u,z^{+}]$ passing through $L_{2}$ and $B^{3}-F_{3}$ has a H-path $P[y^{+},v]$ passing through $L_{3}$. Hence, $P[x,y]\cup P[x^{+},z]\cup P[u,z^{+}]\cup P[y^{+},v]+\{(x,x^{+}),(y,y^{+}),(z,z^{+})\}$ is a H-path of $BH_{n}-F$ passing through $L$.
\end{proof}



\begin{lemma}
If $|E(L_{0})\cup F_{0}|=2n-3$, $u\in V_{i}$, $v\in V_{j}$ for $i,j\in N_{4}$ and $i\neq j$, then $BH_{n}-F$ has a H-path $P[u,v]$ passing through $L$.
\end{lemma}

\begin{proof}
In this case, $|E(L_{k})\cup F_{k}|\leq 1$ for $k\in N_{4}\setminus \{0\}$. By Lemma \ref{le-1}, $B^{0}-F_{0}$ has a H-cycle $C_{0}$ passing through $L_{0}$. By symmetry, it suffices to consider the following two cases.

{\it Case 1.}  $i=0$.

Let $(u,x)\in E(C_{0})\setminus E(L_{0})$. There is a neighbor of $x$ in $B^3$, say $x^+$, incident with none of $E(L_3)$. Thus, $P[u,x]=C_0-(u,x)$ is a H-path passing through $L_{0}$ of $B^0-F_0$.

{\it Case 1.1.}  $j=1$.

Let $y\in V_{1}\cap X$ (resp. $z\in V_{2}\cap X$) such that $y$ (resp. $z$) is incident with none of $E(L_1)$ (resp. $E(L_2)$). By Theorem \ref{th-yang2019}, $B^{1}-F_{1}$, $B^{2}-F_{2}$, $B^{3}-F_{3}$ have H-paths $P[v,y]$, $P[y^{+},z]$, $P[z^{+},x^{+}]$ passing through $L_{1}$, $L_{2}$ and $L_{3}$, respectively. Thus, $P[u,x]\cup P[v,y]\cup P[y^{+},z]\cup P[z^{+},x^{+}]+\{(x,x^{+}),(y,y^{+}),(z,z^{+})\}$ is a H-path of $BH_{n}-F$ passing through $L$.

{\it Case 1.2.}  $j=2$ or $j=3$.

By Lemma \ref{le-3}, there is an edge $(a,b)\in E(P[u,x])\setminus E(L_{0})$ for some $a\in X$ and $b\in Y$ such that $a^{\pm}$ (resp. $b^{\pm}$) are incident with none of $E(L_{1})$ (resp. $E(L_3)$), and $\{a,b\}\cap \{u,x\}=\emptyset$.

Suppose first that $j=2$. By Lemma \ref{le-2}, there is a $z\in V_{3}\cap Y$ such that $z^{\pm}$ are incident with none of $E(L_2)$. Recall that $|E(L_{k})\cup F_{k}|\leq 1$ for $k\in N_{4}\setminus \{0\}$. By Theorem \ref{th-yang2019}, $B^{3}-F_{3}$ has a H-path $P[x^{+},z]$ passing through $L_{3}$. Let $d$ be the neighbor of $b^+$ on the segment of $P[x^+,z]$ between $x^+$ and $b^+$. There is a neighbor of $d$, say $d^+$, incident with none of $E(L_2)$. Theorem \ref{th-yang2019} implies that $B^{2}-F_{2}$ has a H-path $P[z^{+},v]$ passing through $L_{2}$. Let $c$ be the neighbor of $d^+$ on the segment of $P[z^+,v]$ between $z^+$ and $d^+$. Again by Theorem \ref{th-yang2019}, $B^{1}-F_{1}$ has a H-path $P[a^{+},c^{+}]$ passing through $L_{1}$. Hence, $P[u,x]\cup P[a^{+},c^{+}]\cup P[z^{+},v]\cup P[x^{+},z]+\{(a,a^{+}),(b,b^{+}),(c,c^{+}),(d,d^{+}),(x,x^{+}),(z,z^{+})\}-\{(a,b),$ $(d^{+},c),(b^{+},d)\}$ is a H-path of $BH_{n}-F$ passing through $L$.

Suppose now that $j=3$. Let $c\in V_{1}\cap X$. There is a neighbor of $c$ in $B^2$, say $c^+$, incident with none of $E(L_2)$. By Theorem \ref{th-yang2019}, $B^{1}-F_{1}$ has a H-path $P[a^{+},c]$ passing through $L_{1}$, $B^{3}-F_{3}$ has a H-path $P[b^{+},v]$ passing through $L_{3}$. Let $d$ be the neighbor of $x^+$ on the segment of $P[b^+,v]$ between $b^+$ and $x^+$. Theorem \ref{th-yang2019} implies that $B^{2}-F_{2}$ has a H-path $P[c^{+},d^{+}]$ passing through $L_{2}$. Hence, $P[u,x]\cup P[a^{+},c]\cup P[c^{+},d^{+}]\cup P[b^{+},v]+\{(a,a^{+}),(b,b^{+}),(c,c^{+}),(d,d^{+}),$ $(x,x^{+})\}-\{(a,b),(x^{+},d)\}$ is a desired H-path of $BH_{n}-F$.

{\it Case 2.} $i\neq 0$.

Without loss of generality, assume that $j>i$.
Note that $|E(L_{k})\cup F_{k}|\leq 1$ for $k\in N_{4}\setminus \{0\}$. By Lemma \ref{le-3}, there is an edge $(a,b)\in E(C_0)\setminus E(L_{0})$ for some $a\in X$ and $b\in Y$ such that $a^{\pm}$ (resp. $b^{\pm}$) are incident with none of $E(L_1)$ (resp. $E(L_3)$).
Thus, $P[a,b]=C_0-(a,b)$ is a H-path passing through $L_{0}$ of $B^0-F_0$.

{\it Case 2.1.}  $i=1,j=2$.

By Lemma \ref{le-2}, there is a $c\in V_{3}\cap Y$ such that $c^+$ is incident with none of $E(L_2)$. Theorem \ref{th-yang2019} implies that $B^{1}-F_{1}$, $B^{2}-F_{2}$, $B^{3}-F_{3}$ have H-paths $P[a^{+},u]$, $P[c^{+},v]$, $P[b^{+},c]$ passing through $L_{1}$, $L_{2}$ and $L_{3}$, respectively. Thus, $P[a,b]\cup P[a^{+},u]\cup P[c^{+},v]\cup P[b^{+},c]+\{(a,a^{+}),(b,b^{+}),(c,c^{+})\}$ is a H-path of $BH_{n}-F$ passing through $L$.


{\it Case 2.2.}  $i=1,j=3$.

By Lemma \ref{le-3}, there is an edge $(x,y)\in E(P[a,b])\setminus E(L_{0})$ for some $x\in X$ and $y\in Y$ such that $x^{\pm}$ (resp. $y^{\pm}$) are incident with none of $E(L_1)$ (resp. $E(L_3)$) and $\{x,y\}\cap \{a,b\}=\emptyset$.
By Theorem \ref{th-yang2019}, $B^{1}-F_{1}$ has a H-path $P[a^{+},u]$ passing through $L_{1}$, $B^{3}-F_{3}$ has a H-path $P[b^{+},v]$ passing through $L_{3}$. Let $z$ be the neighbor of $x^+$ on the segment of $P[a^+,u]$ between $a^+$ and $x^+$ and let $w$ be the neighbor of $y^+$ on the segment of $P[b^+,v]$ between $b^+$ and $y^+$. There is a neighbor of $z$ in $B^2$, say $z^+$, incident with none of $E(L_2)$. Theorem \ref{th-yang2019} implies that $B^{2}-F_{2}$ has a H-path $P[z^{+},w^{+}]$ passing through $L_{2}$. Hence, $P[a,b]\cup P[a^{+},u]\cup P[z^{+},w^{+}]\cup P[b^{+},v]+\{(a,a^{+}),(b,b^{+}),(w,w^+),(x,x^{+}),(y,y^{+}),(z,z^+)\}-\{(x,y),(x^{+},z),(y^{+},w)\}$ is a H-path of $BH_{n}-F$ passing through $L$.

{\it Case 2.3}  $i=2,j=3$.

Let $c\in V_{1}\cap X$. There is a neighbor of $c$ in $B^2$, say $c^+$, incident with none of $E(L_2)$. Theorem \ref{th-yang2019} implies that $B^{1}-F_{1}$, $B^{2}-F_{2}$, $B^{3}-F_{3}$ have H-paths $P[a^{+},c]$, $P[c^{+},u]$, $P[b^{+},v]$ passing through $L_{1}$, $L_{2}$ and $L_{3}$, respectively. Hence, $P[a,b]\cup P[a^{+},c]\cup P[c^{+},u]\cup P[b^{+},v]+\{(a,a^{+}),(b,b^{+}),(c,c^{+})\}$ is a H-path of $BH_{n}-F$ passing through $L$.
\end{proof}



\begin{lemma}
If $|E(L_{0})\cup F_{0}|=2n-2$, $u\in V_{i}$, $v\in V_{j}$ for $i,j\in N_{4}$, and $i\neq j$, then $BH_{n}-F$ has a H-path $P[u,v]$ passing through $L$.
\end{lemma}

\begin{proof} In this case, $E(L_0)=E(L)\neq \emptyset$, $F_0=F\neq \emptyset$ and $E(L_{k})=F_{k}=\emptyset$ for $k\in N_{4}\setminus \{0\}$. Proposition \ref{pr-1} implies that $B^{0}-F_{0}$ has a H-path $P[a,b]$ passing through $L_{0}$. There is a neighbor of $a$ (resp. $b$) in $B^1$ (resp. $B^3$), say $a^+$ (resp. $b^+$), being not $v$ (resp. $u$). By symmetry, it suffices to consider the following two cases.

{\it Case 1.}  $i=0$.

{\it Case 1.1.}  $u\neq a$.

In this case, there is an edge $(u,x)$ on $P[a,b]$ but not in $L_{0}$. Let $t=x^-$ if $x=b$ and let $t=x^+$ otherwise. Then $t\neq b^+$.

Suppose first that $j=1$. Let $c,y\in V_{1}\cap X$, $d,z\in V_{2}\cap X$ be pair-wise distinct. Theorem \ref{th-cheng2014} implies that $B^{1}$ has two vertex-disjoint paths $P[a^{+},c]$ and $P[v,y]$ such that each vertex of $B^1$ lies on one of the two paths and $B^{2}$ has two node-disjoint paths $P[c^{+},d]$ and $P[y^{+},z]$ such that each vertex of $B^2$ lies on one of the two paths and $B^{3}$ has two vertex-disjoint paths $P[d^{+},b^{+}]$ and $P[z^{+},t]$ such that each vertex of $B^3$ lies on one of the two paths. Thus, $P[a,b]\cup P[a^{+},c]\cup P[v,y]\cup P[c^{+},d]\cup P[y^{+},z]\cup P[d^{+},b^{+}]\cup P[z^{+},t]+\{(a,a^{+}),(b,b^{+}),(c,c^{+}),(d,d^{+}),(x,t),(y,y^{+}),(z,z^{+})\}-(u,x)$ is a H-path of $BH_{n}-F$ passing through $L$.

Suppose second that $j=2$. Let $c\in V_{1}\cap X$ such that $c^+\neq v$ and let $d,z\in V_{2}\cap X$ such that $d\neq z$. Theorem \ref{th-xu2007} implies that $B^{1}$ has a H-path $P[a^{+},c]$. By Theorem \ref{th-cheng2014}, $B^{2}$ (resp. $B^3$) has two vertex-disjoint paths $P[c^{+},d]$ and $P[v,z]$ (resp. $P[d^{+},b^{+}]$ and $P[z^{+},t]$) such that each vertex of $B^2$ (resp. $B^3$) lies on one of the two paths. Hence, $P[a,b]\cup P[a^{+},c]\cup P[c^{+},d]\cup P[v,z]\cup P[d^{+},b^{+}]\cup P[z^{+},t]+\{(a,a^{+}),(b,b^{+}),(c,c^{+}),(d,d^{+}),(x,t),(z,z^{+})\}-(u,x)$ is a H-path of $BH_{n}-F$ passing through $L$.

Suppose now that $j=3$. Let $c\in V_{1}\cap X$, $d\in V_{2}\cap X$ such that $d^+\neq v$. Theorem \ref{th-xu2007} implies that $B^{1}$ has a H-path $P[a^{+},c]$, $B^{2}$ has a H-path $P[c^{+},d]$. By Theorem \ref{th-cheng2014}, $B^{3}$ has two vertex-disjoint paths $P[d^{+},b^{+}]$ and $P[t,v]$ such that each vertex of $B^3$ lies on one of the two paths. Hence, $P[a,b]\cup P[a^{+},c]\cup P[c^{+},d]\cup P[d^{+},b^{+}]\cup P[t,v]+\{(a,a^{+}),(b,b^{+}),(c,c^{+}),(d,d^{+}),(x,t)\}-(u,x)$ is a desired H-path of $BH_{n}-F$.

{\it Case 1.2.}  $u=a$.

If $j=1$, let $c\in V_1\cap X$, $d\in V_2\cap X$. Theorem \ref{th-xu2007} implies that $B^{1}$, $B^{2}$, $B^{3}$ have H-paths $P[v,c]$, $P[c^+,d]$, $P[d^+,b^+]$, respectively. Thus, $P[u,b]\cup P[v,c]\cup P[c^+,d]\cup P[d^+,b^+]+\{(b,b^+),(c,c^+),(d,d^+)\}$ is a H-path of $BH_{n}-F$ passing through $L$.

If $j=2$ or $j=3$, since $|E(P[u,b])\setminus E(L_0)|\geq (4^{n-1}-1)-(2n-2)\geq 11$, there is an edge $(x,y)\in E(P[u,b])\setminus E(L_0)$ for some $x\in X$ and $y\in Y$ such that $\{x,y\}\cap \{u,b\}=\emptyset$.

Suppose first that $j=2$. Let $z\in V_1\cap X$ such that $z^+\neq v$ and $c,w\in V_2\cap X$ such that $c\neq w$. Theorem \ref{th-xu2007} implies that $B^{1}$ has a H-path $P[x^+,z]$. By Theorem \ref{th-cheng2014}, $B^{2}$ (resp. $B^3$) has two vertex-disjoint paths $P[z^{+},w]$ and $P[v,c]$ (resp. $P[c^{+},b^{+}]$ and $P[w^{+},y^+]$) such that each vertex of $B^2$ (resp. $B^3$) lies on one of the two paths. Thus, $P[u,b]\cup P[x^{+},z]\cup P[z^{+},w]\cup P[v,c]\cup P[c^{+},b^{+}]\cup P[w^{+},y^+]+\{(b,b^{+}),(c,c^{+}),(w,w^+),(x,x^+),(y,y^+),(z,z^{+})\}-(x,y)$ is a H-path of $BH_{n}-F$ passing through $L$.

Suppose now that $j=3$. Let $z\in V_1$, $w\in V_2$ such that $w^+\neq v$. By Theorem \ref{th-xu2007}, $B^{1}$ and $B^{2}$ have H-paths $P[x^+,z]$ and $P[z^+,w]$. Theorem \ref{th-cheng2014} implies that $B^{3}$ has two vertex-disjoint paths $P[w^{+},y^+]$ and $P[b^+,v]$ such that each vertex of $B^3$ lies on one of the two paths. Thus, $P[u,b]\cup P[x^{+},z]\cup P[z^{+},w]\cup P[b^{+},v]\cup P[w^{+},y^+]+\{(b,b^{+}),(w,w^+),(x,x^+),(y,y^+),(z,z^{+})\}-(x,y)$ is a H-path of $BH_{n}-F$ passing through $L$.





{\it Case 2.} $i\neq 0$.

Without loss of generality, assume that $j>i$.

{\it Case 2.1.}  $i=1,j=2$.


Let $c\in V_{3}\cap Y$. By Theorem \ref{th-xu2007}, $B^{1}$, $B^{2}$, $B^{3}$ have H-paths $P[u,a^{+}]$, $P[v,c^{+}]$ and $P[c,b^{+}]$, respectively. Thus, $P[a,b]\cup P[u,a^{+}]\cup P[v,c^{+}]\cup P[c,b^{+}]+\{(a,a^{+}),(b,b^{+}),(c,c^{+})\}$ is a desired H-path of $BH_{n}-F$.


{\it Case 2.2.}  $i=1,j=3$.

There is an edge $(x,y)\in E(P[a,b])\setminus E(L_{0})$ for some $x\in X$ and $y\in Y$ such that $\{x,y\}\cap \{a,b\}=\emptyset$. Let $z\in V_{1}\cap X\setminus \{u\}$ and $w\in V_{3}\cap Y\setminus \{v\}$. Theorem \ref{th-cheng2014} implies that $B^{1}$ (resp. $B^3$) has two vertex-disjoint paths $P[a^{+},u]$ and $P[x^{+},z]$ (resp. $P[b^{+},v]$ and $P[y^{+},w]$) such that each vertex of $B^1$ (resp. $B^3$) lies on one of the two paths. By Theorem \ref{th-xu2007}, $B^{2}$ has a H-path $P[z^{+},w^+]$. Thus, $P[a,b]\cup P[a^{+},u]\cup P[x^{+},z]\cup P[z^{+},w^+]\cup P[y^{+},w]\cup P[b^{+},v]+\{(a,a^{+}),(b,b^{+}),(w,w^{+}),(x,x^{+}),(y,y^{+}),(z,z^{+})\}-(x,y)$ is a H-path of $BH_{n}-F$ passing through $L$.

{\it Case 2.3.}  $i=2,j=3$.

Let $c\in V_{1}\cap X$. By Theorem \ref{th-xu2007}, $B^{1}$, $B^{2}$, $B^{3}$ have H-paths $P[a^{+},c]$, $P[c^{+},u]$ and $P[b^{+},v]$, respectively. Thus, $P[a,b]\cup P[a^{+},c]\cup P[c^{+},u]\cup P[b^{+},v]+\{(a,a^{+}),(b,b^{+}),(c,c^{+})\}$ is a desired H-path of $BH_{n}-F$.
\end{proof}






\section{$F^c=\emptyset$ and $|L^c|=1$.}
\label{section4}


In this section, let $(x,x^+)$ be the edge of $L^c$ for some $x\in X$ and $x^+\in Y$, and assume $(x,x^+)\in E_{l,l+1}$ for some $l\in N_4$.

\begin{lemma}\label{le-8}

Let $r\in V_j\cap X$ (resp. $r\in V_j\cap Y$) be incident with at most one edge of $L_j$, and let $y\in X$ and $z\in Y$ such that $\{y,z\}$ is compatible to $L_j$. If $|E(L_{0})\cup F_{0}|\leq 2n-5$, then there is a neighbor $s$ of $r$ in $B^j$ such that

\noindent $(i)$. $(r,s)\notin E(L_j)$; and

\noindent $(ii)$. $L_{j}+(r,s)$ is a linear forest; and

\noindent $(\romannumeral3)$. $\{y,z\}$ is compatible to $L_{j}+(r,s)$; and

\noindent $(\romannumeral4)$. $s^{+}$ or $s^{-}$ is not an internal vertex of $L_{j-1}$ (resp. $L_{j+1}$); and

\noindent {$(\romannumeral5)$. furthermore, if $|E(L_{0})\cup F_{0}|\leq 2n-6$ and $y$ (resp. $z$) is incident with none of $E(L_j)$, $s^+$ or $s^-$ is incident with none of $E(L_{j-1})$ (resp. $E(L_{j+1})$).}

\end{lemma}

\begin{proof}
For $n=3$, $|E(L_{k})\cup F_{k}|\leq |E(L_{0})\cup F_{0}|\leq 2n-5=1$, and it is not hard to verify that the lemma holds. It remains consider that $n\geq 4$.
The proofs for the cases that $r\in V_j\cap X$ and $r\in V_j\cap Y$ are analogous. We here only consider the case that $r\in V_j\cap X$.





There are $|N_{B^{j}}(r)|=2n-2$ vertex candidates.
Clearly, the number of such $s$ that fails $(i)$ does not exceed 1. Since there are at most $\lceil(|E(L_{j})|-1/{2}\rceil$ internal vertices in $L_j$, and there is at most one path between $r$ and $s$ in $L_j$, the number of such $s$ that fails $(ii)$ does not exceed $\lceil(|E(L_{j})|-1)/{2}\rceil+1$. There is at most one path $P[y,a]$ (resp. $P[z,b]$) in $L_j$ taking $y$ (resp. $z$) as an end vertex, and there is no path between $y$ and $z$ in $L_j$.
If an $s$ supports $(i)$ and $(ii)$ but fails $(iii)$ then $\{a,b\}=\{r,s\}$, and so the number of such $s$ does not exceed $1$.

Suppose first that $|E(L_{0})\cup F_{0}|\leq 2n-5$. Let $H$ be the set of internal vertices in $L_{j-1}$. Then $|H|\leq \lceil(|E(L_{j-1})|-1)/{2}\rceil$. For two distinct $w,h\in H$, if $w$ is the shadow vertex of $h$, then the two vertices $w^+$ (i.e., $h^-$) and $w^-$ (i.e., $h^+$) fail $(\romannumeral4)$. Therefore, the $|H|$ vertices in $H$ will make at most $|H|$ vertices of $N_{B^{j}}(r)$ fail $(\romannumeral4)$. Note that $F\neq \emptyset$, and $|E(L)|\leq |E(L)\cup F|-|F|\leq (2n-2)-1=2n-3$.
Thus, the total number of vertex candidates that fail the lemma does not exceed $1+(\lceil(|E(L_{j})|-1)/{2}\rceil+1)+1+|H|\leq 3+\lceil(|E(L_{j})|-1)/{2}\rceil+\lceil(|E(L_{j-1})|-1)/{2}\rceil\leq 3+(|E(L_{j})|+|E(L_{j-1})|)/{2}\leq 3+(|E(L)|-|L^c|)/{2}\leq 3+((2n-3)-1)/{2}=n+1$. Since $|N_{B^{j}}(r)|-(n+1)=(2n-2)-(n+1)>0$ for $n\geq 4$, there is an $s\in N_{B^{j}}(r)$ supporting the lemma.

Suppose now that $|E(L_{0})\cup F_{0}|\leq 2n-6$ and $y$ (resp. $z$) is incident with none of $E(L_j)$. Then the number of such $s$ does not exceed $0$. Let $H$ be the set of even vertices which are not singletons in $L_{j-1}$. Then $|H|\leq |E(L_{j-1})|$. For two distinct $w,h\in H$, if $w$ is the shadow vertex of $h$, then the two vertices $w^+$ (i.e., $h^-$) and $w^-$ (i.e., $h^+$) fail $(\romannumeral5)$. Therefore, the $|H|$ vertices in $H$ will make at most $|H|$ vertices of $N_{B^{j}}(r)$ fail $(\romannumeral4)$. Note that $F\neq \emptyset$, and $|E(L)|\leq |E(L)\cup F|-|F|\leq (2n-2)-1=2n-3$.
Thus, the total number of vertex candidates that fail the lemma does not exceed  $1+(\lceil(|E(L_{j})|-1)/{2}\rceil+0)+1+|H|\leq 2+\lceil(|E(L_{j})|-1)/{2}\rceil+|E(L_{j-1})|\leq 2+(|E(L)|-|L^c|)/{2}+(|E(L_{j-1})|/2)\leq 2+((2n-3)-1)/{2}+(2n-6)/2=2n-3$. Since $|N_{B^{j}}(r)|-(2n-3)=(2n-2)-(2n-3)>0$, there is an $s\in N_{B^{j}}(r)$ supporting the lemma.
\end{proof}

\begin{lemma}\label{le-10}

Given $l\in N_4$, suppose $P[x,r]$ is a maximal path of $L_{l}$. Let $(x,y)\in E(P[x,r])$ and let $z\in V_l\cap X-\{x,r\}$ such that $z$ is incident with at most one edge of $L_{l}$. If $|E(L_{l})\cup F_{l}|\leq 2n-6$, then there are two neighbors $s$ and $t$ of $x$ in $N_{B^{l}}(x)\setminus \{y\}$ such that

\noindent $(i)$. $L_{l}+\{(x,s),(x,t)\}-(x,y)$ is a linear forest, and $\{y,z\}$ is compatible to $L_{l}+\{(x,s),(x,t)\}-(x,y)$; and

\noindent $(\romannumeral2)$. $s^{+}$ or $s^{-}$ is incident with none of $E(L_{l-1})\cup F_{l-1}$; and

\noindent $(\romannumeral3)$. $t^{+}$ or $t^{-}$ is not an internal vertex of $L_{l-1}$; and 

\noindent $(\romannumeral4)$. $t$ is not the shadow vertex of $s$.
\end{lemma}

\begin{proof}


A vertex $s\in N_{B^{l}}(x)$ fails the lemma only if

\noindent ($a$). $s$ is incident with an edge of $L_l$; or

\noindent ($b$). $s^{\pm}$ are incident with an edge of $E(L_{l-1})\cup F_{l-1}$.

There are $|N_{B^{l}}(x)|=2n-2$ vertex candidates.
Since there are at most $|E(L_l)|$ vertices incident with an edge of $L_l$, the number of such $s$ that supports $(a)$ does not exceed $|E(L_l)|$. Let $H$ be the set of even vertices which are not singletons in $L_{l-1}\cup F_{l-1}$. Then $|H|\leq |E(L_{l-1}\cup F_{l-1}|$. For two distinct $w,h\in H$, if $w$ is the shadow vertex of $h$, then the two vertices $w^+$ (i.e., $h^-$) and $w^-$ (i.e., $h^+$) support $(b)$. Therefore, the $|H|$ vertices in $H$ will make at most $|H|$ vertices of $N_{B^{l}}(x)$ support $(b)$.
Thus, the total number of such $s\in N_{B^{l}}(x)$ failing the lemma does not exceed $|E(L_l)|+|H|\leq |E(L_l)\cup F_l|+|E(L_{l-1})\cup F_{l-1}|\leq |E(L)\cup F|-1\leq 2n-3$. Since $|N_{B^{l}}(x)|-(2n-3)=(2n-2)-(2n-3)>0$, there is a vertex $s\in N_{B^{l}}(x)$ supporting the lemma.

A vertex $t\in N_{B^{l}}(x)\setminus \{s\}$ fails the lemma only if

\noindent ($c$). $t$ is incident with an edge of $L_{l}$; or

\noindent ($d$). $t^{\pm}$ are internal vertices of $L_{l-1}$; or

\noindent ($e$). $t$ is the shadow of $s$.

Since there are at most $|E(L_l)|$ vertices incident with an edge of $L_l$, the number of such $t$ that supports $(c)$ does not exceed $|E(L_l)|$.
Let $H$ be the set of internal vertices which are not singletons in $L_{l-1}$. Then $|H|\leq \lceil(|E(L_{l-1})|-1)/{2}\rceil$. For two distinct $w,h\in H$, if $w$ is the shadow vertex of $h$, then the two vertices $w^+$ (i.e., $h^-$) and $w^-$ (i.e., $h^+$) support $(d)$. Therefore, the $|H|$ vertices in $H$ will make at most $|H|$ vertices of $N_{B^{l}}(x)\setminus \{s\}$ support $(d)$. Clearly, the number of such $t$ that supports $(e)$ does not exceed 1. Note that $F\neq \emptyset$, and $|E(L)|\leq |E(L)\cup F|-|F|\leq (2n-2)-1=2n-3$.
Thus, the total number of such $t\in N_{B^{l}}(x)\setminus \{s\}$ failing the lemma does not exceed $|E(L_l)|+\lceil(|E(L_{l-1})|-1)/{2}\rceil+1\leq |E(L_l)|+(|E(L_{l-1})|)/{2}+1\leq \frac{(|E(L_l)|+|E(L_{l-1})|)+|E(L_{l})|}{2}+1\leq \frac{(|E(L)\setminus L^c|)+|E(L_{l})|}{2}+1\leq \frac{(2n-3)-1+(2n-6)}{2}+1=2n-4$. Since $|N_{B^{l}}(x)\setminus \{s\}|-(2n-4)=(2n-3)-(2n-4)>0$, then there is a vertex $t\in N_{B^{l}}(x)\setminus \{s\}$ supporting the lemma.
\end{proof}

\begin{lemma}\label{le-11}

Let $y\in V_1\cap Y$ (resp. $y\in V_{3}\cap Y$) such that $y$ is incident with none of $E(L_1)$ (resp. $E(L_{3})$) if $l=1$ (resp. $l=2$). If $|E(L_{0})\cup F_{0}|\in \{2n-5,2n-4\}$, then there is a $z\in N_{B^1}(x)-\{y\}$ (resp. $z\in N_{B^{3}}(x^+)-\{y\}$) such that

\noindent ($i$). $(x,z)\notin E(L_1)$ (resp. $(x^+,z)\notin E(L_{3})$); and

\noindent ($ii$). $L_1+(x,z)$ (resp. $L_{3}+(x^+,z)$) is a linear forest; and

\noindent ($iii$). $z^+$ or $z^-$ is incident with none of $E(L_0)$.

\end{lemma}

\begin{proof}

The proofs for the cases that $j=1$ and $j=2$ are analogous. We here only consider the case that $j=1$.
There are $|N_{B^1}(x)-\{y\}|\leq 2n-3$ candidates of $z$. Note that $x$ is incident with $L^c$. None of candidate of $z$ fails $(i)$ if $E(L_1)=\emptyset$; and at most one, otherwise. Note that $|E(L_1)|\leq 2$. There is no internal vertex of $L_1$ if $|E(L_1)|\leq 1$; and at most one, otherwise. Therefore none of candidate of $z$ fails $(ii)$ if $|E(L_1)|\leq 1$; and at most one, otherwise. Let $H$ be the set of even vertices which are not singletons in $L_0$. For two distinct $s,h\in H$, if $s$ is the shadow vertex of $t$, then the two vertices $s^+\in V_1\cap Y$ (i.e., $h^-$) and $s^-\in V_j\cap Y$ (i.e., $h^+$) may be not as candidates of $z$. Thus, per each vertex in $H$ fails at most one candidate of $z$, and so at most $|H|\leq |E(L_0)|\leq |E(L)|-|L^c|-|E(L_1)|\leq (2n-3)-1-|E(L_1)|=2n-4-|E(L_1)|$ candidates of $z$ fails $(iii)$.
If $|E(L_1)|=2$, then the total number of such $z$ failing the lemma does not exceed $1+1+|H|\leq 2+(2n-4-2)\leq 2n-4<|N_{B^1}(x)-\{y\}|$.
If $|E(L_1)|=1$, then the total number of such $z$ failing the lemma does not exceed $1+0+|H|\leq 1+(2n-4-1)=2n-4<|N_{B^1}(x)-\{y\}|$.
If $E(L_1)=\emptyset$, then the total number of such $z$ failing the lemma does not exceed $0+0+|H|\leq 2n-4<|N_{B^1}(x)-\{y\}|$.
The lemma holds.
\end{proof}

\begin{lemma}\label{le-13}

{Let $y\in V_0\cap Y$ such that $y$ is incident with none of $E(L_0)$ if $l=0$ and $|E(L_{0})\cup F_{0}|=2n-5$, then there is a $z\in N_{B^0}(x)-\{y\}$ such that

\noindent ($i$). $(x,z)\notin E(L_0)$; and

\noindent ($ii$). $L_0+(x,z)$ is a linear forest; and

\noindent ($iii$). $z^+$ or $z^-$ is incident with none of $E(L_0)$.}

\end{lemma}

\begin{proof}

There are $|N_{B^0}(x)-\{y\}|\leq 2n-3$ candidates of $z$. Note that $x$ is incident with $L^c$. The number candidate of $z$ fails $(i)$ at most $1$. Note that $|E(L_3)|\leq 2$. There is no internal vertex of $L_3$ if $|E(L_3)|\leq 1$; and at most one, otherwise. Since there are at most $\lceil(|E(L_{0})|-1/{2}\rceil$ internal vertices in $L_0$, the number candidate of such $z$ that fails $(ii)$ does not exceed $\lceil(|E(L_{0})|-1)/{2}\rceil$.
Let $H$ be the set of even vertices which are not singletons in $L_3$. For two distinct $s,h\in H$, if $s$ is the shadow vertex of $t$, then the two vertices $s^+\in V_0\cap Y$ (i.e., $h^-$) and $s^-\in V_0\cap Y$ (i.e., $h^+$) may be not as candidates of $z$. Thus, per each vertex in $H$ fails at most one candidate of $z$, and so at most $|H|\leq |E(L_3)|$ candidates of $z$ fails $(iii)$. Note that $F\neq \emptyset$.

If $|E(L_3)|=2$, then $n\geq 4$, $|F_0|=|F|\geq 1$, $|E(L_0)|\leq |E(L_0)\cup F_0|-|F_0|\leq (2n-5)-1=2n-6$, the total number of such $z$ failing the lemma does not exceed $1+\lceil(|E(L_{0})|-1)/{2}\rceil+|H|\leq 1+|E(L_{0})|/{2}+2\leq 1+(2n-6)/2+2=n<|N_{B^1}(x)-\{y\}|$ for $n\geq 4$.

If $|E(L_3)|\leq 1$,
then the total number of such $z$ failing the lemma does not exceed $1+\lceil(|E(L_{0})|-1)/{2}\rceil+|H|\leq 1+(2n-6)/{2}+1=n-1<|N_{B^1}(x)-\{y\}|$.
The lemma holds.
\end{proof}

\begin{lemma}\label{le-12}

Suppose $l=0$ and $|E(L_{0})\cup F_{0}|=2n-5$. Let $P[x,r]$ be a maximal path of $L_{0}$ with $r\neq u$ and let $(x,y)\in E(P[x,r])$. Then there are distinct vertices $s,t\in N_{B^{0}}(x)\setminus \{y\}$ such that

\noindent $(i)$. $L_{0}+\{(x,s),(x,t)\}-(x,y)$ is a linear forest; and

\noindent $(\romannumeral2)$. $s^{\pm}$ are incident with none of $E(L_{3})$.

\end{lemma}

\begin{proof}

There are $|N_{B^{0}}(x)\setminus \{y\}|$ vertex candidates of $s$. An $s\in N_{B^{0}}(x)\setminus \{y\}$ fails $(i)$ only if $s$ is an internal vertex of $L_0$. The number of such $s$ that fails $(i)$ does not exceed $\lceil(|E(L_{0})|-1)/{2}\rceil$ because there are at most $\lceil(|E(L_{0})|-1)/{2}\rceil$ internal vertices in $L_0$. An $s\in N_{B^{0}}(x)\setminus \{y\}$ fails $(ii)$ only if $s^{+}$ or $s^-$ is incident with an edge of $E(L_{3})$. Let $H$ be the set of even vertices which are not singletons in $L_3$. For two distinct $g,h\in H$, if $g$ is the shadow vertex of $h$, then the two vertices $g^+\in V_0\cap Y$ (i.e., $h^-$) and $g^-\in V_0\cap Y$ (i.e., $h^+$) may be not as candidates of $s$. Thus, per each vertex in $H$ fails at most two candidates of $s$, and so at most $|H|\leq 2|E(L_3)|$ candidates of $s$ fails $(ii)$. Then the total number of such $s$ failing the lemma does not exceed $\lceil(|E(L_{0})|-1)/{2}\rceil+2|H|\leq |(2n-6)/2+2|E(L_3)|\leq n-3+4=n+1<2n-3$ for $n>4$. For $n=3$, $|E(L_0)|=|E(L_3)|\leq 1$, the total number of such $s$ failing the lemma does not exceed $0+2|H|\leq 2<2n-3$. We now consider that $n=4$. In this scenario, $|E(L)|\leq |E(L)\cup F|-|F|\leq (2n-2)-1=5$ and $\sum_{k=0}^3|E(L_k)|\leq |E(L)|-|L^c|\leq 4$.

Suppose first that $|E(L_3)|=2$. Then $|E(L_0)|\leq 2$, and $L_0$ has no internal vertex or has exactly one internal vertex (i.e. $y$). No matter which case above, the number of such $s$ that fails $(i)$ is $0$. Then the total number of such $s$ failing the lemma does not exceed $0+2|H|\leq 4<2n-3$.

Suppose now that $|E(L_3)|\leq 1$. Then the total number of such $s$ failing the lemma does not exceed $\lceil(|E(L_{0})|-1)/{2}\rceil+2|H|\leq \lceil(|E(L_{0})|-1)/{2}\rceil+2|E(L_3)|\leq (2n-6)/2+2=3<2n-3$.

Note that $|N_{B^{0}}(x)\setminus \{y\}|=2n-3$. There is a vertex $s\in N_{B^{0}}(x)\setminus \{y\}$ supporting the lemma.

There are $|N_{B^{0}}(x)\setminus \{y,s\}|=2n-4$ vertex candidates of $t$. A $t\in N_{B^{0}}(x)\setminus \{y,s\}$ fails $(i)$ only if $s$ is an internal vertex of $L_0$ or $P[t,s]$ is a maximal path of $L_{0}$ or $t^{\pm}$ are incident with an edge of $E(L_3)$. Since $L_0$ has at most $\lceil(|E(L_{0})|-1)/{2}\rceil$ internal vertices and has at most one maximal path which takes $t$ and $s$ as end vertices, the total number of such $t$ fails the lemma does not exceed $\lceil(|E(L_{0})|-1)/{2}\rceil+1\leq |(2n-6)/2+1=n-2<2n-4$.
Therefore there is a vertex $t\in N_{B^{0}}(x)\setminus \{y,s\}$ supporting the lemma.
\end{proof}

\begin{lemma}

If $|E(L_{0})\cup F_{0}|\leq 2n-4$ and $u,v\in V_i$ for $i\in N_{4}$, then $BH_{n}-F$ has a H-path $P[u,v]$ passing through $L$.
\end{lemma}


\begin{proof} In this case, $|E(L_{k})\cup F_{k}|\leq 2n-5$, for each $k\in N_{4}\setminus \{0\}$. In this scenario, the proofs of the cases $l=0$, $l=1$, $l=2$ and $l=3$ are analogous. We here only consider the case $l=0$.

{\it Case 1.}   $i=0$.

By Lemma \ref{le-2}, there are vertices $z\in V_1\cap X$, $w\in V_2\cap X$ such that $z$ (resp. $z^+$) is incident with none of $E(L_1)$ (resp. $E(L_2)$), and $w$ (resp. $w^+$) is incident with none of $E(L_2)$ (resp. $E(L_3)$).

Suppose first that $|E(L_{0})\cup F_{0}|\leq 2n-5$. By Lemma \ref{le-8}, there is a neighbor $y$ of $x$ in $B^{0}$ such that $(x,y)\notin E(L_0)$, $L_0+(x,y)$ is a linear forest, $\{u,v\}$ is compatible to $L_0+(x,y)$, and $y^+$ or $y^-$, say $y^+$, is not an internal vertex of $L_3$. Note that $|E(L_{0}+\{(x,y)\})\cup F_{0}|\leq 2n-4$. By the induction hypothesis, $B^{0}-F_{0}$ has a H-path $P[u,v]$ passing through $L_0+(x,y)$.

Suppose now that $|E(L_{0})\cup F_{0}|=2n-4$. In this case, $|E(L_{k})\cup F_{k}|\leq 1$ for any $k\in N_{4}\setminus \{0\}$. By the induction hypothesis, $B^{0}-F_{0}$ has a H-path $P[u,v]$ passing through $L_{0}$. Let $(x,y)\in E(P[u,v])\setminus E(L_{0})$.

No matter which case above, by the induction hypothesis, $B^1-F_1$, $B^2-F_2$, $B^3-F_3$ have H-paths $P[x^+,z]$, $P[z^+,w]$, $P[w^+,y^+]$ passing through $L_{1}$, $L_{2}$ and $L_{3}$, respectively. Thus, $P[u,v]\cup P[x^{+},z]\cup P[z^{+},w]\cup P[w^{+},y^{+}]+\{(w,w^{+}),(x,x^{+}),(y,y^{+}),(z,z^{+})\}-(x,y)$ is a desired H-path of $BH_n-F$.

{\it Case 2.}  $i=1$.

By Lemma \ref{le-8}, there is a neighbor $z$ of $x^+$ in $B^{1}$ such that $(x^+,z)\notin E(L_1)$, $L_1+(x^+,z)$ is a linear forest, $\{u,v\}$ is compatible to $L_1+(x^+,z)$, and $z^+$ or $z^-$, say $z^+$, is not an internal vertex of $L_2$. Note that $|E(L_{1}+\{(x^{+},z\})\cup F_{1}|\leq 2n-4$. By the induction hypothesis, $B^{1}-F_{1}$ has a H-path $P[u,v]$ passing through $L_{1}+(x^{+},z)$. Lemma \ref{le-2} implies that there are vertices $y\in V_0\cap Y$ and $w\in V_2\cap X$ such that $y$ (resp. $y^{+}$) is incident with none of $E(L_0)$ (resp. $E(L_3)$), and $w$ (resp. $w^{+}$) is incident with none of $E(L_2)$ (resp. $E(L_3)$). By the induction hypothesis, $B^{0}-F_{0}$, $B^{2}-F_{2}$, $B^{3}-F_{3}$ have H-paths $P[x,y]$, $P[z^{+},w]$, $P[w^{+},y^{+}]$ passing through $L_{0}$, $L_{2}$ and $L_{3}$, respectively. Thus, $P[x,y]\cup P[u,v]\cup P[z^{+},w]\cup P[w^{+},y^{+}]+\{(w,w^{+}),(x,x^{+}),(y,y^{+}),(z,z^{+})\}-(x^{+},z)$ is a H-path of $BH_{n}-F$ passing through $L$.

{\it Case 3.}  $i=2$.

By Lemma \ref{le-2}, there are vertices $y\in V_0\cap Y$ and $z\in V_1\cap X$ such that $y$ (resp. $y^+$) is incident with none of $E(L_0)$ (resp. $E(L_3)$), and $z$ (resp. $z^+$) is incident with none of $E(L_1)$ (resp. $E(L_2)$). By Lemma \ref{le-8}, there is a neighbor $w$ of $z^{+}$ in $B^{2}$ such that $(z^+,w)\notin E(L_2)$, $L_{2}+(z^{+},w)$ is a linear forest, $\{u,v\}$ is compatible to $L_{2}+(z^{+},w)$, and $w^+$ or $w^-$, say $w^+$, is not an internal vertex of $L_3$. Note that $|E(L_{2}+\{(z^{+},w)\})\cup F_{2}|\leq 2n-4$. By the induction hypothesis, $B^{0}-F_{0}$, $B^{1}-F_{1}$, $B^{2}-F_{2}$, $B^{3}-F_{3}$ have H-paths $P[x,y]$, $P[x^{+},z]$, $P[u,v]$, $P[w^{+},y^{+}]$ passing through $L_{0}$, $L_{1}$, $L_{2}+(z^{+},w)$ and $L_{3}$, respectively. Thus, $P[x,y]\cup P[x^{+},z]\cup P[u,v]\cup P[w^{+},y^{+}]+\{(w,w^{+}),(x,x^{+}),(y,y^{+}),(z,z^{+})\}-(z^{+},w)$ is a H-path of $BH_{n}-F$ passing through $L$.

{\it Case 4.}  $i=3$.

By Lemma \ref{le-2}, there are vertices $y\in V_0\cap Y$ and $z\in V_1\cap X$ such that $y$ (resp. $y^{+}$) is incident with none of $E(L_{0})$ (resp. $E(L_{3})$), and $z$ (resp. $z^{+}$) is incident with none of $E(L_{1})$ (resp. $E(L_{2})$). By Lemma \ref{le-8}, there is a neighbor $w$ of $y^{+}$ such that $(y^+,w)\notin E(L_3)$, $L_{3}+(y^{+},w)$ is a linear forest, $\{u,v\}$ is compatible to $L_{3}+(y^{+},w)$, and $w^+$ or $w^-$, say $w^+$, is not an internal vertex of $L_2$. Note that $|E(L_{3}+\{(y^{+},w)\})\cup F_{3}|\leq 2n-4$. By the induction hypothesis, $B^{0}-F_{0}$, $B^{1}-F_{1}$, $B^{2}-F_{2}$, $B^{3}-F_{3}$ have H-paths $P[x,y]$, $P[x^{+},z]$, $P[z^{+},w^{+}]$, $P[u,v]$ passing through $L_{0}$, $L_{1}$, $L_{2}$ and $L_{3}+(y^{+},w)$, respectively. Thus, $P[x,y]\cup P[x^{+},z]\cup P[z^{+},w^{+}]\cup P[u,v]+\{(w,w^{+}),(x,x^{+}),(y,y^{+}),(z,z^{+})\}-(y^{+},w)$ is a H-path of $BH_{n}-F$ passing through $L$.
\end{proof}









\begin{lemma}

If $|E(L_{0})\cup F_{0}|\in \{2n-4,2n-5\}$ and $u\in V_{i}$, $v\in V_{j}$ for $i\in N_{4}$, $j\in N_4\setminus \{i\}$ then $BH_{n}-F$ has a H-path $P[u,v]$ passing through $L$.
\end{lemma}

\begin{proof}
In this case, $|E(L_{k})\cup F_{k}|\leq 2$, for each $k\in N_{4}\setminus \{0\}$. Without loss of generality, assume that $j>i$.

{\it Case 1.}  $l=0$.

{\it Case 1.1.}  $i=0$ and $j=1$.

{\it Case 1.1.1.}  $x$ is incident with none of $E(L_0)$.

By Lemma \ref{le-2}, there is an $a\in V_0\cap Y$ such that $a$ and $a^{\pm}$ are incident with none of $E(L_0)$ and $E(L_3)$, respectively. By the induction hypothesis, $B^0-F_0$ has a H-path $P[u,a]$ passing through $L_0$. Let $y$ be the neighbor of $x$ on $P[u,a]$, if $u=x$; and let $y$ be the neighbor of $x$ on the segment of $P[u,a]$ between $u$ and $x$, otherwise. Then $y\neq a$. Since $|E(L_{3})|\leq 2$, $y^+$ or $y^-$, say $y^+$, is not an internal vertex of $L_3$. 

If $E(L_{1})\cup F_{1}=\emptyset$, $|E(L_{m})\cup F_{m}|\leq 2$ for some $m\in \{2,3\}$. By Lemma \ref{le-2}, there is a $b\in V_3\cap Y$ such that $b$ (resp. $b^+$) is incident with none of $E(L_3)$ (resp. $E(L_2)$). By the induction hypothesis, $B^3-F_3$ has a H-path $P[y^+,b]$ passing through $L_3$. Let $c$ be the neighbor of $a^+$ on the segment of $P[y^+,b]$ between $y^+$ and $a^+$. Since $|E(L_2)|\leq 2$, $c^+$ or $c^-$, say $c^+$, is not an internal vertex of $L_2$. Let $d\in V_2\cap Y$ such that $d$ is incident with none of $E(L_2)$. By the induction hypothesis, $B^2-F_2$ has a H-path $P[c^+,d]$ passing through $L_2$. Let $z$ be the neighbor of $b^+$ on the segment of $P[c^+,d]$ between $b^+$ and $c^+$.

Suppose first that $x^+\neq v$. By Theorem \ref{th-cheng2014}, there exist two vertex-disjoint paths $P[z^{+},x^+]$ and $P[d^+,v]$ in $B^1$ such that each vertex of $B^1$ lies on one of the two paths. Therefore, $P[u,a]\cup P[z^{+},x^+]\cup P[d^+,v]\cup P[c^+,d]\cup P[y^+,b]+\{(a,a^+),(b,b^+),(c,c^+),(d,d^+),(x,x^+),(y,y^+),(z,z^+)\}-\{(x,y),(b^+,z),(a^+,c)\}$ is a H-path of $BH_{n}-F$ passing through $L$.

Suppose now that $x^+=v$. By Theorem \ref{th-lv2014}, $B^1-\{v\}$ has a H-path $P[z^+,d^+]$. Thus, $P[u,a]\cup P[z^{+},d^+]\cup P[c^+,d]\cup P[y^+,b]+\{(a,a^+),(b,b^+),(c,c^+),$ $(d,d^+),(x,v),(y,y^+),(z,z^+)\}-\{(x,y),(b^+,z),(a^+,c)\}$ is a H-path of $BH_{n}-F$ passing through $L$.

If $E(L_{2})\cup F_{2}=\emptyset$, $|E(L_{m})\cup F_{m}|\leq 2$ for some $m\in \{1,3\}$. Let $b\in V_3\cap Y$ such that $b$ is incident with none of $E(L_3)$. By the induction hypothesis, $B^3-F_3$ has a H-path $P[y^+,b]$ passing through $L_3$. Let $c$ be the neighbor of $a^+$ on the segment of $P[y^+,b]$ between $y^+$ and $a^+$.

Suppose first that $x^+$ is incident with none of $E(L_1)$.
Let $d\in V_1\cap X$ such that $d$ is incident with none of $E(L_1)$. By the induction hypothesis, $B^1-F_1$ has a H-path $P[v,d]$ passing through $L_1$. Let $z$ be the neighbor of $x^+$ on $P[v,d]$, if $v=x^+$; and let $z$ be the neighbor of $x^+$ on the segment of $P[v,d]$ between $v$ and $x^+$, otherwise. Then $z\neq d$. By Theorem \ref{th-cheng2014}, there exist two vertex-disjoint paths $P[z^{+},b^+]$ and $P[d^+,c^+]$ in $B^2$ such that each vertex of $B^2$ lies on one of the two paths. Therefore, $P[u,a]\cup P[v,d]\cup P[z^{+},b^+]\cup P[d^+,c^+]\cup P[y^+,b]+\{(a,a^+),(b,b^+),(c,c^+),(d,d^+),(x,x^+),(y,y^+),(z,z^+)\}-\{(x,y),(x^+,z),(a^+,c)\}$ is a H-path of $BH_{n}-F$ passing through $L$.

Suppose second that $L_1$ has a maximal path $P[x^+,w]$ with $w\in X$. Then $v\neq x^+$. By the induction hypothesis, $B^1-F_1$ has a H-path $P[v,w]$ passing through $L_1$. Let $z$ be the neighbor of $x^+$ on $P[v,w]$ such that $z\neq w$. By Theorem \ref{th-cheng2014}, there exist two vertex-disjoint paths $P[z^{+},b^+]$ and $P[w^+,c^+]$ in $B^2$ such that each vertex of $B^2$ lies on one of the two paths. Therefore, $P[u,a]\cup P[v,w]\cup P[z^{+},b^+]\cup P[w^+,c^+]\cup P[y^+,b]+\{(a,a^+),(b,b^+),(c,c^+),(w,w^+),(x,x^+),(y,y^+),(z,z^+)\}-\{(x,y),(x^+,z),(a^+,c)\}$ is a H-path of $BH_{n}-F$ passing through $L$.

Suppose third that $L_1$ has a maximal path $P[x^+,v]$ with $v\neq x^+$. Let $d\in V_1\cap X$ such that $d$ is incident with none of $E(L_1)$. By the induction hypothesis, $B^1-F_1$ has a H-path $P[v,d]$ passing through $L_1$. Let $(x^+,z)\in E(P[v,d])\setminus E(L_1)$. Then $z\neq d$. By Theorem \ref{th-cheng2014}, there exist two vertex-disjoint paths $P[z^{+},b^+]$ and $P[d^+,c^+]$ in $B^2$ such that each vertex of $B^2$ lies on one of the two paths. Therefore, $P[u,a]\cup P[v,d]\cup P[z^{+},b^+]\cup P[d^+,c^+]\cup P[y^+,b]+\{(a,a^+),(b,b^+),(c,c^+),(d,d^+),(x,x^+),(y,y^+),(z,z^+)\}-\{(x,y),(x^+,z),(a^+,c)\}$ is a H-path of $BH_{n}-F$ passing through $L$.

Suppose now that $L_1$ has a maximal path $P[x^+,w]$ with $w\in Y\setminus \{x^+,v\}$. Let $(x^+,h)\in E(P[x^+,w])$. By Theorem \ref{th-yang2019}, $B^1-F_1$ has a H-path $P[v,h]$ passing through $L_1-(x^+,h)$. Let $d,z\in N_B^1(x^+)\setminus \{h\}$ and $d\neq z$. By Theorem \ref{th-cheng2014}, there exist two vertex-disjoint paths $P[z^{+},b^+]$ and $P[d^+,c^+]$ in $B^2$ such that each vertex of $B^2$ lies on one of the two paths. Therefore, $P[u,a]\cup P[v,h]\cup P[z^{+},b^+]\cup P[d^+,c^+]\cup P[y^+,b]+\{(x^+,h),(a,a^+),(b,b^+),(c,c^+),(d,d^+),$ $(x,x^+),(y,y^+),(z,z^+)\}-\{(x,y),(x^+,z),(x^+,d),(a^+,c)\}$ is a H-path of $BH_{n}-F$ passing through $L$.


If $E(L_{3})\cup F_{3}=\emptyset$, $|E(L_{m})\cup F_{m}|\leq 2$ for each $m\in \{1,2\}$. By Lemma \ref{le-2}, there is a $d\in V_1\cap X$ such that $d$ and $d^{\pm}$ are incident with none of $E(L_1)$ and $E(L_2)$, respectively. Let $b\in V_2\in X$ such that $b$ is incident with none of $E(L_2)$.

Suppose first that $x^+$ is incident with none of $E(L_1)$. By the induction hypothesis, $B^1-F_1$ has a H-path $P[v,d]$ passing through $L_1$. Let $z$ be the neighbor of $x^+$ on $P[v,d]$, if $v=x^+$; and let $z$ be the neighbor of $x^+$ on the segment of $P[v,d]$ between $v$ and $x^+$, otherwise. Since $|E(L_2)|\leq 2$, $z^+$ or $z^-$, say $z^+$, is not an internal vertex of $L_2$. Then $z\neq d$. By the induction hypothesis, $B^2-F_2$ has a H-path $P[z^+,b]$ passing through $L_2$.
Let $c$ be the neighbor of $d^+$ on the segment of $P[z^+,b]$ between $d^+$ and $z^+$. By Theorem \ref{th-cheng2014}, there exist two vertex-disjoint paths $P[a^+,c^+]$ and $P[y^+,b^+]$ in $B^3$ such that each vertex of $B^3$ lies on one of the two paths. Thus, $P[u,a]\cup P[v,d]\cup P[z^{+},b]\cup P[a^+,c^+]\cup P[y^+,b^+]+\{(a,a^+),(b,b^+),(c,c^+),(d,d^+),(x,x^+),(y,y^+),(z,z^+)\}-\{(x,y),$ $(x^+,z),(d^+,c)\}$ is a H-path of $BH_{n}-F$ passing through $L$.

Suppose second that $L_1$ has a maximal path $P[x^+,w]$ with $w\in X$. In this scenario, $|E(L_{2})\cup F_{2}|\leq 1$. By the induction hypothesis, $B^1-F_1$ has a H-path $P[v,w]$ passing through $L_1$. Since $|E(L_2)|\leq 1$, $w^+$ or $w^-$, say $w^+$, is incident with none of $E(L_2)$. Let $z$ be the neighbor of $x^+$ on $P[v,w]$ such that $z\neq w$. By Theorem \ref{th-yang2019}, $B^2-F_2$ has a H-path $P[z^+,b]$ passing through $L_2$. Let $c$ be the neighbor of $w^+$ on the segment of $P[z^+,b]$ between $w^+$ and $z^+$. Theorem \ref{th-cheng2014} implies that there exist two vertex-disjoint paths $P[a^+,c^+]$ and $P[y^+,b^+]$ in $B^3$ such that each vertex of $B^3$ lies on one of the two paths. Thus, $P[u,a]\cup P[v,w]\cup P[z^{+},b]\cup P[y^+,b^+]\cup P[a^+,c^+]+\{(a,a^+),(b,b^+),(c,c^+),(w,w^+),(x,x^+),(y,y^+),(z,z^+)\}-\{(x,y),(x^+,z),(w^+,c)\}$ is a H-path of $BH_{n}-F$ passing through $L$.

Suppose third that $L_1$ has a maximal path $P[x^+,w]$ with $w\in Y\setminus \{x^+,v\}$. In this scenario, let $(x^+,h)$ be the edge of $P[x^+,w]$. Then $v\neq x^+$, $F_1=\emptyset$ and $E(L_{2})\cup F_{2}=\emptyset$. Theorem \ref{th-yang2019} implies that $B^1-F_1$ has a H-path $P[v,h]$ passing through $L_1-(x^+,h)$. Let $z$ and $d$ be two neighbors of $x^+$ on $P[v,s]$ such that $z\neq d$. Exactly one of $z$ and $d$, say $z$, lies on the segment of $P[v,h]$ between $x^+$ and $v$. Theorem \ref{th-xu2007} implies that $B^2$ has a H-path $P[z^+,b]$. Let $c$ be the neighbor of $t^+$ on the segment of $P[z^+,b]$ between $t^+$ and $z^+$. Theorem \ref{th-cheng2014} implies that there exist two vertex-disjoint paths $P[a^+,c^+]$ and $P[y^+,b^+]$ in $B^3$ such that each vertex of $B^3$ lies on one of the two paths. Thus, $P[u,a]\cup P[v,h]\cup P[z^{+},b]\cup P[y^+,b^+]\cup P[a^+,c^+]+\{(x^+,h),(a,a^+),(b,b^+),(c,c^+),(d,d^+),(x,x^+),(y,y^+),(z,z^+)\}-\{(x,y),(x^+,z),$ $(x^+,d),(d^+,c)\}$ is a H-path of $BH_{n}-F$ passing through $L$.

Suppose now that $L_1$ has a maximal path $P[x^+,v]$ with $v\neq x^+$. In this scenario, $u\neq x$ and $E(L_{2})\cup F_{2}=\emptyset$. By the induction hypothesis, $B^1-F_1$ has a H-path $P[v,d]$ passing through $L_1$. Let $(x^+,z)\in E(P[v,d])\setminus E(L_1)$. By Theorem \ref{th-xu2007} implies that $B^2$ has a H-path $P[z^+,b]$. Let $c$ be the neighbor of $d^+$ on the segment of $P[z^+,b]$ between $d^+$ and $z^+$. Theorem \ref{th-cheng2014} implies that there exist two vertex-disjoint paths $P[a^+,c^+]$ and $P[y^+,b^+]$ in $B^3$ such that each vertex of $B^3$ lies on one of the two paths. Thus, $P[u,a]\cup P[v,d]\cup P[z^{+},b]\cup P[y^+,b^+]\cup P[a^+,c^+]+\{(a,a^+),(b,b^+),(c,c^+),(d,d^+),(x,x^+),(y,y^+),(z,z^+)\}-\{(x,y),(x^+,z),(d^+,c)\}$ is a H-path of $BH_{n}-F$ passing through $L$.

{\it Case 1.1.2.}  $L_0$ has a maximal path $P[x,r]$ with $r\in Y$ and $|E(L_{0})\cup F_{0}|=2n-4$.

In this scenario, $u\neq x$ and $\{u,r\}$ is compatible to $L_0$. By the induction hypothesis, $B^0-F_0$ has a H-path $P[u,r]$ passing through $L_0$. Let $y$ be the neighbor of $x$ on the segment of $P[u,r]$ between $u$ and $x$. Then $(x,y)\notin E(L_0)$. Since $|E(L_3)|\leq 1$, $r^+$ or $r^-$, say $r^+$, is incident with none of $E(L_3)$.


If $|E(L_{3})\cup F_{3}|=1$, then $E(L_{m})\cup F_{m}=\emptyset$ for each $m\in \{1,2\}$. Let $d\in V_1\cap X$. By Theorem \ref{th-xu2007}, $B^1$ has a H-path $P[v,d]$. Let $(x^+,z)\in E(P[v,d])$. Let $h=d^-$ if $z=d$; and $h=d^+$, otherwise. Then $h\neq z^+$. Let $b\in V_3\cap Y$ such that $b$ is incident with none of $E(L_3)$. By Theorem \ref{th-yang2019}, $B^3-F_3$ has a H-path $P[y^+,b]$ passing through $L_3$. Let $c$ be the neighbor of $r^+$ on the segment of $P[y^+,b]$ between $y^+$ and $r^+$. By Theorem \ref{th-cheng2014}, there exist two vertex-disjoint paths $P[z^{+},b^+]$ and $P[h,c^+]$ in $B^2$ such that each vertex of $B^2$ lies on one of the two paths. Thus, $P[u,r]\cup P[v,d]\cup P[z^{+},b^+]\cup P[h,c^+]\cup P[y^+,b]+\{(b,b^+),(c,c^+),(d,h),(r,r^+),(x,x^+),(y,y^+),(z,z^+)\}-\{(x,y),(x^+,z),(r^+,c)\}$ is a H-path of $BH_{n}-F$ passing through $L$.

If $E(L_{3})\cup F_{3}=\emptyset$, then $|E(L_{m})\cup F_{m}|\leq 1$ for each $m\in \{1,2\}$.

Suppose first that $x^+$ is incident with none of $E(L_1)$. By Lemma \ref{le-2}, there are vertices $d\in V_1\cap X$ and $b\in V_2\cap X$ such that $d$ (resp. $d^+$) is incident with none of $E(L_1)$ (resp. $E(L_2)$) and $b$ is incident with none of $E(L_2)$. By Theorem \ref{th-yang2019}, $B^1-F_1$ has a H-path $P[v,d]$ passing through $L_1$. Let $(x^+,z)\in E(P[v,d])$. Let $h=d^-$ if $z=d$; and $h=d^+$, otherwise. Then $h\neq z^+$.
Theorem \ref{th-yang2019} implies that $B^2-F_2$ has a H-path $P[z^+,b]$ passing through $L_2$.
Let $c$ be the neighbor of $h$ on the segment of $P[z^+,b]$ between $z^+$ and $h$. By Theorem \ref{th-cheng2014}, there exist two vertex-disjoint paths $P[r^+,c^+]$ and $P[y^+,b^+]$ in $B^3$ such that each vertex of $B^3$ lies on one of the two paths. Thus, $P[u,r]\cup P[v,d]\cup P[z^{+},b]\cup P[r^+,c^+]\cup P[y^+,b^+]+\{(b,b^+),(c,c^+),(d,h),(r,r^+),(x,x^+),(y,y^+),(z,z^+)\}-\{(x,y),(x^+,z),(r^+,c)\}$ is a H-path of $BH_{n}-F$ passing through $L$.

Suppose now that $x^+$ is incident with some edge of $E(L_1)$. In this scenario, let $(x^+,w)$ be the edge of $L_1$. Then $v\neq x^+$, $F_1=\emptyset$ and $E(L_{2})\cup F_{2}=\emptyset$. Theorem \ref{th-yang2019} implies that $B^1$ has a H-path $P[v,w]$ passing through $L_1$. Let $z$ be the neighbor of $x^+$ on $P[v,w]$ such that $z\neq w$. Let $b\in V_2\cap X$. Theorem \ref{th-xu2007} implies that $B^2$ has a H-path $P[w^+,b]$. Let $c$ be the neighbor of $z^+$ on the segment of $P[w^+,b]$ between $w^+$ and $z^+$. Theorem \ref{th-cheng2014} implies that there exist two vertex-disjoint paths $P[r^+,b^+]$ and $P[y^+,c^+]$ in $B^3$ such that each vertex of $B^3$ lies on one of the two paths. Therefore, $P[u,r]\cup P[v,w]\cup P[w^{+},b]\cup P[y^+,c^+]\cup P[r^+,b^+]+\{(b,b^+),(c,c^+),(r,r^+),(w,w^+),(x,x^+),(y,y^+),(z,z^+)\}$ $-\{(x,y),(x^+,z),(z^+,c)\}$ is a H-path of $BH_{n}-F$ passing through $L$.

{\it Case 1.1.3.}  $L_0$ has a maximal path $P[x,r]$ with $r\notin \{x,u\}$ and $|E(L_{0})\cup F_{0}|=2n-5$.

Let $(x,s)\in E(P[x,r])$. By Lemma \ref{le-12}, there are two distinct vertices $y,t\in N_{B^0}(x)\setminus \{s\}$ such that $L_0+\{(x,y),(x,t)\}-(x,s)$ is a linear forest and $y^{\pm}$ are incident with none of $E(L_3)$. Note that $\{u,s\}$ is compatible to $L_0+\{(x,y),(x,t)\}-(x,s)$ and $|E(L_0+\{(x,y),(x,t)\}-(x,s))\cup F_0|=2n-4$. By the induction hypothesis, $B^0-F_0$ has a H-path $P[u,s]$ passing through $L_0+\{(x,y),(x,t)\}-(x,s)$. Since $|E(L_3)|\leq 2$, $t^+$ or $t^-$, say $t^+$, is not an internal vertex of $L_3$.


If $E(L_{1})\cup F_{1}=\emptyset$, $|E(L_{m})\cup F_{m}|\leq 2$ for some $m\in \{2,3\}$. By Lemma \ref{le-2}, there is a $b\in V_3\cap Y$ such that $b$ (resp. $b^+$) is incident with none of $E(L_3)$ (resp. $E(L_2)$). By the induction hypothesis, $B^3-F_3$ has a H-path $P[t^+,b]$ passing through $L_3$. Let $c$ be the neighbor of $y^+$ on the segment of $P[t^+,b]$ between $y^+$ and $t^+$. Since $|E(L_2)|\leq 2$, $c^+$ or $c^-$, say $c^+$, is not an internal vertex of $L_2$. By Lemma \ref{le-2}, there is a $d\in V_1\cap X$ such that $d$ and $d^{\pm}$ are incident with none of $E(L_1)$ and $E(L_2)$, respectively. By the induction hypothesis, $B^2-F_2$ has a H-path $P[c^+,d^+]$ passing through $L_2$. Let $z$ be the neighbor of $b^+$ on the segment of $P[c^+,d^+]$ between $b^+$ and $c^+$. Since $z\neq d^+$, $z^-\neq (d^+)^-$ (i.e. $d$).

Suppose first that $x^+\neq v$ and $y$ lies on the segment of $P[u,s]$ between $x$ and $u$. Theorem \ref{th-cheng2014} implies that there exist two vertex-disjoint paths $P[x^+,d]$ and $P[v,z^-]$ in $B^1$ such that each vertex of $B^1$ lies on one of the two paths. Therefore, $P[u,s]\cup P[v,z^-]\cup P[x^{+},d]\cup P[c^+,d^+]\cup P[t^+,b]+\{(x,s),(b,b^+),(c,c^+),(d,d^+),(t,t^+),$ $(x,x^+),(y^+,y),(z,z^-)\}-\{(x,y),$ $(x,t),$ $(b^+,z),(y^+,c)\}$ is a H-path of $BH_{n}-F$ passing through $L$.

Suppose second that $x^+\neq v$ and $t$ lies on the segment of $P[u,s]$ between $x$ and $u$. By Theorem \ref{th-cheng2014}, there exist two vertex-disjoint paths $P[x^+,z^-]$ and $P[v,d]$ in $B^1$ such that each vertex of $B^1$ lies on one of the two paths. Thus, $P[u,s]\cup P[v,d]\cup P[x^{+},z^-]\cup P[c^+,d^+]\cup P[t^+,b]+\{(x,s),(b,b^+),(c,c^+),(d,d^+),(t,t^+),(x,x^+),$ $(y,y^+),(z,z^-)\}-\{(x,y),$ $(x,t),$ $(b^+,z),(y^+,c)\}$ is a H-path of $BH_{n}-F$ passing through $L$.

Suppose now that $x^+=v$. By Theorem \ref{th-lv2014}, $B^1-\{v\}$ has a H-path $P[z^-,d]$. Thus, $P[u,s]\cup P[z^-,d]\cup P[c^+,d^+]\cup P[t^+,b]+\{(x,s),(b,b^+),(c,c^+),$ $(d,d^+),(t,t^+),(x,v),(y,y^+),(z,z^-)\}-\{(x,y),(x,t),(z,$ $b^+),(y^+,c)\}$ is a H-path of $BH_{n}-F$ passing through $L$.

If $E(L_{2})\cup F_{2}=\emptyset$, $|E(L_{m})\cup F_{m}|\leq 2$ for some $m\in \{1,3\}$. Let $b\in V_3\cap Y$ such that $b$ is incident with none of $E(L_3)$. By the induction hypothesis, $B^3-F_3$ has a H-path $P[t^+,b]$ passing through $L_3$. Let $c$ be the neighbor of $y^+$ on the segment of $P[t^+,b]$ between $y^+$ and $t^+$. Let $d\in V_1\cap X$ such that $d$ is incident with none of $E(L_1)$. By the induction hypothesis, $B^1-F_1$ has a H-path $P[v,d]$ passing through $L_1$. Let $(x^+,z)\in E(P[v,d])\setminus E(L_1)$. Let $g=d^-$, if $z=d$; and $g=d^+$, otherwise. Then $g\neq z^+$.

Suppose first that $y$ lies on the segment of $P[u,s]$ between $x$ and $u$. By Theorem \ref{th-cheng2014}, there exist two vertex-disjoint paths $P[z^+,c^+]$ and $P[g,b^+]$ in $B^2$ such that each vertex of $B^2$ lies on one of the two paths. Therefore, $P[u,s]\cup P[v,d]\cup P[z^{+},c^+]\cup P[g,b^+]\cup P[t^+,b]+\{(x,s),(b,b^+),(c,c^+),(d,g),(t^+,$ $t),(x,x^+),(y,y^+),(z,z^+)\}-\{(x,y),(x,t),(z,$ $x^+),(y^+,c)\}$ is a H-path of $BH_{n}-F$ passing through $L$.

Suppose second that $t$ lies on the segment of $P[u,s]$ between $x$ and $u$. By Theorem \ref{th-cheng2014}, there exist two vertex-disjoint paths $P[z^+,b^+]$ and $P[g,c^+]$ in $B^2$ such that each vertex of $B^2$ lies on one of the two paths. Therefore, $P[u,s]\cup P[v,d]\cup P[z^{+},b^+]\cup P[g,c^+]\cup P[t^+,b]+\{(x,s),(b,b^+),(c,c^+),(d,g),(t^+,$ $t),(x,x^+),(y,y^+),(z,z^+)\}-\{(x,y),(x,t),(z,$ $x^+),(y^+,c)\}$ is a H-path of $BH_{n}-F$ passing through $L$.

If $E(L_{3})\cup F_{3}=\emptyset$, $|E(L_{m})\cup F_{m}|\leq 2$ for each $m\in \{1,2\}$.
By Lemma \ref{le-2}, there is a $d\in V_1\cap X$ such that $d$ and $d^{\pm}$ are incident with none of $E(L_1)$ and $E(L_2)$, respectively. By the induction hypothesis, $B^1-F_1$ has a H-path $P[v,d]$ passing through $L_1$. Let $(x^+,z)\in E(P[v,d])\setminus E(L_1)$. Since $|E(L_2)|\leq 2$, $z^+$ or $z^-$, say $z^+$, is not an internal vertex of $L_2$. Let $g=d^-$, if $z=d$; and $g=d^+$, otherwise. Then $g\neq z^+$. Let $b\in V_2\cap X$ such that $b$ is incident with none of $E(L_2)$. By the induction hypothesis, $B^2-F_2$ has a H-path $P[z^+,b]$ passing through $L_2$. Let $c$ be the neighbor of $g$ on the segment of $P[z^+,b]$ between $g$ and $z^+$.

Suppose first that $y$ lies on the segment of $P[u,s]$ between $x$ and $u$. By Theorem \ref{th-cheng2014}, there exist two vertex-disjoint paths $P[y^+,b^+]$ and $P[t^+,c^+]$ in $B^3$ such that each vertex of $B^3$ lies on one of the two paths. Thus, $P[u,s]\cup P[v,d]\cup P[z^{+},b]\cup P[y^+,b^+]\cup P[t^+,c^+]+\{(x,s),(b,b^+),(c,c^+),(d,g),(t,t^+),(x^+,$ $x),(y,y^+),(z,z^+)\}-\{(x,y),(x,t),(z,$ $x^+),(g,c)\}$ is a H-path of $BH_{n}-F$ passing through $L$.

Suppose now that $t$ lies on the segment of $P[u,s]$ between $x$ and $u$. By Theorem \ref{th-cheng2014}, there exist two vertex-disjoint paths $P[y^+,c^+]$ and $P[t^+,b^+]$ in $B^3$ such that each vertex of $B^3$ lies on one of the two paths. Thus, $P[u,s]\cup P[v,d]\cup P[z^{+},b]\cup P[y^+,c^+]\cup P[t^+,b^+]+\{(x,s),(b,b^+),(c,c^+),(d,g),(t,t^+),(x^+,$ $x),(y,y^+),(z,z^+)\}-\{(x,y),(x,t),(z,$ $x^+),(g,c)\}$ is a H-path of $BH_{n}-F$ passing through $L$.

{\it Case 1.1.4.}  $L_0$ has a maximal path $P[x,u]$ with $u\neq x$.

In this scenario, $v\neq x^+$. By Lemma \ref{le-2}, there is a $a\in V_0\cap Y$ such that $a$ and $a^{\pm}$ are incident with $E(L_0)$ and $E(L_3)$, respectively. By the induction hypothesis, $B^0-F_0$ has a H-path $P[u,a]$ passing through $L_0$. Let $(x,y)\in E(P[u,a])\setminus E(L_0)$. Since the length of the segment of $P[u,a]$ between $x$ and $a$ is $|E(P[u,a])|-|E(P[x,u])|\geq (4^{n-1}-1)-(2n-4)\geq 13$, we have $y\neq a$. Since $|E(L_3)|\leq 2$, $y^+$ or $y^-$, say $y^+$, is not an internal vertex of $L_3$. By Lemma \ref{le-2}, there is a $b\in V_3\cap Y$ such that $b$ and $b^{\pm}$ are incident with none of $E(L_3)$ and $E(L_2)$, respectively.

If $E(L_{1})\cup F_{1}=\emptyset$, $|E(L_{m})\cup F_{m}|\leq 2$ for some $m\in \{2,3\}$. By the induction hypothesis, $B^3-F_3$
has a H-path $P[y^+,b]$ passing through $L_3$. Let $c$ be the neighbor of $a^+$ on the segment of $P[y^+,b]$ between $a^+$ and $y^+$. Since $|E(L_2)|\leq 2$, $c^+$ or $c^-$, say $c^+$, is not an internal vertex of $L_2$. Let $d\in V_2\cap Y$ such that $d$ is incident with none of $E(L_2)$. By the induction hypothesis, $B^2-F_2$ has a H-path $P[c^+,d]$ passing through $L_2$. Let $z$ be the neighbor of $b^+$ on the segment of $P[c^+,d]$ between $b^+$ and $z^+$. By Theorem \ref{th-cheng2014}, there exist two vertex-disjoint paths $P[x^+,z^+]$ and $P[v,d^+]$ in $B^1$ such that each vertex of $B^1$ lies on one of the two paths. Thus, $P[u,a]\cup P[v,d^+]\cup P[x^+,z^+]\cup P[c^{+},d]\cup P[y^+,b]+\{(a,a^+),(b,b^+),(c,c^+),(d,d^+),(x,x^+),$ $(y,y^+),(z,z^+)\}-\{(x,y),(b^+,z),(a^+,c)\}$ is a H-path of $BH_{n}-F$ passing through $L$.

If $E(L_{2})\cup F_{2}=\emptyset$, $|E(L_{m})\cup F_{m}|\leq 2$ for some $m\in \{1,3\}$. By the induction hypothesis, $B^3-F_3$
has a H-path $P[y^+,b]$ passing through $L_3$. Let $c$ be the neighbor of $a^+$ on the segment of $P[y^+,b]$ between $a^+$ and $y^+$.

Suppose first that $x^+$ is incident with none of $E(L_1)$.
Let $d\in V_1\cap X$ such that $d$ is incident with none of $E(L_1)$. By the induction hypothesis, $B^1-F_1$ has a H-path $P[v,d]$ passing through $L_1$. Let $z$ be the neighbor of $x^+$ on the segment of $P[v,d]$ between $x^+$ and $v$. By Theorem \ref{th-cheng2014}, there exist two vertex-disjoint paths $P[z^{+},b^+]$ and $P[d^+,c^+]$ in $B^2$ such that each vertex of $B^2$ lies on one of the two paths. Thus, $P[u,a]\cup P[v,d]\cup P[z^{+},b^+]\cup P[d^+,c^+]\cup P[y^+,b]+\{(a,a^+),(b,b^+),(c,c^+),(d,d^+),(x,x^+),(y,y^+),(z,z^+)\}-\{(x,y),(x^+,z),(a^+,c)\}$ is a H-path of $BH_{n}-F$ passing through $L$.

Suppose second that $L_1$ has a maximal path $P[x^+,w]$ with $w\in X$. By the induction hypothesis, $B^1-F_1$ has a H-path $P[v,w]$ passing through $L_1$. Let $z$ be the neighbor of $x^+$ on $P[v,w]$ such that $z\neq w$. By Theorem \ref{th-cheng2014}, there exist two vertex-disjoint paths $P[z^{+},b^+]$ and $P[w^+,c^+]$ in $B^2$ such that each vertex of $B^2$ lies on one of the two paths. Thus, $P[u,a]\cup P[v,w]\cup P[z^{+},b^+]\cup P[w^+,c^+]\cup P[y^+,b]+\{(a,a^+),(b,b^+),(c,c^+),(w,w^+),(x,x^+),(y,y^+),(z,z^+)\}-\{(x,y),(x^+,z),(a^+,c)\}$ is a H-path of $BH_{n}-F$ passing through $L$.

Suppose now that $L_1$ has a maximal path $P[x^+,w]$ with $w\in Y\setminus \{x^+\}$. In this case, $w\neq v$. Let $(x^+,h)\in E(P[x^+,w])$. By Theorem \ref{th-yang2019}, $B^1-F_1$ has a H-path $P[v,h]$ passing through $L_1-(x^+,h)$. Let $z,d\in N_B^1(x^+)$ in $P[v,h]$ and $z\neq d$. By Theorem \ref{th-cheng2014}, there exist two vertex-disjoint paths $P[z^{+},b^+]$ and $P[d^+,c^+]$ in $B^2$ such that each vertex of $B^2$ lies on one of the two paths. Thus, $P[u,a]\cup P[v,h]\cup P[z^{+},b^+]\cup P[d^+,c^+]\cup P[y^+,b]+\{(x^+,h),(a,a^+),(b,b^+),(c,c^+),(d,d^+),(x,x^+),(y,y^+),(z,z^+)\}-\{(x,y),(x^+,z),$ $(x^+,d),(a^+,c)\}$ is a H-path of $BH_{n}-F$ passing through $L$.

If $E(L_{3})\cup F_{3}=\emptyset$, $|E(L_{m})\cup F_{m}|\leq 2$ for each $m\in \{1,2\}$. By Lemma \ref{le-2}, there is a $d\in V_1\cap X$ such that $d$ and $d^{\pm}$ are incident with none of $E(L_1)$ and $E(L_2)$, respectively.

Suppose first that $x^+$ is incident with none of $E(L_1)$.
By the induction hypothesis, $B^1-F_1$ has a H-path $P[v,d]$ passing through $L_1$. Let $z$ be the neighbor of $x^+$ on the segment of $P[v,d]$ between $x^+$ and $v$. Since $|E(L_2)|\leq 2$, $z^+$ or $z^-$, say $z^+$, is not an internal vertex of $L_2$.
By the induction hypothesis, $B^2-F_2$ has a H-path $P[z^+,b^+]$ passing through $L_2$.
Let $c$ be the neighbor of $d^+$ on the segment of $P[z^+,b^+]$ between $d^+$ and $z^+$. Since $c\neq b^+$, $c^-\neq (b^+)^-$ (i.e. $b$). By Theorem \ref{th-cheng2014}, there exist two vertex-disjoint paths $P[a^+,c^-]$ and $P[y^+,b]$ in $B^3$ such that each vertex of $B^3$ lies on one of the two paths. Thus, $P[u,a]\cup P[v,d]\cup P[z^{+},b^+]\cup P[a^+,c^-]\cup P[y^+,b]+\{(a,a^+),(b,b^+),(c,c^-),(d,d^+),(x,x^+),$ $(y,y^+),(z,z^+)\}-\{(x,y),$ $(x^+,z),(d^+,c)\}$ is a H-path of $BH_{n}-F$ passing through $L$.

Suppose second that $L_1$ has a maximal path $P[x^+,w]$ with $w\in X$. In this scenario, $|E(L_{2})\cup F_{2}|\leq 1$. Since $|E(L_2)|\leq 1$, $w^+$ or $w^-$, say $w^+$, is incident with none of $E(L_2)$. Theorem \ref{th-yang2019} implies that $B^1-F_1$ has a H-path $P[v,w]$ passing through $L_1$. Let $z$ be the neighbor of $x^+$ on $P[v,w]$ such that $z\neq w$. By Theorem \ref{th-yang2019}, $B^2-F_2$ has a H-path $P[z^+,b^+]$ passing through $L_2$. Let $c$ be the neighbor of $w^+$ on the segment of $P[z^+,b^+]$ between $w^+$ and $z^+$. And $c^-\neq b$. Theorem \ref{th-cheng2014} implies that there exist two vertex-disjoint paths $P[a^+,c^-]$ and $P[y^+,b]$ in $B^3$ such that each vertex of $B^3$ lies on one of the two paths. Therefore, $P[u,a]\cup P[v,w]\cup P[z^{+},b]\cup P[y^+,b^+]\cup P[a^+,c^-]+\{(a,a^+),(b,b^+),(c,c^-),(w,w^+),(x,x^+),(y,y^+),(z,z^+)\}-\{(x,y),(x^+,z),(w^+,$ $c)\}$ is a H-path of $BH_{n}-F$ passing through $L$.

Suppose now that $L_1$ has a maximal path $P[x^+,w]$ with $w\in Y\setminus \{x^+\}$. Since $\{u,v\}$ is compatible to $L$, $w\neq v$. In this scenario, let $(x^+,h)\in E(P[x^+,w])$. Then $E(L_{2})\cup F_{2}=\emptyset$. By Theorem \ref{th-yang2019}, $B^1-F_1$ has a H-path $P[v,h]$ passing through $L_1-(x^+,h)$. Let $z$ and $d$ be two neighbors of $x^+$ on $P[v,h]$
such that $z\neq d$. By Theorem \ref{th-xu2007}, $B^2$ has a H-path $P[z^+,b^+]$. Let $c$ be the neighbor of $d^+$ on the segment of $P[z^+,b^+]$ between $d^+$ and $z^+$. By Theorem \ref{th-cheng2014}, there exist two vertex-disjoint paths $P[a^+,c^-]$ and $P[y^+,b]$ in $B^3$ such that each vertex of $B^3$ lies on one of the two paths. Therefore, $P[u,a]\cup P[v,h]\cup P[z^{+},b^+]\cup P[y^+,b]\cup P[a^+,c^-]+\{(x^+,h),(a,a^+),(b,b^+),(c,c^-),(d,d^+),(x,x^+),(y,y^+),(z,z^+)\}-\{(x,y),(x^+,z),$ $(x^+,d),(d^+,c)\}$ is a H-path of $BH_{n}-F$ passing through $L$.

{\it Case 1.1.5.}  $L_0$ has a maximal path $P[x,r]$ with $r\in X\setminus \{x,u\}$ and $|E(L_0)\cup F_0|=2n-4$

Let $(x,a)\in E(P[x,r])$. Then $\{u,a\}$ is compatible to $L_0-(x,a)$. By the induction hypothesis, $B^0-F_0$ has a H-path $P[u,a]$ passing through $L_0-(x,a)$. Let $y$ be the neighbor of $x$ on the segment of $P[u,a]$ between $u$ and $x$, and $s$ be the other neighbor of $x$ on $P[u,a]$.

If $|E(L_{3})\cup F_{3}|=1$, then $E(L_{m})\cup F_{m}=\emptyset$ for each $m\in \{1,2\}$, and there is a neighbor of $s$ in $B^3$, say $s^+$, incident with none of $E(L_3)$. Let $d\in V_1\cap X$. By Theorem \ref{th-xu2007}, $B^1$ has a H-path $P[v,d]$. Let $(x^+,z)\in E(v,d)$. Let $h=d^-$ if $z=d$; and $h=d^+$, otherwise. Then $h\neq z^+$. Let $b\in V_3\cap Y$ such that $b$ is incident with none of $E(L_3)$. Theorem \ref{th-yang2019} implies that $B^3-F_3$ has a H-path $P[y^+,b]$ passing through $L_3$. Let $c$ be the neighbor of $s^+$ on the segment of $P[y^+,b]$ between $y^+$ and $s^+$. Theorem \ref{th-cheng2014} implies that there exist two vertex-disjoint paths $P[z^{+},b^+]$ and $P[h,c^+]$ in $B^2$ such that each vertex of $B^2$ lies on one of the two paths. Therefore, $P[u,a]\cup P[v,d]\cup P[z^{+},b^+]\cup P[h,c^+]\cup P[y^+,b]+\{(x,a),(b,b^+),(c,c^+),(d,h),(s,s^+),(x,x^+),(y,y^+),(z,z^+)\}-\{(x,y),$ $(x,s),(x^+,z),(s^+,c)\}$ is a H-path of $BH_{n}-F$ passing through $L$.

If $E(L_{3})\cup F_{3}=\emptyset$, then $|E(L_{m})\cup F_{m}|\leq 1$ for each $m\in \{1,2\}$.

Suppose first that $x^+$ is incident with none of $E(L_1)$. Lemma \ref{le-2} implies that there are vertices $d\in V_1\cap X$ and $b\in V_2\cap X$ such that $d$ (resp. $d^+$) is incident with none of $E(L_1)$ (resp. $E(L_2)$) and $b$ is incident with none of $E(L_2)$. Theorem \ref{th-yang2019} implies that $B^1-F_1$ has a H-path $P[v,d]$ passing through $L_1$. Let $z$ be a neighbor of $x^+$ on $P[v,d]$ such that $(z,x^+)\notin E(L_1)$. Let $h=d^-$ if $z=d$; and $h=d^+$, otherwise. Then $h\neq z^+$. Theorem \ref{th-yang2019} implies that $B^2-F_2$ has a H-path $P[z^+,b]$ passing through $L_2$.
Let $c$ be the neighbor of $h$ on the segment of $P[z^+,b]$ between $z^+$ and $h$. Theorem \ref{th-cheng2014} implies that there exist two vertex-disjoint paths $P[s^+,c^+]$ and $P[y^+,b^+]$ in $B^3$ such that each vertex of $B^3$ lies on one of the two paths. Thus, $P[u,a]\cup P[v,d]\cup P[z^{+},b]\cup P[s^+,c^+]\cup P[y^+,b^+]+\{(x,a),(b,b^+),(c,c^+),(d,h),(s,s^+),(x,x^+),(y,y^+),(z,z^+)\}-\{(x,y),(x,s),(z,$ $x^+),(h,c)\}$ is a H-path of $BH_{n}-F$ passing through $L$.

Suppose now that $x^+$ is incident with some edge of $E(L_1)$. In this scenario, let $(x^+,w)$ be the edge of $L_1$. Then $v\neq x^+$, $F_1=\emptyset$ and $E(L_{2})\cup F_{2}=\emptyset$. Theorem \ref{th-yang2019} implies that $B^1$ has a H-path $P[v,w]$ passing through $L_1$. Let $z$ be the neighbor of $x^+$ on $P[v,w]$ such that $z\neq w$. Let $b\in V_2\cap X$. By Theorem \ref{th-xu2007}, $B^2$ has a H-path $P[w^+,b]$. Let $c$ be the neighbor of $z^+$ on the segment of $P[w^+,b]$ between $w^+$ and $z^+$. By Theorem \ref{th-cheng2014}, there exist two vertex-disjoint paths $P[s^+,b^+]$ and $P[y^+,c^+]$ in $B^3$ such that each vertex of $B^3$ lies on one of the two paths. Thus, $P[u,a]\cup P[v,w]\cup P[w^{+},b]\cup P[y^+,c^+]\cup P[s^+,b^+]+\{(x,a),(b,b^+),(c,c^+),(s,s^+),(w,w^+),(x,x^+),(y,y^+),(z,z^+)\}-\{(x,y),(x,s),(x^+,z),(z^+,c)\}$ is a H-path of $BH_{n}-F$ passing through $L$.

{\it Case 1.2.}  $i=0$ and $j\in \{2,3\}$.

By Lemma \ref{le-2}, there is an $a\in V_0\cap Y$ such that $a$ and $a^{\pm}$ are incident with none of $E(L_0)$ and $E(L_3)$, respectively. By the induction hypothesis, $B^0-F_0$ has a H-path $P[u,a]$ passing through $L_0$. Let $(x,y)\in E(P[u,a])\setminus E(L_0)$. Since $|E(L_3)|\leq 2$, $y^+$ or $y^-$, say $y^+$, is not an internal vertex of $L_3$. Let $g=a^-$ if $y=a$; and $g=a^+$, otherwise. Then $g\neq y^+$. 

{\it Case 1.2.1.} $j=2$.

Suppose first that $E(L_{3})\cup F_{3}\neq \emptyset$. Then $|E(L_{m})\cup F_{m}|\leq 1$ for each $m\in \{1,2\}$. By Lemma \ref{le-2}, there is a $b\in V_3\cap Y$ such that $b$ (resp. $b^+$) is incident with none of $E(L_3)$ (resp. $E(L_2)$). By the induction hypothesis, $B^3-F_3$ has a H-path $P[y^+,b]$ passing through $L_3$. Let $c$ be the neighbor of $g$ on the segment of $P[y^+,b]$ between $g$ and $y^+$. Since $|E(L_2)|\leq 1$, $c^+$ or $c^-$, say $c^+$, is incident with none of $E(L_2)$. By Theorem \ref{th-yang2019}, $B^2-F_2$ has a H-path $P[c^+,v]$ passing through $L_2$. Let $z$ be the neighbor of $b^+$ on the segment of $P[c^+,v]$ between $b^+$ and $c^+$. Since $|E(L_1)|\leq 1$, $z^+$ or $z^-$, say $z^+$, is incident with none of $E(L_1)$. By Theorem \ref{th-yang2019}, $B^1-F_1$ has a H-path $P[x^+,z^+]$ passing through $L_1$. Thus, $P[u,a]\cup P[x^+,z^+]\cup P[c^{+},v]\cup P[y^+,b]+\{(a,g),(b,b^+),(c,c^+),(x,x^+),(y,y^+),(z,z^+)\}-\{(x,y),(b^+,z),(g,c)\}$ is a H-path of $BH_{n}-F$ passing through $L$.

Suppose now that $E(L_{3})\cup F_{3}=\emptyset$. Then $|E(L_{m})\cup F_{m}|\leq 2$ for each $m\in \{1,2\}$. By Lemma \ref{le-2}, there are vertices $d\in V_1\cap X$ and $b\in V_2\cap X$ such that $d$ and $d^{\pm}$ are incident with none of $E(L_1)$ and $E(L_2)$, respectively, and $b$ is incident with none of $E(L_2)$. There is a neighbor of $d$ in $B^2$, say $d^+$, being not $v$. By the induction hypothesis, $B^1-F_1$, $B^2-F_2$ have H-paths $P[x^+,d]$, $P[v,b]$ passing through $L_1$ and $L_2$, respectively. Let $c$ be the neighbor of $d^+$ on the segment of $P[v,b]$ between $d^+$ and $v$. By Theorem \ref{th-cheng2014}, there exist two vertex-disjoint paths $P[g,c^+]$ and $P[y^+,b^+]$ in $B^3$ such that each vertex of $B^3$ lies on one of the two paths. Thus, $P[u,a]\cup P[x^+,d]\cup P[v,b]\cup P[g,c^+]\cup P[y^+,b^+]+\{(a,g),(b,b^+),(c,c^+),(d,d^+),(x,x^+),(y,y^+)\}-\{(x,y),(d^+,c)\}$ is a H-path of $BH_{n}-F$ passing through $L$.

{\it Case 1.2.3.}  $j=3$.

Suppose first that $|E(L_{3})\cup F_{3}|\leq 1$. In this case, $|E(L_3)|\leq 1$, $y^+$ or $y^-$, say $y^+$, is incident with none of $E(L_3)$. By Lemma \ref{le-2}, there is a $d\in V_1\cap X$ such that $d$ (resp. $d^+$) is incident with $E(L_1)$ (resp. $E(L_2)$). Recall that $|E(L_{k})\cup F_{k}|\leq 2$ for each $k\in N_4\setminus \{0\}$. By the induction hypothesis, $B^1-F_1$, $B^3-F_3$ have H-paths $P[x^+,d]$, $P[y^+,v]$ passing through $L_1$ and $L_3$, respectively. Let $c$ be the neighbor of $g$ on the segment of $P[y^+,v]$ between $y^+$ and $g$. By the induction hypothesis, $B^2-F_2$ has a H-path $P[d^+,c^+]$ passing through $L_2$. Thus, $P[u,a]\cup P[x^+,d]\cup P[d^+,c^+]\cup P[y^+,v]+\{(a,g),(c,c^+),(d,d^+),(x,x^+),(y,y^+)\}-\{(x,y),(g,c)\}$ is a H-path of $BH_{n}-F$ passing through $L$.

Suppose now that $|E(L_{3})\cup F_{3}|=2$. In this case, $|E(L_{0})\cup F_{0}|\leq 2n-5$, and $n\geq 4$. By Lemma \ref{le-13} and Lemma \ref{le-8}, there is a $z\in N_B^0(x)$ such that $L_0+(x,z)$ is a linear forest and $z^+$ or $z^-$, say $z^+$, is incident with none of $E(L_3)$. By the induction hypothesis, $B^0-F_0$, $B^3-F_3$ have H-paths $P[u,a]$, $P[z^+,v]$ passing through $L_0+(x,z)$ and $L_3$, respectively. Let $c$ be the neighbor of $g$ on the segment of $P[z^+,v]$ between $z^+$ and $g$, and let $d\in V_1\cap X$. By Theorem \ref{th-xu2007}, $B^1$, $B^2$ have H-paths $P[x^+,d]$ and $P[d^+,c^+]$, respectively. Thus, $P[u,a]\cup P[x^+,d]\cup P[d^+,c^+]\cup P[z^+,v]+\{(a,g),(c,c^+),(d,d^+),(x,x^+),(z,z^+)\}-\{(x,y),(g,c)\}$ is a H-path of $BH_{n}-F$ passing through $L$.



{\it Case 1.3.}  $i\neq 0$.

By Lemma \ref{le-2}, there is an $a\in V_0\cap Y$ such that $a$ and $a^{\pm}$ are incident with none of $E(L_0)$ and $E(L_3)$, respectively. By the induction hypothesis, $B^0-F_0$ has a H-path $P[x,a]$ passing through $L_0$.

{\it Case 1.3.1.}  $i=1$, $j=2$.

If $x^+$ is not adjacent to $u$ or $(x^+,u)\notin E(L_1)$, $\{x^+,u\}$ is compatible to $L_1$. By Lemma \ref{le-2}, there is a $b\in V_2\cap X$ such that $b$ (resp. $b^+$) is incident with none of $E(L_2)$ (resp. $E(L_3)$). Recall that $|E(L_{k})\cup F_{k}|\leq 2$ for $k\in N_4\setminus \{0\}$.
By the induction hypothesis, $B^1-F_1$, $B^2-F_2$, $B^3-F_3$ have H-paths $P[x^+,u]$, $P[v,b]$, $P[a^+,b^+]$ passing through $L_1$, $L_2$ and $L_3$, respectively. Thus, $P[x,a]\cup P[x^+,u]\cup P[v,b]\cup P[a^+,b^+]+\{(a,a^+),(b,b^+),(x,x^+)\}$ is a H-path of $BH_{n}-F$ passing through $L$.

If $(x^+,u)\in E(L_1)$, $E(L_{m})\cup F_{m}=\emptyset$ for some $m\in \{2,3\}$. By Lemma \ref{le-2}, there is a $d\in V_1\cap Y$ such that $d$ and $d^{\pm}$ are incident with none of $E(L_1)$ and $E(L_0)$, respectively. Then $d\neq x^+$, $d^-\neq (x^+)^-$ (i.e. $x$). Let $t$ be neighbor of $d^-$ on the segment of $P[x,a]$ between $d^-$ and $a$. Let $h=a^-$, if $t=a$; and $h=a^+$, otherwise. By the induction hypothesis, $B^1-F_1$ has a H-path $P[u,d]$ passing through $L_1$. Let $z$ be the neighbor of $x^+$ on $P[u,d]$ such that $z\neq u$.

Suppose first that $m=2$. There is a neighbor of $z$ in $B^2$, say $z^+$, being not $v$. Let $b\in V_3\cap Y$ such that $b$ is incident with none of $E(L_3)$. By Theorem \ref{th-yang2019}, $B^3-F_3$ has a H-path $P[t^+,b]$ passing through $L_3$. Let $c$ be the neighbor of $h$ on the segment of $P[t^+,b]$ between $h$ and $t^+$. By \ref{th-cheng2014}, there exist two vertex-disjoint paths $P[z^+,b^+]$ and $P[c^+,v]$ in $B^2$ such that each vertex of $B^2$ lies on one of the two paths. Thus, $P[x,a]\cup P[u,d]\cup P[c^+,v]\cup P[z^+,b^+]\cup P[t^+,b]+\{(a,h),(b,b^+),(c,c^+),(d,d^-),(t,t^+),(x,x^+),(z^+,$ $z)\}-\{(d^-,t),(x^+,z),(h,c)\}$ is a H-path of $BH_{n}-F$ passing through $L$.

Suppose now that $m=3$. Since $|E(L_2)|\leq 1$, $z^+$ or $z^-$, say $z^+$, is incident with none of $E(L_2)$. Let $b\in V_2\cap X$ such that $b$ is incident with none of $E(L_2)$. By Theorem \ref{th-yang2019}, $B^2-F_2$ has a H-path $P[v,b]$ passing through $L_2$. Let $c$ be the neighbor of $z^+$ on the segment of $P[v,b]$ between $z^+$ and $v$. By \ref{th-cheng2014}, there exist two vertex-disjoint paths $P[h,c^+]$ and $P[t^+,b^+]$ in $B^3$ such that each vertex of $B^3$ lies on one of the two paths. Thus, $P[x,a]\cup P[u,d]\cup P[v,b]\cup P[h,c^+]\cup P[t^+,b^+]+\{(a,h),(b,b^+),(c,c^+),(d,d^-),(t,t^+),(x,x^+),(z,z^+)\}-\{(d^-,t),(x^+,z),(z^+,c)\}$ is a H-path of $BH_{n}-F$ passing through $L$.

{\it Case 1.3.2.}  $i=1$, $j=3$.

There are $\lfloor{|E(P[x,a])|}/{2}\rfloor=\lfloor(4^{n-1}-1)/2\rfloor$ edges each of which has the form $(s,t)$ with $s\in X$ and $t\in Y$ such that $t$ lies on the segment of $P[x,a]$ between $x$ and $s$. Since $\lfloor{|E(P[x,a])|}/{2}\rfloor-|E(L_0)|\geq \lfloor(4^{n-1}-1)/2\rfloor-(2n-4)\geq 5$, there are at least such $5$ edges $(s,t)$ on $P[x,a]$ that meats above requirements and furthermore $(s,t)\notin E(L_0)$. Since $|E(L_1)|+|E(L_3)|\leq 2$, there are at most $4$ ($<5$) such edges $(s,t)$ that meats above requirements and $s^+$ or $s^-$ (resp. $t^+$ or $t^-$) is incident with some edge of $E(L_1)$ (resp. $E(L_3)$). Thus, there is an edge $(s,t)\in E(P[x,a])\setminus E(L_0)$ such that $s^{\pm}$ (resp. $t^{\pm}$) are incident with none of $E(L_1)$ (resp. $E(L_3)$). By the induction hypothesis, $B^3-F_3$ has a H-path $P[t^+,v]$ passing through $L_3$. Let $c$ be the neighbor of $a^+$ on the segment of $P[t^+,v]$ between $a^+$ and $t^+$.

Suppose first that $|E(L_{2})\cup F_{2}|\leq 1$. Then $|E(L_{1})\cup F_{1}|\leq 2$. By the induction hypothesis, $B^1-F_1$ has a H-path $P[s^+,u]$ passing through $L_1$. Let $(x^+,z)\in E(P[s^+,u])\setminus E(L_1)$. Since $|E(L_2)|\leq 1$, $z^+$ or $z^-$, say $z^+$, is incident with none of $E(L_2)$. By Theorem \ref{th-yang2019}, $B^2-F_2$ has a H-path $P[z^+,c^+]$ passing through $L_2$. Thus, $P[x,a]\cup P[u,s^+]\cup P[z^+,c^+]\cup P[t^+,v]+\{(a,a^+),(c,c^+),(s,s^+),(t,t^+),(x,x^+),(z,z^+)\}-\{(s,t),(x^+,z),(a^+,c)\}$ is a H-path of $BH_{n}-F$ passing through $L$.

Suppose now that $|E(L_{2})\cup F_{2}|=2$. Then $E(L_{1})\cup F_{1}=\emptyset$. Since $|E(L_2)|\leq 2$, $c^+$ or $c^-$, say $c^+$, is not an internal vertex of $L_2$. Let $z\in V_2\cap Y$ such that $z$ is incident with none of $E(L_2)$. By the induction hypothesis, $B^2-F_2$ has a H-path $P[c^+,z]$ passing through $L_2$. There is a neighbor of $z$ in $B^1$, say $z^+$, being not $u$. By \ref{th-cheng2014}, there exist two vertex-disjoint paths $P[x^+,u]$ and $P[s^+,z^+]$ in $B^1$ such that each vertex of $B^1$ lies on one of the two paths. Thus, $P[x,a]\cup P[x^+,u]\cup P[s^+,z^+]\cup P[c^+,z]\cup P[t^+,v]+\{(a,a^+),(c,c^+),(s,s^+),(t,t^+),(x,x^+),(z,z^+)\}-\{(s,t),(a^+,c)\}$ is a H-path of $BH_{n}-F$ passing through $L$.

{\it Case 1.3.3.}  $i=2$, $j=3$.

By Lemma \ref{le-2}, there is a $b\in V_1\cap X$ such that $b$ (resp. $b^+$) is incident with none of $E(L_1)$ (resp. $E(L_2)$). By the induction hypothesis, $B^1-F_1$, $B^2-F_2$, $B^3-F_3$ have H-paths $P[x^+,b]$, $P[b^+,u]$, $P[a^+,v]$ passing through $L_1$, $L_2$ and $L_3$, respectively. Thus, $P[x,a]\cup P[x^+,b]\cup P[b^+,u]\cup P[a^+,v]+\{(a,a^+),(b,b^+),(x,x^+)\}$ is a H-path of $BH_{n}-F$ passing through $L$.

{\it Case 2.}  $l=1$.

{\it Case 2.1.}  $i=0$.

By Lemma \ref{le-2}, there is an $a\in V_0\cap Y$ such that $a$ and $a^{\pm}$ are incident with none of $E(L_0)$ and $E(L_3)$, respectively. By the induction hypothesis, $B^0-F_0$ has a H-path $P[u,a]$ passing through $L_0$.

{\it Case 2.1.1.}  $j=1$.

If $x$ is not adjacent to $v$ or $(x,v)\notin E(L_1)$, $\{v,x\}$ is compatible to $L_1$. By Lemma \ref{le-2}, there is a $b\in V_2\cap X$ such that $b$ (resp. $b^+$) is incident with none of $E(L_2)$ (resp. $E(L_3)$). By the induction hypothesis, $B^1-F_1$, $B^2-F_2$, $B^3-F_3$ have H-paths $P[v,x]$, $P[x^+,b]$, $P[a^+,b^+]$ passing through $L_1$, $L_2$ and $L_3$, respectively. Thus, $P[u,a]\cup P[v,x]\cup P[x^+,b]\cup P[a^+,b^+]+\{(a,a^+),(b,b^+),(x,x^+)\}$ is a desired H-path in $BH_n-F$.

If $(x,v)\in E(L_1)$, $E(L_{m})\cup F_{m}=\emptyset$ for some $m\in \{2,3\}$.
By Lemma \ref{le-2}, there is a $d\in V_1\cap X$ such that $d$ and $d^{\pm}$ are incident with none of $E(L_1)$ and $E(L_2)$, respectively. Then $d\neq x$. By Lemma \ref{le-11}, there is a $z\in N_{B^1}(x)\setminus \{v\}$ such that $z^+$ or $z^-$, say $z^+$, is incident with none of $E(L_0)$. 
Note that $L_1+(x,z)$ is a linear forest and $\{v,d\}$ is compatible to $L_1+(x,z)$. For $n=3$, $|E(L_1+(x^+,z))\cup F_1|\leq 2$; and $|E(L_1+(x^+,z))\cup F_1|\leq 2n-4$, otherwise. By the induction hypothesis, $B^1-F_1$ has a H-path $P[v,d]$ passing through $L_1+(x,z)$. 
Let $t$ be the neighbor of $z^+$ on the segment of $P[u,a]$ between $z^+$ and $a$. Since $|E(L_3)|\leq 1$, $t^+$ or $t^-$, say $t^+$, is incident with none of $E(L_3)$. Let $g=a^-$ if $t=a$; and $g=a^+$, otherwise.

Suppose first that $m=2$. Let $b\in V_3\cap Y$. By Theorem \ref{th-yang2019}, $B^3-F_3$ has a H-path $P[g,b]$ passing through $L_3$. Let $c$ be the neighbor of $t^+$ on the segment of $P[g,b]$ between $t^+$ and $g$. By Theorem \ref{th-cheng2014}, there exist two vertex-disjoint paths $P[d^+,c^+]$ and $P[x^+,b^+]$ in $B^2$ such that each vertex of $B^2$ lies on one of the two paths. Thus, $P[u,a]\cup P[v,d]\cup P[d^{+},c^+]\cup P[x^+,b^+]\cup P[g,b]+\{(a,g),(b,b^+),(c,c^+),(d,d^+),(t,t^+),(x,x^+),(z,z^+)\}-\{(z^+,t),(x,z),(t^+,c)\}$ is a H-path of $BH_{n}-F$ passing through $L$.

Suppose now that $m=3$. Let $b\in V_2\cap X$ such that $b$ is incident with none of $E(L_2)$. By Theorem \ref{th-yang2019}, $B^2-F_2$ has a H-path $P[x^+,b]$ passing through $L_2$. Let $c$ be the neighbor of $d^+$ on the segment of $P[x^+,b]$ between $d^+$ and $x^+$. By Theorem \ref{th-cheng2014}, there exist two vertex-disjoint paths $P[g,c^+]$ and $P[t^+,b^+]$ in $B^3$ such that each vertex of $B^3$ lies on one of the two paths. Thus, $P[u,a]\cup P[v,d]\cup P[x^{+},b]\cup P[g,c^+]\cup P[t^+,b^+]+\{(a,g),(b,b^+),(c,c^+),(d,d^+),(t,t^+),(x,x^+),(z,z^+)\}-\{(z^+,t),(x,z),(d^+,c)\}$ is a H-path of $BH_{n}-F$ passing through $L$.

{\it Case 2.1.2.}  $j=2$.

There are $\lfloor{|E(P[u,a])|}/{2}\rfloor=\lfloor(4^{n-1}-1)/2\rfloor$ edges each of which has the form $(s,t)$ with $s\in X$ and $t\in Y$ such that $t$ lies on the segment of $P[u,a]$ between $u$ and $s$. Since $\lfloor{|E(P[u,a])|}/{2}\rfloor-|E(L_0)|\geq \lfloor(4^{n-1}-1)/2\rfloor-(2n-4)\geq 5$, there are at least such $5$ edges $(s,t)$ on $P[u,a]$ that meats above requirements and furthermore $(s,t)\notin E(L_0)$. Since $|E(L_1)|+|E(L_3)|\leq 2$, there are at most $4$ ($<5$) such edges $(s,t)$ that meats above requirements and $s^+$ or $s^-$ (resp. $t^+$ or $t^-$) is incident with some edge of $E(L_1)$ (resp. $E(L_3)$). Thus, there is an edge $(s,t)\in E(P[u,a])\setminus E(L_0)$ such that $s^{\pm}$ (resp. $t^{\pm}$) are incident with none of $E(L_1)$ (resp. $E(L_3)$). By the induction hypothesis, $B^1-F_1$ has a H-path $P[s^+,x]$ passing through $L_1$.

Suppose first that $E(L_{3})\cup F_{3}\neq \emptyset$, then $|E(L_{m})\cup F_{m}|\leq 1$ for each $m\in \{1,2\}$. Let $(x^+,z)\in E(B^2)\setminus E(L_2)$. Since $|E(L_3)|\leq 2$, $z^+$ or $z^-$, say $z^+$, is not an internal vertex of $L_3$. By the induction hypothesis, $B^3-F_3$ has a H-path $P[t^+,z^+]$ passing through $L_3$. Let $c$ be the neighbor of $a^+$ on the segment of $P[t^+,z^+]$ between $a^+$ and $t^+$. Since $|E(L_2)|\leq 1$, $c^+$ or $c^-$, say $c^+$, is incident with none of $E(L_2)$. By the induction hypothesis, $B^2-F_2$ has a H-path $P[c^+,v]$ passing through $L_2+(x^+,z)$.
Thus, $P[u,a]\cup P[s^+,x]\cup P[c^+,v]\cup P[t^+,z^+]+\{(a,a^+),(c,c^+),(s,s^+),(t,t^+),(x,x^+),(z,z^+)\}-\{(s,t),(x^+,z),(a^+,$ $c)\}$ is a H-path of $BH_{n}-F$ passing through $L$.

Suppose now that $E(L_{3})\cup F_{3}=\emptyset$, then $|E(L_{m})\cup F_{m}|\leq 2$ for each $m\in \{1,2\}$. Let $b\in V_2\cap X$ such that $b$ is incident with none of $E(L_2)$. By the induction hypothesis, $B^2-F_2$ has a H-path $P[v,b]$ passing through $L_2$. Let $(x^+,z)\in E(P[v,b])\setminus E(L_2)$.
By Theorem \ref{th-cheng2014}, there exist two vertex-disjoint paths $P[a^+,z^+]$ and $P[t^+,b^+]$ in $B^3$ such that each vertex of $B^3$ lies on one of the two paths. Thus, $P[u,a]\cup P[s^+,x]\cup P[v,b]\cup P[a^+,z^+]\cup P[t^+,b^+]+\{(a,a^+),(b,b^+),(s,s^+),(t,t^+),(x,x^+),(z,z^+)\}-\{(s,t),(x^+,z)\}$ is a H-path of $BH_{n}-F$ passing through $L$.

{\it Case 2.1.3.}  $j=3$.

By Lemma \ref{le-3}, there is an edge $(s,t)\in E(P[u,a])\setminus E(L_0)$ for some $s\in X$ and $t\in Y$ such that $s^+$ or $s^-$ (resp. $t^+$ or $t^-$), say $s^+$ (resp. $t^+$), is incident with none of $E(L_1)$ (resp. $L_3$) and $\{s,t\}\cap \{u,a\}=\emptyset$. By the induction hypothesis, $B^1-F_1$ has a H-path $P[s^+,x]$ passing through $L_1$.

Suppose first that $E(L_{3})\cup F_{3}\neq \emptyset$, then $|E(L_{2})\cup F_{2}|\leq 1$. By the induction hypothesis, $B^3-F_3$ has a H-path $P[t^+,v]$ passing through $L_3$.
Let $b$ be the neighbor of $a^+$ on the segment of $P[t^+,v]$ between $t^+$ and $a^+$. Since $|E(L_2)|\leq 1$, $b^+$ or $b^-$, say $b^+$, is incident with none of $E(L_2)$. By Theorem \ref{th-yang2019}, $B^2-F_2$ has a H-path $P[x^+,b^+]$ passing through $L_2$. Thus, $P[u,a]\cup P[s^+,x]\cup P[x^+,b^+]\cup P[t^+,v]+\{(a,a^+),(b,b^+),(s,s^+),(t,t^+),(x,x^+)\}-\{(s,t),(a^+,b)\}$ is a H-path of $BH_{n}-F$ passing through $L$.

Suppose now that $E(L_{3})\cup F_{3}=\emptyset$, then $|E(L_{2})\cup F_{2}|\leq 2$. Let $b\in V_2\cap X$ such that $b$ is incident with none of $E(L_2)$. By the induction hypothesis, $B^2-F_2$ has a H-path $P[x^+,b]$ passing through $L_2$. By Theorem \ref{th-cheng2014}, there exist two vertex-disjoint paths $P[a^+,v]$ and $P[t^+,b^+]$ in $B^3$ such that each vertex of $B^3$ lies on one of the two paths. Thus, $P[u,a]\cup P[s^+,x]\cup P[x^+,b]\cup P[a^+,v]\cup P[t^+,b^+]+\{(a,a^+),(b,b^+),(s,s^+),(t,t^+),(x,x^+)\}-(s,t)$ is a H-path of $BH_{n}-F$ passing through $L$.

{\it Case 2.2.}  $i\neq 0$.

By Lemma \ref{le-2}, there is a $b\in V_0\cap Y$ such that $b$ (resp. $b^+$) is incident with none of $E(L_0)$ (resp. $E(L_3)$).

{{\it Case 2.2.1.}  $i=1$, $j=2$.}


{\it Case 2.2.1.1.} $E(L_{1})\cup F_{1}=\emptyset$.

In this case, $|E(L_{m})\cup F_{m}|\leq 2$ for some $m\in \{2,3\}$. By Lemma \ref{le-2}, there is an $a\in V_0\cap X$ such that $a$ and $a^{\pm}$ are incident with none of $E(L_0)$ and $E(L_1)$, respectively, and a $d\in V_2\cap X$ such that $d$ and $d^{\pm}$ are incident with none of $E(L_2)$ and $E(L_3)$, respectively. Lemma \ref{le-2} implies that there is a $w\in V_3\cap X$ such that $w$ and $w^{\pm}$ are incident with none of $E(L_3)$ and $E(L_0)$, respectively.

If $x\neq u$, 
by the induction hypothesis, $B^2-F_2$ has a H-path $P[v,d]$ passing through $L_2$. Let $(x^+,z)\in E(P[v,d])\setminus E(L_2)$. Since $|E(L_3)|\leq 2$, $z^+$ or $z^-$, say $z^+$, is not an internal vertex of $L_3$.

Suppose first that $n=3$. In this case, $|E(L_0)|\leq 2$. By Theorem \ref{th-yang2019}, $B^3-F_3$ has a H-path $P[z^+,w]$ passing through $L_3$. Let $c$ be the neighbor of $d^+$ on the segment of $P[z^+,w]$ between $z^+$ and $d^+$. Since $|E(L_0)|\leq 2$, $c^+$ or $c^-$, say $c^+$, is not an internal vertex of $E(L_0)$. 

Suppose now that $n\geq 4$. By Lemma \ref{le-9}, there are two neighbors $c$ and $s$ of $d^+$ such that $c^+$ or $c^-$ and $s^+$ or $s^-$ are incident with none of $E(L_0)$. We claim that there is a $w\in V_3\cap X\setminus \{c,s\}$ such that $w$ and $w^{\pm}$ are incident with none of $E(L_3)$ and $E(L_0)$, respectively. The reason is follows. There are $|V_3\cap X\setminus \{c,s\}|-|E(L_3)|=4^{n-1}/2-4$ candidates of $w$. Since $E(L_0)$ has at most $|E(L_0)|$ even end vertices, each of which fails at most two candidates of such $w$. Since $|V_3\cap X\setminus \{c,s\}|-|E(L_3)|-2|E(L_{0})|\geq (4^{n-1}/2-4)-2(2n-4)>0$, the claim holds. Note that $L_3+\{(d^+,c),(d^+,s)\}$ is a linear forest and $|E(L_3+\{(d^+,c),(d^+,s)\})\cup F_3|\leq 4\leq 2n-4$. By the induction hypothesis, $B^3-F_3$ has a H-path $P[z^+,w]$ passing through $L_3+\{(d^+,c),(d^+,s)\}$. Exactly one of $c$ and $t$, say $c$, lies on the segment of $P[z^+,w]$ between $z^+$ and $d^+$. Note that $c^+$ or $c^-$, say $c^+$, is incident with none of $E(L_0)$.

No matter which case above, by the induction hypothesis, $B^0-F_0$ has a H-path $P[a,c^+]$ passing through $L_0$. Let $y$ be the neighbor of $w^+$ on the segment of $P[a,c^+]$ between $w^+$ and $c^+$. By Theorem \ref{th-cheng2014}, there exist two vertex-disjoint paths $P[a^+,u]$ and $P[y^+,x]$ in $B^1$ such that each vertex of $B^1$ lies on one of the two paths. Thus, $P[a,c^+]\cup P[y^+,x]\cup P[a^+,u]\cup P[v,d]\cup P[z^+,w]+\{(a,a^+),(c,c^+),(d,d^+),(w,w^+),(x,x^+),(y,y^+),(z,z^+)\}-\{(w^+,y),(x^+,z),(d^+,c)\}$ is a H-path of $BH_{n}-F$ passing through $L$.

If $x=u$ and $x^+$ is incident with none of $E(L_2)$, then $x^+\neq v$. By Theorem \ref{th-yang2019}, $B^2-F_2$
has a H-path $P[v,d]$ passing through $L_2$. Let $z$ be the neighbor of $x^+$ on the segment of $P[v,d]$ between $x^+$ and $v$. Since $|E(L_3)|\leq 2$, $z^+$ or $z^-$, say $z^+$, is not an internal vertex of $L_3$.

Suppose first that $n=3$. In this case, $|E(L_3)|\leq 1$ and $|E(L_0)|\leq 2$. By Theorem \ref{th-yang2019}, $B^3-F_3$ has a H-path $P[z^+,w]$ passing through $L_3$. Let $c$ be the neighbor of $d^+$ on the segment of $P[z^+,w]$ between $d^+$ and $z^+$. Since $|E(L_0)|\leq 2$, $c^+$ or $c^-$, say $c^+$, is not an internal vertex of $E(L_0)$. 

Suppose now that $n\geq 4$. By Lemma \ref{le-9}, there are two neighbors $c$ and $s$ of $d^+$ such that $c^+$ or $c^-$ and $s^+$ or $s^-$ are incident with none of $E(L_0)$. We claim that there is an $w\in V_3\cap X\setminus \{c,s\}$ such that $w$ and $w^{\pm}$ are incident with none of $E(L_3)$ and $E(L_0)$, respectively. The reason is follows. There are $|V_3\cap X\setminus \{c,s\}|-|E(L_3)|=4^{n-1}/2-4$ candidates of $w$. Since $E(L_0)$ has at most $|E(L_0)|$ even end vertices, each of which fails at most two candidates of such $w$. Since $|V_3\cap X\setminus \{c,s\}|-|E(L_3)|-2|E(L_{0})|\geq (4^{n-1}/2-4)-2(2n-4)>0$, the claim holds. Note that $L_3+\{(d^+,c),(d^+,s)\}$ is a linear forest and $|E(L_3+\{(d^+,c),(d^+,s)\})\cup F_3|\leq 4\leq 2n-4$. By the induction hypothesis, $B^3-F_3$ has a H-path $P[z^+,w]$ passing through $L_3+\{(d^+,c),(d^+,s)\}$. Exactly one of $c$ and $t$, say $c$, lies on the segment of $P[z^+,w]$ between $d^+$ and $z^+$. Note that $c^+$ or $c^-$, say $c^+$, is incident with none of $E(L_0)$.

No matter which case above, by the induction hypothesis, $B^0-F_0$ has a H-path $P[a,c^+]$ passing through $L_0$. Let $y$ be the neighbor of $w^+$ on the segment of $P[a,c^+]$ between $w^+$ and $c^+$. By Theorem \ref{th-lv2014}, $B^1-\{u\}$ has a H-path $P[a^+,y^+]$. Thus, $P[a,c^+]\cup P[a^+,y^+]\cup P[v,d]\cup P[z^+,w]+\{(a,a^+),(c,c^+),(d,d^+),(w,w^+),(u,x^+),$ $(y,y^+),(z,z^+)\}-\{(w^+,y),(x^+,z),(d^+,c)\}$ is a H-path of $BH_{n}-F$ passing through $L$.

If $x=u$ and $L_2$ has a maximal path $P[x^+,r]$ with $r\neq x^+$, therefore, $|E(L_3)|\leq 1$ and $v\neq r$.

Suppose first that $r\in Y$. Then $E(L_3)\cup F_3=\emptyset$. Let $(x^+,z)\in E(P[x^+,r])$. Note that $\{v,z\}$ is compatible to $L_2-(x^+,z)$. By Theorem \ref{th-yang2019}, $B^2-F_2$ has a H-path $P[v,z]$ passing through $L_2-(x^+,z)$. Let $s,t$ be two distinct neighbors of $x^+$ on $P[v,z]$. Exactly one of $s$ and $t$, say $s$, lies on the segment of $P[v,z]$ between $x^+$ and $v$. By Lemma \ref{le-9}, there are two neighbors $c$ and $h$ of $s^+$ such that $c^+$ or $c^-$ and $h^+$ or $h^-$ are incident with none of $E(L_0)$. We claim that there is a $g\in V_3\cap X\setminus \{c,s\}$ such that $g$ and $g^{\pm}$ are incident with none of $E(L_3)$ and $E(L_0)$, respectively. The reason is follows. There are $|V_3\cap X\setminus \{c,h\}|-|E(L_3)|=4^{n-1}/2-4$ candidates of $g$. Since $E(L_0)$ has at most $|E(L_0)|$ even end vertices, each of which fails at most two candidates of such $g$. Since $|V_3\cap X\setminus \{c,h\}|-|E(L_3)|-2|E(L_{0})|\geq (4^{n-1}/2-4)-2(2n-4)>0$, the claim holds. Note that $L_3+\{(s^+,c),(s^+,h)\}$ is a linear forest and $|E(L_3+\{(s^+,c),(s^+,h)\})\cup F_3|\leq 4\leq 2n-4$. By the induction hypothesis, $B^3-F_3$ has a H-path $P[t^+,g]$ passing through $L_3+\{(s^+,c),(s^+,h)\}$. Exactly one of $c$ and $h$, say $c$, lies on the segment of $P[t^+,g]$ between $t^+$ and $s^+$. Note that $c^+$ or $c^-$, say $c^+$, is incident with none of $E(L_0)$. By the induction hypothesis, $B^0-F_0$ has a H-path $P[a,c^+]$ passing through $L_0$. Let $y$ be the neighbor of $g^+$ on the segment of $P[a,c^+]$ between $g^+$ and $c^+$. By Theorem \ref{th-lv2014}, $B^1-\{u\}$ has a H-path $P[a^+,y^+]$. Thus, $P[a,c^+]\cup P[a^+,y^+]\cup P[v,z]\cup P[t^+,g]+\{(x^+,z),(a,a^+),(c,c^+),(g,g^+),(s,s^+),(t,t^+),(u,x^+),(y,y^+)\}-\{(g^+,y),(x^+,s),$ $(x^+,t),(s^+,c)\}$ is a H-path of $BH_{n}-F$ passing through $L$.

Suppose now that $r\in X$. By the induction hypothesis, $B^2-F_2$ has a H-path $P[v,r]$ passing through $L_2$. Since $|E(L_3)|\leq 1$, $r^+$ or $r^-$, say $r^+$, is incident with none of $E(L_3)$. Let $(x^+,z)\in E(P[v,r])\setminus E(L_2)$. For $n=3$, $|E(L_0)|\leq 2$, $E(L_3)\cup F_3=\emptyset$. By Lemma \ref{le-2}, there is a $t\in V_3\cap X$ such that $t$ and $t^{\pm}$ are incident with none of $E(L_3)$ and $E(L_0)$, respectively. By Theorem \ref{th-xu2007}, $B^3$ has a H-path $P[z^+,t]$. Let $c$ be the neighbor of $r^+$ on the segment of $P[z^+,t]$ between $r^+$ and $z^+$. Since $|E(L_0)|\leq 2$, $c^+$ or $c^-$, say $c^+$, is not an internal vertex of $E(L_0)$. For $n\geq 4$, By Lemma \ref{le-9}, there are two neighbors $c$ and $s$ of $r^+$ such that $c^+$ or $c^-$ and $s^+$ or $s^-$ are incident with none of $E(L_0)$. We claim that there is an $t\in V_3\cap X\setminus \{c,s\}$ such that $t$ and $t^{\pm}$ are incident with none of $E(L_3)$ and $E(L_0)$, respectively. The reason is follows. There are $|V_3\cap X\setminus \{c,s\}|-|E(L_3)|=4^{n-1}/2-4$ candidates of $t$. Since $E(L_0)$ has at most $|E(L_0)|$ even end vertices, each of which fails at most two candidates of such $t$. Since $|V_3\cap X\setminus \{c,s\}|-|E(L_3)|-2|E(L_{0})|\geq (4^{n-1}/2-4)-2(2n-4)>0$, the claim holds. Note that $L_3+\{(r^+,c),(r^+,s)\}$ is a linear forest and $|E(L_3+\{(r^+,c),(r^+,s)\})\cup F_3|\leq 4\leq 2n-4$. By the induction hypothesis, $B^3-F_3$ has a H-path $P[z^+,t]$ passing through $L_3+\{(r^+,c),(r^+,s)\}$. Exactly one of $c$ and $t$, say $c$, lies on the segment of $P[z^+,t]$ between $r^+$ and $z^+$. Note that $c^+$ or $c^-$, say $c^+$, is incident with none of $E(L_0)$. By the induction hypothesis, $B^0-F_0$ has a H-path $P[a,c^+]$ passing through $L_0$. Let $y$ be the neighbor of $t^+$ on the segment of $P[a,c^+]$ between $t^+$ and $c^+$. By Theorem \ref{th-lv2014}, $B^1-\{u\}$ has a H-path $P[a^+,y^+]$. Thus, $P[a,c^+]\cup P[a^+,y^+]\cup P[v,r]\cup P[z^+,t]+\{(a,a^+),(c,c^+),(r,r^+),(t,t^+),(u,x^+),(y,y^+),(z,z^+)\}-\{(t^+,y),(x^+,z),(r^+,c)\}$ is a H-path of $BH_{n}-F$ passing through $L$.

{\it Case 2.2.1.2.} $E(L_{2})\cup F_{2}=\emptyset$.

In this case, $|E(L_{m})\cup F_{m}|\leq 2$ for each $m\in \{1,3\}$. Let $d\in V_3\cap Y$ such that $d$ is incident with none of $E(L_3)$.

If $x^+\neq v$, by Lemma \ref{le-2}, there is an $a\in V_1\cap Y$ such that $a$ and $a^{\pm}$ are incident with none of $E(L_1)$ and $E(L_0)$, respectively. By Lemma \ref{le-11},
there is a $z\in N_B^1(x)-\{a\}$ such that $(x,z)\notin E(L_1)$, and $z^+$ or $z^-$, say $z^+$, is incident with none of $E(L_0)$. Note that $L_1+(x,z)$ is a linear forest and $\{u,a\}$ is compatible to $L_1+(x,z)$. For $n=3$, $|E(L_1+(x,z))\cup F_1|\leq 2$; and for $n\geq 4$, $|E(L_1+(x,z))\cup F_1|\leq 2n-4$. By the induction hypothesis, $B^0-F_0$, $B^1-F_1$ have H-paths $P[a^+,b]$, $P[u,a]$ passing through $L_0$ and $L_1+(x,z)$, respectively. Let $c$ be the neighbor of $z^+$ on the segment of $P[a^+,b]$ between $a^+$ and $z^+$. Since $|E(L_3)|\leq 2$, $c^+$ or $c^-$, say $c^+$, is not an internal vertex of $L_3$. By the induction hypothesis, $B^3-F_3$ has a H-path $P[c^+,d]$ passing through $L_3$. Let $y$ be the neighbor of $b^+$ on the segment of $P[c^+,d]$ between $c^+$ and $b^+$. By Theorem \ref{th-cheng2014}, there exist two vertex-disjoint paths $P[x^+,d^+]$ and $P[y^+,v]$ in $B^2$ such that each vertex of $B^2$ lies on one of the two paths. Thus, $P[a^+,b]\cup P[u,a]\cup P[x^+,d^+]\cup P[y^+,v]\cup P[c^+,d]+\{(a,a^+),(b,b^+),(c,c^+),(d,d^+),(x,x^+),(y,y^+),(z,z^+)\}-\{(z^+,c),(x,z),(b^+,y)\}$ is a H-path of $BH_{n}-F$ passing through $L$.

If $x^+=v$, in this case, $u\neq x$, $|E(L_0)|\leq 2$ and $|E(L_m)|\leq 1$ for $m\in \{1,3\}$.

Suppose first that $x$ is incident with none of $E(L_1)$. For $n=3$, $|E(L_0)|\leq 2$ and $|E(L_m)|\leq 1$ for $m\in \{1,3\}$. By Lemma \ref{le-2}, there is an $a\in V_1\cap Y$ such that $a$ and $a^{\pm}$ are incident with none of $E(L_1)$ and $E(L_0)$, respectively. By Theorem \ref{th-yang2019}, $B^1-F_1$ has a H-path $P[u,a]$ passing through $L_1$. Let $z$ be the neighbor of $x$ on the segment of $P[u,a]$ between $x$ and $u$. Since $|E(L_0)|\leq 2$, $z^+$ or $z^-$, say $z^+$, is not an internal vertex of $L_0$. For $n\geq 4$, $|E(L_1)|\leq E(L_1)\cup F_1\leq 2\leq 2n-6$. By Lemma \ref{le-9}, there are two neighbors $z$ and $s$ of $x$ such that $z^+$ or $z^-$ and $s^+$ or $s^-$ are incident with none of $E(L_0)$. We claim that there is an $a\in V_1\cap Y\setminus \{z,s\}$ such that $a$ and $a^{\pm}$ are incident with none of $E(L_1)$ and $E(L_0)$, respectively. The reason is follows. There are $|V_1\cap Y\setminus \{z,s\}|-|E(L_1)|=4^{n-1}/2-3$ candidates of $a$. Since $E(L_0)$ has at most $|E(L_0)|$ even end vertices, each of which fails at most two candidates of such $a$. Since $|V_1\cap Y\setminus \{z,s\}|-|E(L_1)|-2|E(L_{0})|\geq (4^{n-1}/2-3)-2(2n-4)>0$, the claim holds. Note that $L_1+\{(x,z),(x,s)\}$ is a linear forest and $|E(L_1+\{(x,z),(x,s)\})\cup F_1|\leq 4<2n-4$. By the induction hypothesis, $B^1-F_1$ has a H-path $P[u,a]$ passing through $L_1+\{(x,z),(x,s)\}$. Exactly one of $z$ and $s$, say $z$, lies on the segment of $P[u,a]$ between $u$ and $x$. Note that $z^+$ or $z^-$, say $z^+$, is incident with none of $E(L_0)$. No matter which cases above, by the induction hypothesis, $B^0-F_0$ has a H-path $P[z^+,b]$ passing through $L_0$. Let $c$ be the neighbor of $a^+$ on the segment of $P[z^+,b]$ between $a^+$ and $z^+$. Since $|E(L_3)|\leq 2$, $c^+$ or $c^-$, say $c^+$, is not an internal vertex of $L_3$. By the induction hypothesis, $B^3-F_3$ has a H-path $P[c^+,d]$ passing through $L_3$. Let $y$ be the neighbor of $b^+$ on the segment of $P[c^+,d]$ between $b^+$ and $c^+$. By Theorem \ref{th-lv2014}, $B^2-\{v\}$ has a H-path $P[d^+,y^+]$. Thus, $P[z^+,b]\cup P[u,a]\cup P[d^+,y^+]\cup P[c^+,d]+\{(a,a^+),(b,b^+),(c,c^+),(d,d^+),(x,v),(y,y^+),(z,z^+)\}-\{(a^+,c),(x,z),(b^+,y)\}$ is a H-path of $BH_{n}-F$ passing through $L$.

Suppose second that $L_1$ has a maximal path $P[x,r]$ with $r\neq x$ and $n=3$. Then $r\in Y$. There is a $z\in N_B^1(x)\setminus \{r\}$, such that $z$ is not the shadow vertex of $r$. Since $|E(L_0)|\leq |E(L_0)\cup F_0|\leq 2$, there is at least one of $\{r^+,r^-,z^+,z^-\}$, say $z^+$, incident with none of $L_0$ and $r^+$ or $r^-$, say $r^+$, not an internal vertex of $L_0$. By the induction hypothesis, $B^0-F_0$, $B^1-F_1$ have H-path $P[r^+,b]$, $P[u,r]$ passing through $L_0$ and $L_1+(x,z)$, respectively. Let $c$ be the neighbor of $z^+$ on the segment of $P[r^+,b]$ between $z^+$ and $r^+$. By Theorem \ref{th-yang2019}, $B^3-F_3$ has a H-path $P[c^+,d]$ passing through $L_3$. Let $y$ be the neighbor of $b^+$ on the segment of $P[c^+,d]$ between $b^+$ and $c^+$. By Theorem \ref{th-lv2014}, $B^2-\{v\}$ has a H-path $P[d^+,y^+]$. Thus, $P[r^+,b]\cup P[u,r]\cup P[d^+,y^+]\cup P[c^+,d]+\{(b,b^+),(c,c^+),(d,d^+),(r,r^+),(x,v),(y,y^+),(z,z^+)\}-\{(z^+,c),(x,z),(b^+,y)\}$ is a H-path of $BH_{n}-F$ passing through $L$.

Suppose now that $L_1$ has a maximum path $P[x,r]$ with $r\neq x$ and $n\geq 4$. In this case, $r\neq u$ and $|E(L_3)|\leq 1$. Let $(x,h)\in E(P[x,r])$. By Lemma \ref{le-10}, there are two neighbors $z$ and $s$ of $x$ such that $h\notin \{z,s\}$, $z$ is not the shadow vertex of $s$, $z^+$ or $z^-$, say $z^+$, is incident with none of $E(L_0)$, and $s^+$ or $s^-$, say $s^+$, is not an internal vertex of $L_0$. Note that $\{u,h\}$ is compatible to $L_1+\{(x,z),(x,s)\}-(x,h)$ and $|E(L_1+\{(x,z),(x,s)\}-(x,h))\cup F_1|\leq 2<2n-4$. By the induction hypothesis, $B^1-F_1$ has a H-path $P[u,h]$ passing through $L_1+\{(x,z),(x,s)\}-(x,h)$. By the induction hypothesis, $B^0-F_0$ has a H-path $P[s^+,b]$ passing through $L_0$. Let $c$ be the neighbor of $z^+$ on the segment of $P[s^+,b]$ between $z^+$ and $s^+$. By Theorem \ref{th-yang2019}, $B^3-F_3$ has a H-path $P[c^+,d]$ passing through $L_3$. Let $y$ be the neighbor of $b^+$ on the segment of $P[c^+,d]$ between $c^+$ and $b^+$. By Theorem \ref{th-lv2014}, $B^2-\{v\}$ has a H-path $P[d^+,y^+]$. Thus, $P[s^+,b]\cup P[u,h]\cup P[d^+,y^+]\cup P[c^+,d]+\{(x,h),(b,b^+),(c,c^+),(d,d^+),(s,s^+),(x,v),(y,y^+),(z,z^+)\}-\{(z^+,c),(x,z),$ $(x,s),(b^+,y)\}$ is a H-path of $BH_{n}-F$ passing through $L$.


{\it Case 2.2.1.3.} $E(L_{3})\cup F_{3}=\emptyset$.

In this scenario, $|E(L_{m})\cup F_{m}|\leq 2$ for $m\in \{1,2\}$. The proofs for the cases that $|E(L_{1})\cup F_{1}|=2$ (resp. $|E(L_{2})\cup F_{2}|=2$) is similarly to the case that $E(L_{2})\cup F_{2}=\emptyset$ (resp. $E(L_{1})\cup F_{1}=\emptyset$). We here only consider the case that $|E(L_{1})\cup F_{1}|\leq 1$ and $|E(L_{2})\cup F_{2}|\leq 1$. By Lemma \ref{le-2}, there is an $a\in V_1\cap Y$ such that $a$ and $a^{\pm}$ are incident with none of $E(L_1)$ and $E(L_0)$, respectively. Let $d\in V_2\cap X$ such that $d$ is incident with none of $E(L_2)$. 


If $x^+\neq v$, by Lemma \ref{le-8}, there is a neighbor $z$ of $x$ such that $L_1+(x,z)$ is a linear forest and $z^+$ or $z^-$, say $z^+$, is not an internal vertex of $L_0$. By the induction hypothesis, $B^1-F_1$ has a H-path $P[u,a]$ passing through $L_1+(x,z)$. Let $g=a^-$, if $z=a$; and $g=a^+$, otherwise. By the induction hypothesis, $B^0-F_0$ has a H-path $P[z^+,b]$ passing through $L_0$. Let $c$ be the neighbor of $g$ on the segment of $P[z^+,b]$ between $g$ and $z^+$.

Suppose first that $x^+$ is incident with none of $E(L_2)$. By Theorem \ref{th-yang2019}, $B^2-F_2$ has a H-path $P[v,d]$ passing through $L_2$. Let $y$ be the neighbor of $x^+$ on the segment of $P[v,d]$ between $x^+$ and $v$. By Theorem \ref{th-cheng2014}, there exist two vertex-disjoint paths $P[b^+,y^+]$ and $P[c^+,d^+]$ in $B^3$ such that each vertex of $B^3$ lies on one of the two paths. Thus, $P[z^+,b]\cup P[u,a]\cup P[v,d]\cup P[b^+,y^+]\cup P[c^+,d^+]+\{(a,g),(b,b^+),(c,c^+),(d,d^+),(x,x^+),(y,y^+),$ $(z,z^+)\}-\{(g,c),$ $(x,z),(x^+,y)\}$ is a H-path of $BH_{n}-F$ passing through $L$.

Suppose now that $x^+$ is incident with an edge of $E(L_2)$. In this scenario, let $(x^+,r)$ be the edge of $L_2$. By Theorem \ref{th-yang2019}, $B^2-F_2$ has a H-path $P[v,r]$ passing through $L_2$. Let $y$ be the neighbor of $x^+$ on $P[v,r]$ such that $y\neq r$. By Theorem \ref{th-cheng2014}, there exist two vertex-disjoint paths $P[b^+,y^+]$ and $P[c^+,r^+]$ in $B^3$ such that each vertex of $B^3$ lies on one of the two paths. Thus, $P[z^+,b]\cup P[u,a]\cup P[v,r]\cup P[b^+,y^+]\cup P[c^+,r^+]+\{(a,g),(b,b^+),(c,c^+),(r,r^+),(x,x^+),(y,y^+),(z,z^+)\}-\{(g,c),(x,z),(x^+,y)\}$ is a H-path of $BH_{n}-F$ passing through $L$.

If $x^+=v$, $u\neq x$.

Suppose first that $x$ is incident with none of $E(L_1)$ and $n=3$. By Theorem \ref{th-yang2019}, $B^1-F_1$ has a H-path $P[u,a]$ passing through $L_1$. Let $z$ be the neighbor of $x$ on the segment of $P[u,a]$ between $x$ and $u$. Since $|E(L_0)|\leq 2$, $z^+$ or $z^-$, say $z^+$, is not an internal vertex of $L_0$. By the induction hypothesis, $B^0-F_0$ has a H-path $P[z^+,b]$ passing through $L_0$. Let $c$ be the neighbor of $a^+$ on the segment of $P[z^+,b]$ between $a^+$ and $z^+$. By Theorem \ref{th-yang2019}, $B^2-F_2$ has a H-path $P[v,d]$ passing through $L_2$. Let $(v,y)\in E(P[v,d])$. By Theorem \ref{th-cheng2014}, there exist two vertex-disjoint paths $P[b^+,y^+]$ and $P[c^+,d^+]$ in $B^3$ such that each vertex of $B^3$ lies on one of the two paths. Thus, $P[z^+,b]\cup P[u,a]\cup P[v,d]\cup P[b^+,y^+]\cup P[c^+,d^+]+\{(a,a^+),(b,b^+),(c,c^+),(d,d^+),(x,v),(y,y^+),(z,z^+)\}-\{(a^+,c),(x,z),(v,y)\}$ is a H-path of $BH_{n}-F$ passing through $L$.

Suppose second that $x$ is incident with none of $E(L_1)$ and $n\geq 4$. By Lemma \ref{le-9}, there are two neighbors $z$ and $s$ of $x$ such that $z^+$ or $z^-$ and $s^+$ or $s^-$ are incident with none of $E(L_0)$. We claim that there is an $t\in V_1\cap Y\setminus \{z,s\}$ such that $t$ and $t^{\pm}$ are incident with none of $E(L_1)$ and $E(L_0)$, respectively. The reason is follows. There are $|V_1\cap Y\setminus \{z,s\}|-|E(L_1)|=4^{n-1}/2-3$ candidates of $t$. Since $E(L_0)$ has at most $|E(L_0)|$ even end vertices, each of which fails at most two candidates of such $t$. Since $|V_1\cap Y\setminus \{z,s\}|-|E(L_1)|-2|E(L_{0})|\geq (4^{n-1}/2-3)-2(2n-4)>0$, the claim holds. Note that $L_1+\{(x,z),(x,s)\}$ is a linear forest and $|E(L_1+\{(x,z),(x,s)\})\cup F_1|\leq 4<2n-4$. By the induction hypothesis, $B^1-F_1$ has a H-path $P[u,t]$ passing through $L_1+\{(x,z),(x,s)\}$. Exactly one of $z$ and $s$, say $z$, lies on the segment of $P[u,t]$ between $u$ and $x$. Note that $z^+$ or $z^-$, say $z^+$, is incident with none of $E(L_0)$. By the induction hypothesis, $B^0-F_0$ has a H-path $P[z^+,b]$ passing through $L_0$. Let $c$ be the neighbor of $t^+$ on the segment of $P[z^+,b]$ between $t^+$ and $z^+$. By Theorem \ref{th-yang2019}, $B^2-F_2$ has a H-path $P[v,d]$ passing through $L_2$. Let $(v,y)\in E(P[v,d])$. By Theorem \ref{th-cheng2014}, there exist two vertex-disjoint paths $P[b^+,y^+]$ and $P[c^+,d^+]$ in $B^3$ such that each vertex of $B^3$ lies on one of the two paths. Thus, $P[z^+,b]\cup P[u,t]\cup P[v,d]\cup P[b^+,y^+]\cup P[c^+,d^+]+\{(b,b^+),(c,c^+),(d,d^+),(t,t^+),(x,v),(y,y^+),(z,z^+)\}-\{(t^+,c),(x,z),(v,y)\}$ is a H-path of $BH_{n}-F$ passing through $L$.

Suppose third that $L_1$ has a maximum path $P[x,r]$. Since $|E(L_1)|\leq 1$, $r\in Y$. For $n=3$. In this case, $|E(L_0)|\leq 2$. Let $z\in N_B^1(x)\setminus \{r\}$ such that $z$ is not the shadow vertex of $r$. Thus, there is at least one of $\{r^+,r^-,z^+,z^-\}$, say $r^+$, incident with none of $E(L_0)$, $z^+$ or $z^+$, say $z^+$, not an internal vertex of $L_0$. For $n\geq 4$. By Lemma \ref{le-10}, there are two neighbors $z$ and $s$ of $x$ such that $r\notin \{z,s\}$, $L_1+\{(x,z),(x,s)\}-(x,r)$ is a linear forest, $z^+$ or $z^-$, say $z^+$, is incident with none of $E(L_0)$ and $s^+$ or $s^-$, say $s^+$, is not an internal vertex of $L_0$. Note that $\{u,r\}$ is compatible $L_1+\{(x,z),(x,s)\}-(x,r)$.
By the induction hypothesis, $B^0-F_0$, $B^1-F_1$ have H-paths $P[z^+,b]$, $P[u,r]$ passing through $L_0$ and $L_1$, respectively. Let $c$ be the neighbor of $r^+$ on the segment of $P[z^+,b]$ between $r^+$ and $z^+$. By Theorem \ref{th-yang2019}, $B^2-F_2$ has a H-path $P[v,d]$ passing through $L_2$. Let $(v,y)\in E(P[v,d])$. By Theorem \ref{th-cheng2014}, there exist two vertex-disjoint paths $P[b^+,y^+]$ and $P[c^+,d^+]$ in $B^3$ such that each vertex of $B^3$ lies on one of the two paths. Thus, $P[z^+,b]\cup P[u,r]\cup P[v,d]\cup P[b^+,y^+]\cup P[c^+,d^+]+\{(b,b^+),(c,c^+),(d,d^+),(r,r^+),(x,v),(y,y^+),(z,z^+)\}-\{(r^+,c),(x,z),(v,y)\}$ is a H-path of $BH_{n}-F$ passing through $L$.

{\it Case 2.2.2.}  $i=1$, $j=3$.

By Lemma \ref{le-2}, there is a $a\in V_1\cap Y$ such that $a$ and $a^{\pm}$ are incident with none of $E(L_1)$ and $E(L_0)$, respectively. By Lemma \ref{le-11}, there is a $z\in N_{B^1}(x)-\{a\}$ such that $L_1+(x,z)$ is a linear forest, $(x,z)\notin E(L_1)$, and $z^{+}$ or $z^-$, say $z^+$, is incident with none of $E(L_0)$. Note that $\{u,a\}$ is compatible to $L_1+(x,z)$. For $n=3$, $|E(L_1+(x,z)))\cup F_1|\leq 2$; and $|E(L_1+(x,z)))\cup F_1|\leq 2n-4$, otherwise. By the induction hypothesis, $B^0-F_0$, $B^1-F_1$ have H-paths $P[a^+,b]$, $P[u,a]$ passing through $L_0$ and $L_1+(x,z)$, respectively. Let $c$ be the neighbor of $z^+$ on the segment of $P[g,b]$ between $z^+$ and $g$.


If $E(L_{3})\cup F_{3}\neq \emptyset$, then $|E(L_{m})\cup F_{m}|\leq 1$ for each $m\in \{1,2\}$. Since $|E(L_3)|\leq 2$, $c^+$ or $c^-$, say $c^+$, is not an internal vertex of $L_3$. By the induction hypothesis, $B^3-F_3$ has a H-path $P[b^+,v]$ passing through $L_3$. Let $y$ be the neighbor of $c^+$ on the segment of $P[b^+,v]$ between $c^+$ and $b^+$. Since $|E(L_2)|\leq 1$, $y^+$ or $y^-$, say $y^+$, is incident with none of $E(L_2)$. By Theorem \ref{th-yang2019}, $B^2-F_2$ has a H-path $P[x^+,y^+]$ passing through $L_2$. Thus, $P[g,b]\cup P[u,a]\cup P[x^+,y^+]\cup P[b^+,v]+\{(a,g),(b,b^+),(c,c^+),(x,x^+),(y,y^+),(z,z^+)\}-\{(z^+,c),(x,z),(c^+,y)\}$ is a H-path of $BH_{n}-F$ passing through $L$.

If $E(L_{3})\cup F_{3}=\emptyset$, then $|E(L_{m})\cup F_{m}|\leq 2$ for each $m\in \{1,2\}$. Let $y\in V_2\cap X$ such that $y$ is incident with none of $E(L_2)$. By the induction hypothesis, $B^2-F_2$ has a H-path $P[x^+,y]$ passing through $L_2$. There is a neighbor of $y$ in $B^3$, say $y^+$, being not $v$. By Theorem \ref{th-cheng2014}, there exist two vertex-disjoint paths $P[y^+,b^+]$ and $P[c^+,v]$ in $B^3$ such that each vertex of $B^3$ lies on one of the two paths. Thus, $P[z^+,b]\cup P[u,a]\cup P[x^+,y]\cup P[y^+,b^+]\cup P[c^+,v]+\{(a,g),(b,b^+),(c,c^+),(x,x^+),(y,y^+),(z,z^+)\}-\{(z^+,c),(x,z)\}$ is a H-path of $BH_{n}-F$ passing through $L$.

{\it Case 2.2.3.}  $i=2$, $j=3$.

By Lemma \ref{le-2}, there is an $a\in V_0\cap X$ such that $a$ and $a^{\pm}$ are incident with none of $E(L_0)$ and $E(L_1)$, respectively. By the induction hypothesis, $B^0-F_0$ has a H-path $P[a,b]$ passing through $L_0$.

If $x^+$ is not adjacent to $u$ or $(x^+,u)\notin E(L_2)$. In this scenario, $\{u,x^+\}$ is compatible to $L_2$. By the induction hypothesis, $B^1-F_1$, $B^2-F_2$, $B^3-F_3$ have H-paths $P[a^+,x]$, $P[x^+,u]$, $P[b^+,v]$ passing through $L_1$, $L_2$ and $L_3$, respectively.
Thus, $P[a,b]\cup P[a^+,x]\cup P[x^+,u]\cup P[b^+,v]+\{(a,a^+),(b,b^+),(x,x^+)\}$ is a H-path of $BH_{n}-F$ passing through $L$.

If $(x^+,u)\in E(L_2)$. In this case, $E(L_{m})\cup F_{m}=\emptyset$ for some $m\in \{1,3\}$.

Suppose first that $m=1$. Let $y\in V_2\cap Y$ such that $y$ is incident with none of $E(L_2)$. Since $y\neq x^+$, then $y^-\neq x$. By the induction hypothesis, $B^2-F_2$ has a H-path $P[u,y]$ passing through $L_2$. Let $c$ be the neighbor of $x^+$ on $P[u,y]$ such that $c\neq u$. Since $|E(L_3)|\leq 1$, $c^+$ or $c^-$, say $c^+$, is incident with none of $E(L_3)$. By Lemma \ref{le-2}, there is a $t\in V_3\cap X$ such that $t$ and $t^{\pm}$ are incident with none of $E(L_3)$ and $E(L_0)$, respectively. By Lemma \ref{le-11}, there is a $z\in N_B^3(c^+)-\{t\}$ such that $z^+$ or $z^-$, say $z^+$, incident with none of $E(L_0)$. Note that $\{v,t\}$ is compatible to $L_3+(c^+,z)$. By the induction hypothesis, $B^3-F_3$ has a H-path $P[v,t]$ passing through $L_3+(c^+,z)$. Let $s\in V_0\cap X$ such that $s$ is incident with none of $E(L_0)$. By the induction hypothesis, $B^0-F_0$ has a H-path $P[s,t^+]$ passing through $L_0$. Let $d$ be the neighbor of $z^+$ on the segment of $P[s,t^+]$ between $z^+$ and $t^+$. By Theorem \ref{th-cheng2014}, there exist two vertex-disjoint paths $P[s^+,y^-]$ and $P[d^+,x]$ in $B^1$ such that each vertex of $B^1$ lies on one of the two paths. Thus, $P[s,t^+]\cup P[d^+,x]\cup P[s^+,y^-]\cup P[u,y]\cup P[v,t]+\{(c,c^+),(d,d^+),(s,s^+),(t,t^+),(x,x^+),(y,y^-),(z,z^+)\}-\{(z^+,d),(x^+,c),(c^+,z)\}$ is a H-path of $BH_{n}-F$ passing through $L$.

Suppose second that $m=3$. By Lemma \ref{le-3}, there is an edge $(s,t)\in E(P[a,b])\setminus E(L_0)$ for some $s\in X$ and $t\in Y$ such that $s^+$ or $s^-$, say $s^+$, is incident with none of $E(L_1)$ and $\{s,t\}\cap \{a,b\}=\emptyset$. By Theorem \ref{th-yang2019}, $B^1-F_1$ has a H-path $P[a^+,x]$. Let $y$ be the neighbor of $s^+$ on the segment of $P[a^+,x]$ between $a^+$ and $s^+$. Note that $\{u,y^+\}$ is compatible to $L_2$. By the induction hypothesis, $B^2-F_2$ has a H-path $P[y^+,u]$ passing through $L_2$. Let $c$ be the neighbor of $x$ on $P[y^+,u]$ such that $c\neq u$. There is a neighbor of $c$ in $B^3$, say $c^+$, being not $v$. By Theorem \ref{th-cheng2014}, there exist two vertex-disjoint paths $P[b^+,c^+]$ and $P[t^+,v]$ in $B^3$ such that each vertex of $B^3$ lies on one of the two paths. Thus, $P[a,b]\cup P[a^+,x]\cup P[y^+,u]\cup P[b^+,c^+]\cup P[t^+,v]+\{(a,a^+),(b,b^+),(c,c^+),(s,s^+),(t,t^+),(x,x^+),(y,y^+)\}-\{(s,t),(s^+,y),(x^+,c)\}$ is a H-path of $BH_{n}-F$ passing through $L$.

{\it Case 3.}  $l=2$.

{\it Case 3.1.}  $i=0$.

By Lemma \ref{le-2}, there is an $a\in V_0\cap Y$ such that $a$ and $a^{\pm}$ are incident with none of $E(L_0)$ and $E(L_3)$, respectively. By the induction hypothesis, $B^0-F_0$ has a H-path $P[u,a]$ passing through $L_0$.

{\it Case 3.1.1.}  $j=1$.

By Lemma \ref{le-2}, there is a $b\in V_1\cap X$ such that $b$ (resp. $b^{+}$) is incident with none of $E(L_1)$ (resp. $E(L_2)$). By the induction hypothesis, $B^1-F_1$, $B^2-F_2$, $B^3-F_3$ have H-paths $P[v,b]$, $P[b^+,x]$, $P[a^+,x^+]$ passing through $L_1$, $L_2$ and $L_3$, respectively. Thus, $P[u,a]\cup P[v,b]\cup P[b^+,x]\cup P[a^+,x^+]+\{(a,a^+),(b,b^+),(x,x^+)\}$ is a H-path of $BH_{n}-F$ passing through $L$.

{\it Case 3.1.2.}  $j=2$.

There are $\lfloor{|E(P[u,a])|}/{2}\rfloor=\lfloor(4^{n-1}-1)/2\rfloor$ edges each of which has the form $(s,t)$ with $s\in X$ and $t\in Y$ such that $t$ lies on the segment of $P[u,a]$ between $u$ and $s$. Since $\lfloor{|E(P[u,a])|}/{2}\rfloor-|E(L_0)|\geq \lfloor(4^{n-1}-1)/2\rfloor-(2n-4)\geq 5$, there are at least such $5$ edges $(s,t)$ on $P[u,a]$ that meats above requirements and furthermore $(s,t)\notin E(L_0)$. Since $|E(L_1)|+|E(L_3)|\leq 2$, there are at most $4$ ($<5$) such edges $(s,t)$ that meats above requirements and $s^+$ or $s^-$ (resp. $t^+$ or $t^-$) is incident with some edge of $E(L_1)$ (resp. $E(L_3)$). Thus, there is an edge $(s,t)\in E(P[u,a])\setminus E(L_0)$ such that $s^{\pm}$ (resp. $t^{\pm}$) are incident with none of $E(L_1)$ (resp. $E(L_3)$).

If $E(L_{2})\cup F_{2}=\emptyset$, then $|E(L_{m})\cup F_{m}|\leq 2$ for each $m\in \{1,3\}$. Let $b\in V_1\cap X$ such that $b$ is incident with none of $E(L_1)$. By the induction hypothesis, $B^1-F_1$, $B^3-F_3$ have paths $P[s^+,b]$, $P[a^+,x^+]$ passing through $L_1$ and $L_3$, respectively. There is a neighbor of $b$ in $B^2$, say $b^+$, being not $v$. Let $y$ be the neighbor of $t^+$ on the segment of $P[a^+,x^+]$ between $t^+$ and $a^+$. Since $y\neq x^+$, $y^-\neq x$. By Theorem \ref{th-cheng2014}, there exist two vertex-disjoint paths $P[y^-,v]$ and $P[b^+,x]$ in $B^2$ such that each vertex of $B^2$ lies on one of the two paths. Thus, $P[u,a]\cup P[s^+,b]\cup P[y^-,v]\cup P[b^+,x]\cup P[a^+,x^+]+\{(a,a^+),(b,b^+),(s,s^+),(t,t^+),(x,x^+),(y,y^-)\}-\{(s,t),(t^+,y)\}$ is a H-path of $BH_{n}-F$ passing through $L$.

If $E(L_{2})\cup F_{2}\neq \emptyset$, then $|E(L_{m})\cup F_{m}|\leq 1$ for each $m\in \{1,3\}$.

Suppose first that $x$ is not adjacent to $v$ or $(x,v)\notin E(L_2)$. Then $\{v,x\}$ is compatible to $L_2$. By the induction hypothesis, $B^2-F_2$, $B^3-F_3$ have H-paths $P[v,x]$, $P[a^+,x^+]$ passing through $L_2$ and $L_3$, respectively. Let $y$ be the neighbor of $t^+$ on the segment of $P[a^+,x^+]$ between $a^+$ and $t^+$. Then $y\neq x^+$ and $y^-\neq (x^+)^-$ (i.e. $x$). Let $(y^-,b)\in E(P[v,x])\setminus E(L_2)$. By Theorem \ref{th-yang2019}, $B^1-F_1$ has a H-path $P[s^+,b^+]$ passing through $L_1$. Thus, $P[u,a]\cup P[s^+,b^+]\cup P[v,x]\cup P[a^+,x^+]\cup +\{(a,a^+),(b,b^+),(s,s^+),(t,t^+),(x,x^+),(y,y^-)\}-\{(s,t),$ $(y^-,b),(y,t^+)\}$ is a H-path of $BH_{n}-F$ passing through $L$.

Suppose now that $(x,v)\in E(L_2)$. By Theorem \ref{th-yang2019}, $B^3-F_3$ has a H-path $P[a^+,x^+]$ passing through $L_3$. Let $y$ be the neighbor of $t^+$ on the segment of $P[a^+,x^+]$ between $t^+$ and $x^+$. Let $g=y^-$, if $y\neq x^+$; and $g=y^+$, otherwise. Then $g\neq x$. By the induction hypothesis, $B^2-F_2$ has a H-path $P[g,v]$ passing through $L_2$. Let $b$ be the neighbor of $x$ on $P[g,v]$ such that $b\neq v$. By Theorem \ref{th-yang2019}, $B^1-F_1$ has a H-path $P[s^+,b^+]$ passing through $L_1$. Thus, $P[u,a]\cup P[s^+,b^+]\cup P[v,g]\cup P[a^+,x^+]\cup +\{(a,a^+),(b,b^+),(s,s^+),(t,t^+),(x,x^+),(y,g)\}-\{(s,t),(x,b),(y,$ $t^+)\}$ is a H-path of $BH_{n}-F$ passing through $L$.

{\it Case 3.1.3.}  $j=3$.

By Lemma \ref{le-2}, there is a $d\in V_3\cap X$ such that $d$ and $d^{\pm}$ are incident with none of $E(L_3)$ and $E(L_0)$, respectively. By Lemma \ref{le-11}, there is a $z\in N_{B^3}(x^+)-\{d\}$ such that $(x^+,z)\notin E(L_3)$ and $z^+$ or $z^-$, say $z^+$, is incident with none of $E(L_0)$. Note that $L_3+(x^+,z)$ is a linear forest and $\{v,d\}$ is compatible to $L_3+(x^+,z)$. For $n=3$, $|E(L_3+(x^+,z))\cup F_3|\leq 2$; and $|E(L_3+(x^+,z))\cup F_3|\leq 2n-4$, otherwise. By the induction hypothesis, $B^0-F_0$, $B^3-F_3$ have H-paths $P[u,z^+]$, $P[d,v]$ passing through $L_0$ and $L_3+(x^+,z)$, respectively. Let $y$ be the neighbor of $d^+$ on the segment of $P[u,z^+]$ between $d^+$ and $z^+$. By Lemma \ref{le-2}, there is a $w\in V_1\cap X$ such that $w$ (resp. $w^+$) is incident with none of $E(L_1)$ (resp. $E(L_2)$). By the induction hypothesis, $B^1-F_1$, $B^2-F_2$ have H-paths $P[y^+,w]$, $P[w^+,x]$ passing through $L_1$ and $L_2$, respectively. Thus, $P[u,z^+]\cup P[y^+,w]\cup P[w^+,x]\cup P[d,v]+\{(d,d^+),(w,w^+),(x,x^+),(y,y^+),(z,z^+)\}-\{(d^+,y),(x^+,z)\}$ is a H-path of $BH_{n}-F$ passing through $L$.

{{\it Case 3.2.}  $i\neq 0$.}

By Lemma \ref{le-2}, there is an $a\in V_0\cap X$ such that $a$ and $a^{\pm}$ are incident with none of $E(L_0)$ and $E(L_1)$, respectively.  

{\it Case 3.2.1.}  $i=1$, $j=2$.

By Lemma \ref{le-2}, there is a $b\in V_0\cap Y$ such that $b$ and $b^{\pm}$ are incident with none of $E(L_0)$ and $E(L_3)$, respectively.

If $x$ is not adjacent to $v$ or $(v,x)\notin E(L_2)$, $\{v,x\}$ is compatible to $L_2$. By the induction hypothesis, $B^0-F_0$ has H-path $P[a,b]$ passing through $L_0$.
By the induction hypothesis, $B^1-F_1$, $B^2-F_2$, $B^3-F_3$ have H-paths $P[a^+,u]$, $P[v,x]$, $P[x^+,b^+]$ passing through $L_1$, $L_2$ and $L_3$, respectively. Thus, $P[a,b]\cup P[a^+,u]\cup P[v,x]\cup P[x^+,b^+]+\{(a,a^+),(b,b^+),(x,x^+)\}$ is a H-path of $BH_{n}-F$ passing through $L$.

If $(v,x)\in E(L_2)$, $E(L_{m})\cup F_{m}=\emptyset$ for some $m\in \{1,3\}$. By Lemma \ref{le-2}, there is a $d\in V_2\cap X$ such that $d$ and $d^{\pm}$ are incident with none of $E(L_2)$ and $E(L_3)$, respectively. Then $d\neq x$. By the induction hypothesis, $B^2-F_2$ has a H-path $P[v,d]$ passing through $L_2$. Let $z$ be the neighbor of $x$ on $P[v,d]$ such that $z\neq v$. There is a neighbor $z$ in $B^1$, say $z^+$, being not $u$.

Suppose first that $m=1$. By Lemma \ref{le-2}, there is a $w\in V_3\cap X$ such that $w$ and $w^{\pm}$ are incident with none of $E(L_3)$ and $E(L_0)$, respectively. By Lemma \ref{le-11}, there is a $c\in N_B^3(d^+)-\{w\}$ such that $c^+$ or $c^-$, say $c^+$, is incident with none of $E(L_0)$. By the induction hypothesis, $B^0-F_0$, $B^3-F_3$ have H-paths $P[a,w^+]$, $P[x^+,w]$ passing through $L_0$ and $L_3+(d^+,c)$, respectively. Let $y$ be the neighbor of $c^+$ on the segment of $P[a,w^+]$ between $c^+$ and $w^+$. By Theorem \ref{th-cheng2014}, there exist two vertex-disjoint paths $P[a^+,z^+]$ and $P[y^+,u]$ in $B^1$ such that each vertex of $B^1$ lies on one of the two paths. Thus, $P[a,w^+]\cup P[a^+,z^+]\cup P[y^+,u]\cup P[v,d]\cup P[x^+,w]+\{(a,a^+),(c,c^+),(d,d^+),(w,w^+),(x,x^+),(y,y^+),(z,z^+)\}-\{(c^+,y),(x,z),(d^+,c)\}$ is a H-path of $BH_{n}-F$ passing through $L$.

Suppose now that $m=3$. By Lemma \ref{le-2}, there is a $h\in V_1\cap Y$ such that $h$ and $h^{\pm}$ are incident with none of $E(L_1)$ and $E(L_0)$, respectively. By Lemma \ref{le-11}, there is a $y\in N_B^1(z^+)-\{h\}$ such that $y^+$ or $y^-$, say $y^+$, is incident with none of $E(L_0)$. By the induction hypothesis, $B^0-F_0$, $B^1-F_1$ have H-paths $P[h^+,b]$, $P[u,h]$ passing through $L_0$ and $L_1+(z^+,y)$, respectively. Let $c$ be the neighbor of $y^+$ on the segment of $P[h^+,b]$ between $y^+$ and $h^+$. By Theorem \ref{th-cheng2014}, there exist two vertex-disjoint paths $P[b^+,d^+]$ and $P[x^+,c^+]$ in $B^3$ such that each vertex of $B^3$ lies on one of the two paths. Thus, $P[h^+,b]\cup P[u,h]\cup P[v,d]\cup P[b^+,d^+]\cup P[x^+,c^+]+\{(b,b^+),(c,c^+),(d,d^+),(h,h^+),(x,x^+),(y,y^+),(z,z^+)\}-\{(y^+,c),(z^+,y),(x,z)\}$ is a H-path of $BH_{n}-F$ passing through $L$.


{\it Case 3.2.2.}  $i=1$, $j=3$.

By Lemma \ref{le-2}, there is a $b\in V_3\cap X$ such that $b$ and $b^{\pm}$ are incident with none of $E(L_3)$ and $E(L_0)$, respectively. By Lemma \ref{le-11}, there is a neighbor $z\in N_{B^3}(x^+)-\{b\}$ such that $(x^+,z)\notin E(L_3)$, and $z^+$ or $z^-$, say $z^+$, is incident with none of $E(L_0)$. Then $\{v,b\}$ is compatible to $L_3+(x^+,z)$. For $n=3$, $|E(L_3+(x^+,z))\cup F_3|\leq 2$; and $|E(L_3+(x^+,z))\cup F_3|\leq 2n-4$, otherwise. By the induction hypothesis, $B^0-F_0$, $B^3-F_3$ have H-paths $P[a,z^+]$, $P[v,b]$ passing through $L_0$ and $L_3+(x^+,z)$, respectively. Let $c$ be the neighbor of $b^+$ on the segment of $P[a,z^+]$ between $b^+$ and $z^+$.

Suppose first that $E(L_1)\cup F_1=\emptyset$.
Let $y\in V_2\cap Y$ such that $y$ is incident with none of $E(L_2)$. By the induction hypothesis, $B^2-F_2$ has a H-path $P[x,y]$ passing through $L_2$. There is a neighbor of $y$ in $B^1$, say $y^+$, being not $u$. By Theorem \ref{th-cheng2014}, there exist two vertex-disjoint paths $P[a^+,u]$ and $P[y^+,c^+]$ in $B^1$ such that each vertex of $B^1$ lies on one of the two paths. Thus, $P[a,z^+]\cup P[y^+,c^+]\cup P[a^+,u]\cup P[x,y]\cup P[v,b]+\{(a,a^+),(b,b^+),(c,c^+),(x,x^+),(y,y^+),(z,z^+)\}-\{(b^+,c),(x^+,z)\}$ is a H-path of $BH_{n}-F$ passing through $L$.

Suppose second that $|E(L_1)\cup F_1|=1$. Since $|E(L_1)|\leq 1$, $c^+$ or $c^-$, say $c^+$, is incident with none of $E(L_1)$. By Theorem \ref{th-yang2019}, $B^1-F_1$ has a H-path $P[c^+,u]$ passing through $L_1$. Let $y$ be the neighbor of $a^+$ on the segment of $P[c^+,u]$ between $a^+$ and $c^+$. Since $|E(L_2)|\leq 1$, $y^+$ or $y^-$, say $y^+$, is incident with none of $E(L_2)$. By Theorem \ref{th-yang2019}, $B^2-F_2$ has a H-path $P[x,y^+]$ passing through $L_2$. Thus, $P[a,z^+]\cup P[c^+,u]\cup P[x,y^+]\cup P[v,b]+\{(a,a^+),(b,b^+),(c,c^+),(x,x^+),(y,y^+),(z,z^+)\}-\{(b^+,c),(a^+,y),(x^+,z)\}$ is a H-path of $BH_{n}-F$ passing through $L$.

Suppose now that $|E(L_1)\cup F_1|=2$. Then $E(L_m)\cup F_m=\emptyset$ for $m\in \{2,3\}$. By Lemma \ref{le-2}, there is a $d\in V_0\cap Y$ such that $d$ is incident with none of $E(L_0)$. By the induction hypothesis, $B^0-F_0$ has a H-path $P[a,d]$ passing through $L_0$.
There are $\lfloor{|E(P[a,d])|}/{2}\rfloor=\lfloor(4^{n-1}-1)/2\rfloor$ edges each of which has the form $(s,t)$ with $s\in X$ and $t\in Y$ such that $t$ lies on the segment of $P[a,d]$ between $a$ and $s$. Since $\lfloor{|E(P[a,d])|}/{2}\rfloor-|E(L_0)|\geq \lfloor(4^{n-1}-1)/2\rfloor-(2n-4)\geq 5$,
there are at least such $5$ edges $(s,t)$ on $P[u,a]$ that meats above requirements and furthermore $(s,t)\notin E(L_0)$. Since $|E(L_1)|\leq 2$, there are at most $4$ ($<5$) such edges $(s,t)$ that meats above requirements and $s^+$ or $s^-$ is incident with some edge of $E(L_1)$. Thus, there is an edge $(s,t)\in E(P[a,d])\setminus E(L_0)$ such that $s^{\pm}$ are incident with none of $E(L_1)$. By the induction hypothesis, $B^1-F_1$ has a H-path $P[s^+,u]$ passing through $L_1$. Let $y$ be the neighbor of $a^+$ on the segment of $P[s^+,u]$ between $a^+$ and $s^+$. By Theorem \ref{th-xu2007}, $B^2$ has a H-path $P[x,y^+]$.

If $x^+=v$, by Theorem \ref{th-lv2014}, $B^3-\{v\}$ has a H-path $P[t^+,d^+]$. Thus, $P[a,d]\cup P[s^+,u]\cup P[x,y^+]\cup P[t^+,d^+]+\{(a,a^+),(d,d^+),(s,s^+),(t,t^+),(x,v),(y^+,$ $y)\}-\{(s,t),(a^+,y)\}$ is a H-path of $BH_{n}-F$ passing through $L$.

If $x^+\neq v$, By Theorem \ref{th-cheng2014}, there exist two vertex-disjoint paths $P[d^+,v]$ and $P[x^+,t^+]$ in $B^3$ such that each vertex of $B^3$ lies on one of the two paths. Thus, $P[a,d]\cup P[s^+,u]\cup P[x,y^+]\cup P[d^+,v]\cup P[x^+,t^+]+\{(a,a^+),(d,d^+),(s,s^+),$ $(t,t^+),(x,x^+),(y,y^+)\}-\{(s,t),(a^+,y)\}$ is a H-path of $BH_{n}-F$ passing through $L$.

{\it Case 3.2.3.}  $i=2$, $j=3$.

{\it Case 3.2.3.1.} $|E(L_2\cup F_2)|=2$.

In this case, $E(L_m)\cup F_m=\emptyset$ for $m\in \{1,3\}$. In this case, $n\geq 4$. By Lemma \ref{le-2}, there is a $b\in B^0\cap Y$ such that $b$ is incident with none of $E(L_0)$. By the induction hypothesis, $B^0-F_0$ has a H-path $P[a,b]$ passing through $L_0$. There are $\lfloor{|E(P[a,b])|}/{2}\rfloor=\lfloor(4^{n-1}-1)/2\rfloor$ edges each of which has the form $(s,t)$ with $s\in X$ and $t\in Y$ such that $t$ lies on the segment of $P[a,b]$ between $a$ and $s$. Since $\lfloor{|E(P[a,b])|}/{2}\rfloor-|E(L_0)|\geq \lfloor(4^{n-1}-1)/2\rfloor-(2n-4)>0$,
there are at least such one edge $(s,t)$ on $P[u,a]$ that meats above requirements and furthermore $(s,t)\notin E(L_0)$.

Suppose first that $v\neq x^+$. Let $c\in V_2\cap Y$ such that $c$ is incident with none of $E(L_2)$. By the induction hypothesis, $B^2-F_2$ has a H-path $P[u,c]$ passing through $L_2$. Let $(x,y)\in E(P[u,c])\setminus E(L_2)$. Let $g=c^-$, if $y=c$; and $g=c^+$, otherwise. Then $g\neq y^+$. By Theorem \ref{th-cheng2014}, there exist two vertex-disjoint paths $P[a^+,y^+]$ and $P[s^+,g]$ (resp. $P[b^+,v]$ and $P[x^+,t^+]$) in $B^1$ (resp. $B^3$) such that each vertex of $B^1$ (resp. $B^3$) lies on one of the two paths. Thus, $P[a,b]\cup P[a^+,y^+]\cup P[s^+,g]\cup P[u,c]\cup P[x^+,t^+]\cup P[b^+,v]+\{(a,a^+),(b,b^+),(c,g),(s,s^+),(t,t^+),(x,x^+),(y,y^+)\}-\{(s,t),(x,y)\}$ is a H-path of $BH_{n}-F$ passing through $L$.

Suppose second that $v=x^+$ and $x$ is incident with none of $E(L_2)$. In this case, $u\neq x$. Let $c\in V_2\cap Y$ such that $c$ is incident with none of $E(L_2)$. By the induction hypothesis, $B^2-F_2$ has a H-path $P[u,c]$ passing through $L_2$. Let $y$ be the neighbor of $x$ on the segment of $P[u,c]$ between $x$ and $u$. By Theorem \ref{th-cheng2014}, there exist two vertex-disjoint paths $P[a^+,y^+]$ and $P[s^+,c^+]$ in $B^1$ such that each vertex of $B^1$ lies on one of the two paths. By Theorem \ref{th-lv2014}, $B^3-\{v\}$ has a H-path $P[t^+,b^+]$. Thus, $P[a,b]\cup P[a^+,y^+]\cup P[s^+,c^+]\cup P[u,c]\cup P[t^+,b^+]+\{(a,a^+),(b,b^+),(c,c^+),(s,s^+),(t,t^+),(x,v),(y,y^+)\}-\{(s,t),(x,y)\}$ is a H-path of $BH_{n}-F$ passing through $L$.

Suppose now that $v=x^+$ and $L_2$ has a maximal path $P[x,r]$ with $r\neq x$. In this case, $u\neq r$.
Let $(x,w)\in E(P[x,r])$. Recall that $n\geq 4$. By Lemma \ref{le-10}, there are two distinct vertices $y,c\in N_B^2(x)\setminus \{w\}$ such that $L_2+\{(x,y),(x,c)\}-(x,w)$ is a linear forest. Note that $\{u,w\}$ is compatible to $L_2+\{(x,y),(x,c)\}-(x,w)$ and $|E(L_2+\{(x,y),(x,c)\}-(x,w))\cup F_2|\leq 2n-4$, by the induction hypothesis, $B^2-F_2$ has a H-path $P[u,w]$ passing through $L_2+\{(x,y),(x,c)\}-(x,w)$. By Theorem \ref{th-cheng2014}, there exist two vertex-disjoint paths $P[a^+,y^+]$ and $P[s^+,c^+]$ in $B^1$ such that each vertex of $B^1$ lies on one of the two paths. By Theorem \ref{th-lv2014}, $B^3-\{v\}$ has a H-path $P[t^+,b^+]$. Thus, $P[a,b]\cup P[a^+,y^+]\cup P[s^+,c^+]\cup P[u,w]\cup P[t^+,b^+]+\{(x,w),(a,a^+),(b,b^+),(c,c^+),(s,s^+),(t,t^+),(x,v),(y,y^+)\}-\{(s,t),(x,y),(x,c)\}$ is a H-path of $BH_{n}-F$ passing through $L$.

{\it Case 3.2.3.2.} $|E(L_2\cup F_2)|=1$.

For $n=3$, $E(L_m)\cup F_m=\emptyset$ for each $m\in \{1,3\}$ is similarly to the case that $|E(L_2)\cup F_2|=2$, we can construct a H-path of $BH_n-F$ passing through $L$. It remains to consider $|E(L_2)\cup F_2|=1$, $n\geq 4$. In this case, $E(L_m)\cup F_m=\emptyset$ for some $m\in \{1,3\}$. The proofs for the cases that $m=1$ and $m=3$ are analogous. We here only consider that $m=1$. Let $c\in V_2\cap Y$ such that $c$ is incident with none of $E(L_2)$. By Theorem \ref{th-yang2019}, $B^2-F_2$ has a H-path $P[u,c]$ passing through $L_2$. Let $(x,y)\in E(P[u,c])\setminus E(L_2)$.

If $x^+\neq v$ and $x^+$ is incident with none of $E(L_3)$,
by Lemma \ref{le-9}, there are two neighbors $z$ and $d$ of $x^+$ such that $z^+$ or $z^-$ and $d^+$ or $d^-$ are incident with none of $E(L_0)$. We claim that there is an $t\in V_3\cap X\setminus \{z,d\}$ such that $t$ and $t^{\pm}$ are incident with none of $E(L_3)$ and $E(L_0)$, respectively. The reason is follows. There are $|V_3\cap X\setminus \{z,d\}|-|E(L_3)|=4^{n-1}/2-3$ candidates of $t$. Since $E(L_0)$ has at most $|E(L_0)|$ even end vertices, each of which fails at most two candidates of such $t$. Since $|V_3\cap X\setminus \{z,d\}|-|E(L_3)|-2|E(L_{0})|\geq (4^{n-1}/2-3)-2(2n-4)>0$, the claim holds. Note that $L_3+\{(x^+,z),(x^+,d)\}$ is a linear forest and $|E(L_3+\{(x^+,z),(x^+,d)\})\cup F_3|\leq 3\leq 2n-4$. By the induction hypothesis, $B^3-F_3$ has a H-path $P[v,t]$ passing through $L_3+\{(x^+,z),(x^+,d)\}$. Exactly one of $z$ and $d$, say $z$, lies on the segment of $P[v,t]$ between $x^+$ and $v$. Note that $z^+$ or $z^-$, say $z^+$, is incident with none of $E(L_0)$. By the induction hypothesis, $B^0-F_0$ has a H-path $P[a,t^+]$ passing through $L_0$. Let $s$ be the neighbor of $z^+$ on the segment of $P[a,t^+]$ between $z^+$ and $t^+$. By Theorem \ref{th-cheng2014}, there exist two vertex-disjoint paths $P[a^+,c^+]$ and $P[s^+,y^+]$ in $B^1$ such that each vertex of $B^1$ lies on one of the two paths. Thus, $P[a,t^+]\cup P[a^+,c^+]\cup P[s^+,y^+]\cup P[u,c]\cup P[v,t]+\{(a,a^+),(c,c^+),(s,s^+),(t,t^+),(x,x^+),(y,y^+),(z,z^+)\}-\{(z^+,s),(x,y),(x^+,z)\}$ is a H-path of $BH_{n}-F$ passing through $L$.

If $x^+\neq v$ and $x^+$ is incident with an edge of $E(L_3)$, let $(x^+,w)\in E(L_3)$. By Lemma \ref{le-10}, there are two distinct vertices $z,d\in N_B^3(x^+)\setminus \{w\}$ such that $L_3+\{(x^+,z),(x^+,d)\}-(x^+,w)$ is a linear forest, $z^+$ or $z^-$, say $z^+$, is incident with none of $E(L_0)$ and $d^+$ or $d^-$, say $d^+$, is not an internal vertex of $L_0$. Note that $\{v,w\}$ is compatible to $L_3+\{(x^+,z),(x^+,d)\}-(x^+,w)$. By the induction hypothesis, $B^0-F_0$, $B^3-F_3$ have H-paths $P[a,d^+]$, $P[v,w]$ passing through $L_0$ and $L_3+\{(x^+,z),(x^+,d)\}-(x^+,w)$, respectively. Let $s$ be the neighbor of $z^+$ on the segment of $P[a,d^+]$ between $z^+$ and $d^+$.

Suppose first that $z$ lies on the segment of $P[v,w]$ between $x^+$ and $v$. By Theorem \ref{th-cheng2014}, there exist two vertex-disjoint paths $P[a^+,c^+]$ and $P[s^+,y^+]$ in $B^1$ such that each vertex of $B^1$ lies on one of the two paths. Thus, $P[a,d^+]\cup P[a^+,c^+]\cup P[s^+,y^+]\cup P[u,c]\cup P[v,w]+\{(x^+,w),(a,a^+),(c,c^+),(d,d^+),(s,s^+),$ $(x,x^+),(y,y^+),(z,z^+)\}-\{(z^+,s),(x,y),(x^+,z),(x^+,d)\}$ is a H-path of $BH_{n}-F$ passing through $L$.

Suppose now that $d$ lies on the segment of $P[v,w]$ between $x^+$ and $v$. By Theorem \ref{th-cheng2014}, there exist two vertex-disjoint paths $P[a^+,y^+]$ and $P[s^+,c^+]$ in $B^1$ such that each vertex of $B^1$ lies on one of the two paths. Thus, $P[a,d^+]\cup P[a^+,y^+]\cup P[s^+,c^+]\cup P[u,c]\cup P[v,w]+\{(x^+,w),(a,a^+),(c,c^+),(d,d^+),$ $(s,s^+),(x,x^+),(y,y^+),(z,z^+)\}-\{(z^+,s),(x,y),(x^+,z),(x^+,d)\}$ is a H-path of $BH_{n}-F$ passing through $L$.

If $x^+=v$, in this case, $u\neq x$. By Lemma \ref{le-2}, there is a $t\in V_3\cap X$ such that $t$ and $t^{\pm}$ are incident with none of $E(L_3)$ and $E(L_0)$, respectively. By Lemma \ref{le-11}, there is a $z\in N_B^3(x^+)-\{t\}$ such that $z^+$ or $z^-$, say $z^+$, is incident with none of $E(L_0)$. By the induction hypothesis, $B^0-F_0$, $B^3-F_3$ have H-paths $P[a,t^+]$, $P[v,t]$ passing through $L_0$ and $L_3+(x^+,z)$, respectively. Let $s$ be the neighbor of $z^+$ on the segment of $P[a,t^+]$ between $z^+$ and $t^+$.

Suppose first that $x$ is incident with none of $E(L_2)$.
Let $y$ be the neighbor of $x$ on the segment of $P[u,c]$ between $x$ and $u$. By Theorem \ref{th-cheng2014}, there exist two vertex-disjoint paths $P[a^+,c^+]$ and $P[s^+,y^+]$ in $B^1$ such that each vertex of $B^1$ lies on one of the two paths. Thus, $P[a,t^+]\cup P[a^+,c^+]\cup P[s^+,y^+]\cup P[u,c]\cup P[v,t]+\{(a,a^+),(c,c^+),(s,s^+),(t,t^+),(x,v),(y,y^+),(z,z^+)\}-\{(z^+,s),$ $(x,y),(x^+,z)\}$ is a H-path of $BH_{n}-F$ passing through $L$.

Suppose now that $x$ is incident with an edge of $E(L_2)$.
Let $(x,r)\in E(L_2)$. By Lemma \ref{le-10}, there are two distinct vertices $c,y\in N_B^2(x)\setminus \{r\}$ such that $L_2+\{(x,c),(x,y)\}-(x,r)$ is a linear forest. By the induction hypothesis, $B^2-F_2$ has a H-path $P[u,r]$ passing through $L_2+\{(x,c),(x,y)\}-(x,r)$. By Theorem \ref{th-cheng2014}, there exist two vertex-disjoint paths $P[a^+,c^+]$ and $P[s^+,y^+]$ in $B^1$ such that each vertex of $B^1$ lies on one of the two paths. Thus, $P[a,t^+]\cup P[a^+,c^+]\cup P[s^+,y^+]\cup P[u,r]\cup P[v,t]+\{(x,r),(a,a^+),(c,c^+),(s,s^+),(t,t^+),(x,$ $v),(y,y^+),(z,z^+)\}-\{(z^+,s),(x,c),(x,y),(x^+,z)\}$ is a H-path of $BH_{n}-F$ passing through $L$.

{\it Case 3.2.3.3.} $E(L_2)\cup F_2=\emptyset$.

In this case, $|E(L_m)\cup F_m|\leq 2$ for $m\in \{1,3\}$ is similarly to the case that $E(L_m)\cup F_m=\emptyset$ for some $m\in \{1,3\}$.

{\it Case 4.}  $l=3$.

{\it Case 4.1.}  $i=0$.

By Lemma \ref{le-2}, there is an $a\in V_0\cap Y\setminus \{x^+\}$ such that $a$ and $a^{\pm}$ are incident with none of $E(L_0)$ and $E(L_3)$, respectively. Since $a\neq x^+$, $a^-\neq (x^+)^-$ (i.e. $x$).


{\it Case 4.1.1.}  $j=1$

If $\{u,x^+\}$ is compatible to $L_0$, by the induction hypothesis, $B^0-F_0$ has a H-path $P[u,x^+]$ passing through $L_0$. By Lemma \ref{le-2}, there are vertices $z\in V_1\cap X$, $y\in V_2\cap X$ such that $z$ (resp. $z^+$) is incident with none of $E(L_1)$ (resp. $E(L_2)$) and $y$ (resp. $y^+$) is incident with none of $E(L_2)$ (resp. $E(L_3)$). By the induction hypothesis, $B^1-F_1$, $B^2-F_2$, $B^3-F_3$ have H-paths $P[v,z]$, $P[z^+,y]$, $P[y^+,x]$ passing through $L_1$, $L_2$ and $L_3$, respectively. Thus, Thus, $P[u,x^+]\cup P[v,z]\cup P[z^+,y]\cup P[y^+,x]+\{(x,x^+),(y,y^+),(z,z^+)\}$ is a H-path of $BH_{n}-F$ passing through $L$.

If $L_0$ has a maximum path $P[u,x^+]$, by the induction hypothesis, $B^0-F_0$ has a H-path $P[u,a]$ passing through $L_0$. Let $(x^+,y)\in E(P[u,a])\setminus E(L_0)$.

Suppose first that $E(L_1)\cup F_1=\emptyset$. By Lemma \ref{le-2}, there is a $w\in V_3\cap Y$ such that $w$ and $w^{\pm}$ are incident with none of $E(L_3)$ and $E(L_2)$, respectively. By the induction hypothesis, $B^3-F_3$ has a H-path $P[x,w]$ passing through $L_3$. Let $b$ be the neighbor of $a^-$ on the segment of $P[x,w]$ between $a^-$ and $x$. Since $|E(L_2)|\leq 2$, $b^+$ or $b^-$, say $b^+$ is not an internal vertex of $L_2$. Let $z\in V_2\cap Y$ such that $z$ is incident with none of $E(L_2)$. By the induction hypothesis, $B^2-F_2$ has a H-path $P[b^+,z]$ passing through $L_2$. Let $c$ be the neighbor of $w^+$ on the segment of $P[b^+,z]$ between $w^+$ and $b^+$. By Theorem \ref{th-cheng2014}, there exist two vertex-disjoint paths $P[y^+,c^+]$ and $P[z^+,v]$ in $B^1$ such that each vertex of $B^1$ lies on one of the two paths. Thus, $P[u,a]\cup P[y^+,c^+]\cup P[z^+,v]\cup P[b^+,z]\cup P[x,w]+\{(a,a^-),(b,b^+),(c,c^+),(w,w^+),(x,x^+),(y,y^+),(z,z^+)\}-\{(x^+,y),(w^+,c),(a^-,b)\}$ is a H-path of $BH_{n}-F$ passing through $L$.

Suppose second that $E(L_2)\cup F_2=\emptyset$. Since $|E(L_1)|\leq 2$, $y^+$ or $y^-$, say $y^+$, is not an internal vertex of $L_1$. Let $z\in V_1\cap X$ such that $z$ is incident with none of $E(L_1)$. By the induction hypothesis, $B^1-F_1$ has a H-path $P[v,z]$ passing though $L_1$. Let $(y^+,c)\in E(P[v,z])\setminus E(L_1)$. Let $g=z^-$, if $c=z$; and $g=z^+$, otherwise. Then $g\neq c^+$. Let $w\in V_3\cap Y$ such that $w$ is incident with none of $E(L_3)$. By the induction hypothesis, $B^3-F_3$ has a H-path $P[x,w]$ passing through $L_3$. Let $b$ be the neighbor of $a^-$ on the segment of $P[x,w]$ between $a^-$ and $x$. By Theorem \ref{th-cheng2014}, there exist two vertex-disjoint paths $P[w^+,c^+]$ and $P[g,b^+]$ in $B^2$ such that each vertex of $B^2$ lies on one of the two paths. Thus, $P[u,a]\cup P[v,z]\cup P[w^+,c^+]\cup P[g,b^+] P[x,w]+\{(a,a^-),(b,b^+),(c,c^+),(w,w^+),(x,x^+),(y,y^+),(z,g)\}-\{(x^+,y),(y^+,c),(a^-,b)\}$ is a H-path of $BH_{n}-F$ passing through $L$.

Suppose now that $E(L_3)\cup F_3=\emptyset$. Since $|E(L_1)|\leq 2$, $y^+$ or $y^-$, say $y^+$, is not an internal vertex of $L_1$. By Lemma \ref{le-2}, there is a $z\in V_1\cap X$ such that $z$ and $z^{\pm}$ are incident with none of $E(L_1)$ and $E(L_2)$, respectively. By the induction hypothesis, $B^1-F_1$ has a H-path $P[v,z]$ passing though $L_1$. Let $(y^+,c)\in E(P[v,z])\setminus E(L_1)$. Since $|E(L_2)|\leq 2$, $c^+$ or $c^-$, say $c^+$, is not an internal vertex of $L_2$. Let $g=z^-$, if $c=z$; and $g=z^+$, otherwise. Let $w\in V_2\cap X$ such that $w$ is incident with none of $E(L_2)$. By the induction hypothesis, $B^2-F_2$ has a H-path $P[c^+,w]$ passing though $L_2$. Let $b$ be the neighbor of $g$ on the segment of $P[c^+,w]$ between $g$ and $c^+$. By Theorem \ref{th-cheng2014}, there exist two vertex-disjoint paths $P[a^-,b^+]$ and $P[x,w^+]$ in $B^3$ such that each vertex of $B^3$ lies on one of the two paths. Thus, $P[u,a]\cup P[v,z]\cup P[c^+,w]\cup P[a^-,b^+]\cup P[x,w^+]+\{(a,a^-),(b,b^+),(c,c^+),(w,w^+),(x,x^+),(y,y^+),(z,g)\}-\{(x^+,y),(y^+,c),(g,b)\}$ is a H-path of $BH_{n}-F$ passing through $L$.

{\it Case 4.1.2.}  $j=2$.

By the induction hypothesis, $B^0-F_0$ has a H-path $P[u,a]$ passing through $L_0$. Let $(x^+,y)\in E(P[u,a])\setminus E(L_0)$.

Suppose first that $E(L_3)\cup F_3=\emptyset$. Since $|E(L_1)|\leq 2$, $y^+$ or $y^-$, say $y^+$, is not an internal vertex of $L_1$. By Lemma \ref{le-2}, there is a $z\in V_1\cap X$ such that $z$ and $z^{\pm}$ are incident with none of $E(L_1)$ and $E(L_2)$, respectively. By the induction hypothesis, $B^1-F_1$ has a H-path $P[y^+,z]$ passing through $L_1$. There is a neighbor of $z$ in $B^2$, say $z^+$, being not $v$. Let $w\in V_2\cap X$ such that $w$ is incident with none of $E(L_2)$. By the induction hypothesis, $B^2-F_2$ has a H-path $P[v,w]$ passing through $L_2$. Let $b$ be the neighbor of $z^+$ on the segment of $P[v,w]$ between $z^+$ and $v$. By Theorem \ref{th-cheng2014}, there exist two vertex-disjoint paths $P[a^-,b^+]$ and $P[x,w^+]$ in $B^3$ such that each vertex of $B^3$ lies on one of the two paths. Thus, $P[u,a]\cup P[y^+,z]\cup P[v,w]\cup P[a^-,b^+]\cup P[x,w^+]+\{(a,a^-),(b,b^+),(w,w^+),(x,x^+),(y,y^+),(z,z^+)\}-\{(x^+,y),(z^+,b)\}$ is a H-path of $BH_{n}-F$ passing through $L$.

Suppose now that $E(L_3)\cup F_3\neq \emptyset$. In this case, $|E(L_m)\cup F_m|\leq 1$ for each $m\in \{1,2\}$. Since $|E(L_1)|\leq 1$, $y^+$ or $y^-$, say $y^+$, is incident with none of $E(L_1)$. By Lemma \ref{le-2}, there is a $w\in V_3\cap Y$ such that $w$ and $w^{\pm}$ are incident with none of $E(L_3)$ and $E(L_2)$, respectively. By the induction hypothesis, $B^3-F_3$ has a H-path $P[x,w]$ passing through $L_3$. Let $b$ be the neighbor of $a^-$ on the segment of $P[x,w]$ between $a^-$ and $x$. Since $|E(L_2)|\leq 1$, $b^+$ or $b^-$, say $b^+$, is incident with none of $E(L_2)$. By Theorem \ref{th-yang2019}, $B^2-F_2$ has a H-path $P[b^+,v]$ passing through $L_2$. Let $z$ be the neighbor of $w^+$ on the segment of $P[b^+,v]$ between $w^+$ and $b^+$. By Theorem \ref{th-yang2019}, $B^1-F_1$ has a H-path $P[y^+,z^+]$ passing through $L_1$. Thus, $P[u,a]\cup P[y^+,z^+]\cup P[b^+,v]\cup P[x,w]+\{(a,a^-),(b,b^+),(w,w^+),(x,x^+),(y,y^+),(z,z^+)\}-\{(x^+,y),(w^+,z),(a^-,b)\}$ is a H-path of $BH_{n}-F$ passing through $L$.


{\it Case 4.1.2.}  $j=3$.

Suppose first that $x$ is not adjacent to $v$ or $(x,v)\notin E(L_3)$. In this case, $\{v,x\}$ is compatible to $L_3$. By the induction hypothesis, $B^0-F_0$, $B^3-F_3$ have H-paths $P[u,a]$, $P[v,x]$ passing through $L_0$ and $L_3$, respectively. Let $(x^+,y)\in E(P[u,a])\setminus E(L_0)$ and let $b$ be the neighbor of $a^-$ on the segment of $P[x,v]$ between $a^-$ and $x$. Since $|E(L_2)|\leq 2$ (resp. $|E(L_1)|\leq 2$), $b^+$ or $b^-$ (resp. $y^+$ or $y^-$), say $b^+$ (resp. $y^+$), is not an internal vertex of $L_2$ (resp. $L_1$). By Lemma \ref{le-2}, there is a $z\in V_2\cap Y$ such that $z$ (resp. $z^+$) is incident with none of $E(L_2)$ (resp. $E(L_1)$). By the induction hypothesis, $B^1-F_1$, $B^2-F_2$ have H-paths $P[y^+,z^+]$, $P[b^+,z]$ passing through $L_1$ and $L_2$, respectively. Thus, $P[u,a]\cup P[y^+,z^+]\cup P[b^+,z]\cup P[v,x]+\{(a,a^-),(b,b^+),(x,x^+),(y,y^+),(z,z^+)\}-\{(x^+,y),(a^-,b)\}$ is a H-path of $BH_{n}-F$ passing through $L$.

Suppose now that $(x,v)\in E(L_3)$. In this scenario, $|E(L_{m})\cup F_{m}|\leq 1$ for each $m\in \{1,2\}$. Since $\{u,v\}$ is compatible to $L$, none of the paths in $L_0$ has both $u$ and $x^+$ as end vertices. Then $\{u,x^+\}$ is compatible to $L_0$. By the induction hypothesis, $B^0-F_0$, $B^3-F_3$ have H-paths $P[u,x^+]$, $P[v,a^-]$ passing through $L_0$ and $L_3$, respectively. Let $y$ be the neighbor of $a$ on the segment of $P[u,x^+]$ between $a$ and $x^+$ and let $b$ be the neighbor of $x$ on $P[v,a^-]$ such that $b\neq v$. By Lemma \ref{le-2}, there is a $z\in V_2\cap Y$ such that $z$ (resp. $z^+$) is incident with none of $E(L_2)$ (resp. $E(L_1)$). By Theorem \ref{th-yang2019}, $B^1-F_1$, $B^2-F_2$ have H-paths $P[y^+,z^+]$, $P[b^+,z]$ passing through $L_1$ and $L_2$, respectively. Thus, $P[u,x^+]\cup P[y^+,z^+]\cup P[b^+,z]\cup P[v,a^-]+\{(a,a^-),(b,b^+),(x,x^+),(y,y^+),(z,z^+)\}-\{(a,y),(x,b)\}$ is a H-path of $BH_{n}-F$ passing through $L$.

{\it Case 4.2.}  $i\neq 0$.

By Lemma \ref{le-2}, there is an $a\in V_0\cap X$ such that $a$ and $a^{\pm}$ are incident with none of $E(L_0)$ and $E(L_1)$, respectively. By the induction hypothesis, $B^0-F_0$ has a H-path $P[a,x^+]$ passing through $L_0$.

{\it Case 4.2.1.}  $i=1$, $j=2$.

By Lemma \ref{le-2}, there is a $b\in V_3\cap Y$ such that $b$ (resp. $b^+$) is incident with none of $E(L_3)$ (resp. $E(L_2)$). By the induction hypothesis, $B^1-F_1$, $B^2-F_2$, $B^3-F_3$ have H-paths $P[a^+,u]$, $P[v,b^+]$, $P[x,b]$ passing through $L_1$, $L_2$ and $L_3$, respectively. Thus, $P[a,x^+]\cup P[a^+,u]\cup P[v,b^+]\cup P[x,b]+\{(a,a^+),(b,b^+),(x,x^+)\}$ is a H-path of $BH_{n}-F$ passing through $L$.

{\it Case 4.2.2.}  $i=1$, $j=3$.

There are $\lfloor{|E(P[a,x^+])|}/{2}\rfloor=\lfloor(4^{n-1}-1)/2\rfloor$ edges each of which has the form $(s,t)$ with $s\in X$ and $t\in Y$ such that $t$ lies on the segment of $P[a,x^+]$ between $a$ and $s$. Since $\lfloor{|E(P[a,x^+])|}/{2}\rfloor-|E(L_0)|\geq \lfloor(4^{n-1}-1)/2\rfloor-(2n-4)\geq 5$, there are at least such $5$ edges $(s,t)$ on $P[a,x^+]$ that meats above requirements and furthermore $(s,t)\notin E(L_0)$. Since $|E(L_1)|+|E(L_3)|\leq 2$, there are at most $4$ ($<5$) such edges $(s,t)$ that meats above requirements and $s^+$ or $s^-$ (resp. $t^+$ or $t^-$) is incident with some edge of $E(L_1)$ (resp. $E(L_3)$). Thus, there is an edge $(s,t)\in E(P[a,x^+])\setminus E(L_0)$ such that $s^{\pm}$ (resp. $t^{\pm}$) are incident with none of $E(L_1)$ (resp. $E(L_3)$). Since $t\neq x^+$, $t^-\neq (x^+)^-$
(i.e. $x$).

Suppose first that $|E(L_2)\cup F_2|\leq 1$. By the induction hypothesis, $B^1-F_1$, $B^3-F_3$ have H-paths $P[s^+,u]$, $P[t^-,v]$ passing through $L_1$ and $L_3$, respectively. Let $z$ be the neighbor of $a^+$ on the segment of $P[s^+,u]$ between $a^+$ and $s^+$. Let $(x,b)\in E(P[t^-,v])\setminus E(L_3)$. Since $|E(L_2)|\leq 1$, $z^+$ or $z^-$, say $z^+$, is incident with none of $E(L_2)$. By Theorem \ref{th-yang2019}, $B^2-F_2$ has a H-path $P[z^+,b^+]$ passing through $L_2$. Thus, $P[a,x^+]\cup P[s^+,u]\cup P[z^+,b^+]\cup P[t^-,v]+\{(a,a^+),(b,b^+),(s,s^+),(t,t^-),(x,x^+),(z,z^+)\}-\{(s,t),(a^+,z),(x,b)\}$ is a H-path of $BH_{n}-F$ passing through $L$.

Suppose now that $|E(L_2)\cup F_2|=2$. In this case, $E(L_m)\cup F_m=\emptyset$ for $m\in \{1,3\}$. Let $z\in V_2\cap Y$ such that $z$ is incident with none of $E(L_2)$ and let $b\in V_2\cap X$ such that $b$ is incident with none of $E(L_2)$. By the induction hypothesis, $B^2-F_2$ has a H-path $P[b,z]$ passing through $L_2$. There is a neighbor of $z$ (resp. $b$) in $B^1$ (resp. $B^3$), say $z^+$ (resp. $b^+$), being not $u$ (resp. $v$). By Theorem \ref{th-cheng2014}, there exist two vertex-disjoint paths $P[a^+,u]$ and $P[s^+,z^+]$ (resp. $P[x,v]$ and $P[t^-,b^+]$) in $B^1$ (resp. $B^3$) such that each vertex of $B^1$ (resp. $B^3$) lies on one of the two paths. Thus, $P[a,x^+]\cup P[a^+,u]\cup P[s^+,z^+]\cup P[b,z]\cup P[x,v]\cup P[t^-,b^+]+\{(a,a^+),(b,b^+),(s,s^+),(t,t^-),(x,x^+),(z,z^+)\}-(s,t)$ is a H-path of $BH_{n}-F$ passing through $L$.

{\it Case 4.2.3.}  $i=2$, $j=3$.

If $x$ is not adjacent to $v$ or $(x,v)\notin E(L_3)$, by Lemma \ref{le-2}, there is a $b\in V_1\cap X$ such that $b$ (resp. $b^+$) is incident with none of $E(L_1)$ (resp. $E(L_2)$). By the induction hypothesis, $B^1-F_1$, $B^2-F_2$, $B^3-F_3$ have H-paths $P[a^+,b]$, $P[b^+,u]$, $P[x,v]$ passing through $L_1$, $L_2$ and $L_3$, respectively. Thus, $P[a,x^+]\cup P[a^+,b]\cup P[b^+,u]\cup P[x,v]+\{(a,a^+),(b,b^+),(x,x^+)\}$ is a H-path of $BH_{n}-F$ passing through $L$.

If $(x,v)\in E(L_3)$, then $E(L_{m})\cup F_{m}=\emptyset$ for some $m\in \{1,2\}$. According to the Case 4.2.2, there is an edge $(s,t)\in E(P[a,x^+])\setminus E(L_0)$ for $s\in X$ and $t\in Y$ such that $t$ lies on the segment of $P[a,x^+]$ between $a$ and $s$, and $s^{\pm}$ (resp. $t^{\pm}$) are incident with none of $E(L_1)$ (resp. $E(L_3)$). Since $t\neq x^+$, $t^-\neq (x^+)^-$ (i.e. $x$). By the induction hypothesis, $B^3-F_3$ has a H-path $P[t^-,v]$ passing through $L_3$. Let $b$ be the neighbor of $x$ on $P[t^-,v]$ such that $b\neq v$.

Suppose first that $m=1$. Since $|E(L_2)|\leq 1$, $b^+$ or $b^-$, say $b^+$, is incident with none of $E(L_2)$. Let $z\in V_2\cap Y$. By Theorem \ref{th-yang2019}, $B^2-F_2$ has a H-path $P[u,z]$ passing through $L_2$. Let $c$ be the neighbor of $b^+$ on the segment of $P[u,z]$ between $b^+$ and $u$, if $u\neq b^+$; and let $c$ be the neighbor of $b^+$ on the segment of $P[u,z]$ between $b^+$ and $z$, otherwise. Then $c\neq z$. By Theorem \ref{th-cheng2014}, there exist two vertex-disjoint paths $P[a^+,c^+]$ and $P[s^+,z^+]$ in $B^1$ such that each vertex of $B^1$ lies on one of the two paths. Thus, $P[a,x^+]\cup P[a^+,c^+]\cup P[s^+,z^+]\cup P[u,z]\cup P[t^-,v]+\{(a,a^+),(b,b^+),(c,c^+),(s,s^+),(t,t^-),(x,x^+),(z^+,$ $z)\}-\{(s,t),(b^+,c),(x,b)\}$ is a H-path of $BH_{n}-F$ passing through $L$.

Suppose now that $m=2$. There is a neighbor of $b$ in $B^2$, say $b^+$, being not $u$. Let $z\in V_1\cap X$. By Theorem \ref{th-yang2019}, $B^1-F_1$ has a H-path $P[a^+,z]$ passing through $L_1$. Let $c$ be the neighbor of $s^+$ on the segment of $P[a^+,z]$ between $s^+$ and $a^+$. By Theorem \ref{th-cheng2014}, there exist two vertex-disjoint paths $P[c^+,u]$ and $P[z^+,b^+]$ in $B^2$ such that each vertex of $B^2$ lies on one of the two paths. Thus, $P[a,x^+]\cup P[a^+,z]\cup P[c^+,u]\cup P[z^+,b^+]\cup P[t^-,v]+\{(a,a^+),(b,b^+),(c,c^+),(s,s^+),(t,t^-),(x,x^+),(z,z^+)\}-\{(s,t),(s^+,c),(x,b)\}$ is a H-path of $BH_{n}-F$ passing through $L$.
\end{proof}

\begin{lemma}

If $|E(L_{0})\cup F_{0}|\leq 2n-6$ and $u\in V_{i}$, $v\in V_{j}$ for $i,j\in N_{4}$, and $i\neq j$, then $BH_{n}-F$ has a H-path $P[u,v]$ passing through $L$.
\end{lemma}

\begin{proof}
In this case, $|E(L_{k})\cup F_{k}|\leq 2n-6$, for each $k\in N_{4}$. In this scenario, the proofs of the cases $l=0$, $l=1$, $l=2$ and $l=3$ are analogous. We here only consider the case $l=0$.

{\it Case 1.}   $i=0$.

{\it Case 1.1.}  $j=1$.

{\it Case 1.1.1.}  $x$ (resp. $x^+$) is incident with none of $E(L_0)$ (resp. $E(L_1)$).

Suppose first that $u=x$. In this case, $v\neq x^+$. By Lemma \ref{le-2}, there is an $a\in V_0\cap Y$ such that $a$ and $a^{\pm}$ are incident with none of $E(L_{0})$ and $E(L_{3})\cup F_3$, respectively. By Lemma \ref{le-8}, there is a neighbor $y\in N_B^0(x)$ such that $(x,y)\notin E(L_0)$, $L_{0}+(x,y)$ is a linear forest, $\{u,a\}$ is compatible to $L_{0}+(x,y)$ and $y^+$ or $y^-$, say $y^+$, is incident with none of $E(L_3)$. Note that $|E(L_{0}+(x,y))\cup F_{0}|\leq 2n-5$. By the induction hypothesis, $B^{0}-F_{0}$ has a H-path $P[u,a]$ passing through $L_{0}+(x,y)$.
By Lemma \ref{le-9}, there are two neighbors $z$ and $s$ of $x^{+}$ in $B^{1}$ such that $z^+$ or $z^-$, and $s^+$ or $s^-$ are incident with none of $E(L_2)$ and $L_{1}+\{(x^{+},z),(x^{+},s)\}$ is a linear forest. We claim that there is a $d\in V_1\cap X\setminus \{z,s\}$ such that $d$ and $d^{\pm}$ are incident with none of $E(L_1)$ and $E(L_2)\cup F_2$, respectively. The reason is follows. There are |$V_1\cap X\setminus \{z,s\}|-|E(L_1)|\geq 4^{n-1}/2-(2n-6)$ candidates of $d$. Since $E(L_2)\cup F_2$ has at most $|E(L_2)\cup F_2|$ odd end vertices, each of which fails at most two candidates of such $d$. Since |$V_1\cap X\setminus \{z,s\}|-|E(L_1)|-2|E(L_0)\cup F_0|\geq 4^{n-1}/2-(2n-6)-2(2n-6)>0$, the claim holds.
Note that $\{v,d\}$ is compatible to $L_{1}+\{(x^{+},z),(x^{+},s)\}$, and $|E(L_{1}+\{(x^{+},z),(x^{+},s)\})\cup F_1|\leq 2n-4$. By the induction hypothesis, $B^{1}-F_{1}$ has a H-path $P[v,d]$ passing through $L_{1}+\{(x^{+},z),(x^{+},s)\}$. Exactly one of $z$ and $s$, say $z$, lies on the segment of $P[v,d]$ between $v$ and $x^+$. Recall that $z^+$ or $z^-$, say $z^+$, is incident with none of $E(L_2)$.
By Lemma \ref{le-9}, $a^+$ has two neighbors $c$ and $t$ in $B^{3}$ such that $c^+$ or $c^-$ (resp. $t^+$ or $t^-$), say $c^+$ (resp. $t^+$), is incident with none of $E(L_2)$, and $L_{3}+\{(a^+,c),(a^+,t)\}$ is a linear forest. Again by Lemma \ref{le-9}, there are two neighbors $b$ and $h$ of $d^+$ in $B^2$ such that $b^+$ or $b^-$ (resp. $h^+$ or $h^-$), say $b^+$ (resp. $h^+$), is incident with none of $E(L_3)$ and $L_{2}+\{(d^+,b),(d^+,h)\}$ is a linear forest. For any $g\in \{b^+,h^+\}$,
$\{y^+,g\}$ is compatible to $L_{3}+\{(a^+,c),(a^+,t)\}$ and $|E(L_{3}+\{(a^+,c),(a^+,t)\})\cup F_3|\leq 2n-4$. By the induction hypothesis, $B^{3}-F_{3}$ has a H-path $P[y^{+},g]$ passing through $L_{3}+\{(a^+,c),(a^+,t)\}$. Exactly one of $c$ and $t$, say $c$, lies on the segment of $P[y^+,g]$ between $y^+$ and $a^+$. Note that $\{z^+,c^+\}$ is compatible to $L_{2}+\{(d^+,b),(d^+,h)\}$ and $|E(L_{2}+\{(d^+,b),(d^+,h)\})\cup F_2|\leq 2n-4$. By the induction hypothesis, $B^{2}-F_{2}$ has a H-path $P[z^{+},c^+]$ passing through $L_{2}+\{(d^+,b),(d^+,h)\}$. Exactly one of $b$ and $h$, say $b$, lies on the segment of $P[z^+,c^+]$ between $z^+$ and $d^+$. Thus, $P[u,a]\cup P[v,d]\cup P[z^{+},c^+]\cup P[y^+,g]+\{(a,a^+),(b,g),(c,c^+),(d,d^+),(u,x^+),(y,y^+),(z,z^+)\}-\{(x,y),(x^+,z),(d^+,b),$ $(a^+,c)\}$ is a H-path of $BH_{n}-F$ passing through $L$.

Suppose now that $u\neq x$.
By Lemma \ref{le-9}, there are two neighbors $y$ and $w$ of $x$ in $B^{0}$ such that $y^+$ or $y^-$ (resp. $w^+$ or $w^-$), say $y^+$ (resp. $w^+$), is incident with none of $E(L_3)$ and $L_{0}+\{(x,y),(x,w)\}$ is a linear forest. We claim that there is an $a\in V_0\cap Y\setminus \{y,w\}$ such that $a$ and $a^{\pm}$ are incident with none of $E(L_0)$ and $E(L_3)\cup F_3$, respectively. The reason is follows. There are |$V_0\cap Y\setminus \{y,w\}|-|E(L_0)|\geq 4^{n-1}/2-(2n-6)$ candidates of $a$. Since $E(L_3)\cup F_3$ has at most $|E(L_3)\cup F_3|$ even end vertices, each of which fails at most two candidates of such $a$. Since |$V_0\cap Y\setminus \{y,w\}|-|E(L_0)|-2|E(L_3)\cup F_3|\geq 4^{n-1}/2-(2n-6)-2(2n-6)>0$, the claim holds. Note that $\{u,a\}$ is compatible to $L_{0}+\{(x,y),(x,w)\}$, and $|E(L_{0}+\{(x,y),(x,w)\})\cup F_0|\leq 2n-4$. By the induction hypothesis, $B^{0}-F_{0}$ has a H-path $P[u,a]$ passing through $L_{0}+\{(x,y),(x,w)\}$. Exactly one of $y$ and $w$, say $y$, lies on the segment of $P[u,a]$ between $u$ and $x$. By Lemma \ref{le-2}, there is a $d\in V_1\cap X$ such that $d$ and $d^{\pm}$ are incident with none of $E(L_1)$ and $E(L_2)\cup F_2$, respectively. By Lemma \ref{le-8}, there is a $z\in N_B^1(x^+)$ such that $(x^+,z)\notin E(L_1)$, $L_{1}+(x^+,z)$ is a linear forest, $\{v,d\}$ is compatible to $L_{1}+(x^+,z)$ and $z^+$ or $z^-$, say $z^+$, is incident with none of $E(L_2)$. Note that $|E(L_{1}+\{(x^+,z)\})\cup F_{1}|\leq 2n-5$. By the induction hypothesis, $B^{1}-F_{1}$ has a H-path $P[v,d]$ passing through $L_{1}+(x^+,z)$. Let $w=d^-$, if $z=d$; and $w=d^+$, otherwise. Then $w\neq z^+$. By Lemma \ref{le-9}, $a^+$ has two neighbors $c$ and $t$ in $B^{3}$ such that $c^+$ or $c^-$ (resp. $t^+$ or $t^-$), say $c^+$ (resp. $t^+$), is incident with none of $E(L_2)$, and $L_{3}+\{(a^+,c),(a^+,t)\}$ is a linear forest. Again by Lemma \ref{le-9}, there are two neighbors $b$ and $h$ of $w$ in $B^2$ such that $b^+$ or $b^-$ (resp. $h^+$ or $h^-$), say $b^+$ (resp. $h^+$), is incident with none of $E(L_3)$ and $L_{2}+\{(w,b),(w,h)\}$ is a linear forest. For any $g\in \{b^+,h^+\}$,
$\{y^+,g\}$ is compatible to $L_{3}+\{(a^+,c),(a^+,t)\}$ and $|E(L_{3}+\{(a^+,c),(a^+,t)\})\cup F_3|\leq 2n-4$. By the induction hypothesis, $B^{3}-F_{3}$ has a H-path $P[y^{+},g]$ passing through $L_{3}+\{(a^+,c),(a^+,t)\}$. Exactly one of $c$ and $t$, say $c$, lies on the segment of $P[y^+,g]$ between $y^+$ and $a^+$. Note that $\{z^+,c^+\}$ is compatible to $L_{2}+\{(w,b),(w,h)\}$ and $|E(L_{2}+\{(w,b),(w,h)\})\cup F_2|\leq 2n-4$. By the induction hypothesis, $B^{2}-F_{2}$ has a H-path $P[z^{+},c^+]$ passing through $L_{2}+\{(w,b),(w,h)\}$. Exactly one of $b$ and $h$, say $b$, lies on the segment of $P[z^+,c^+]$ between $z^+$ and $w$. Thus, $P[u,a]\cup P[v,d]\cup P[z^{+},c^+]\cup P[y^+,g]+\{(a,a^+),(b,g),(c,c^+),(d,w),(x,x^+),(y,y^+),(z,z^+)\}-\{(x,y),(x^+,z),(w,b),$ $(a^+,c)\}$ is a H-path of $BH_{n}-F$ passing through $L$.

{\it Case 1.1.2.}  $x$ is incident with none of $E(L_0)$ and $L_1$ has a maximal path $P[x^+,w]$ with $w\neq x^+$.

In this case, $v\neq x^{+}$. By Lemma \ref{le-2}, there is an $a\in V_0\cap Y$ such that $a$ (resp. $a^{\pm}$) is incident with none of $E(L_0)$ (resp. $E(L_3)\cup F_3$).

Suppose first that $w\neq v$.
By Lemma \ref{le-8}, there is a neighbor $y$ of $x$ in $B^0$ such that $(x,y)\notin E(L_0)$, $L_{0}+(x,y)$ is a linear forest, $\{u,a\}$ is compatible to $L_{0}+(x,y)$, and $y^+$ or $y^-$, say $y^+$, is incident with none of $E(L_3)$. Note that $|E(L_{0}+\{(x,y)\})\cup F_{0}|\leq 2n-5$. By the induction hypothesis, $B^{0}-F_{0}$ has a H-path $P[u,a]$ passing through $L_{0}+(x,y)$. Let $g=a^-$ if $y=a$ and let $g=a^+$ otherwise. Then $g\neq y^+$.
Let $(x^+,s)\in E(P[x^+,w])$. By Lemma \ref{le-10}, there are two distinct vertices $z,t\in N_B^1(x^+)\setminus \{s\}$ such that $L_{1}+\{(x^+,z),(x^+,t)\}-(x^+,s)$ is a linear forest, $z$ is not the shadow vertex of $t$, $z^+$ or $z^-$, say $z^+$, is incident with none of $E(L_2)\cup F_2$ and $t^+$ or $t^-$, say $t^+$, is not an internal vertex of $L_2$.
Note that $\{v,s\}$ is compatible to $L_{1}+\{(x^+,z),(x^+,t)\}-(x^+,s)$ and $|E(L_{1}+\{(x^+,z),(x^+,t)\}-(x^+,s))\cup F_1|\leq 2n-5$. By the induction hypothesis, $B^1-F_1$ has a H-path $P[v,s]$ passing through $L_{1}+\{(x^+,z),(x^+,t)\}-(x^+,s)$. By Lemma \ref{le-9}, $g$ has two neighbors $c$ and $r$ in $B^{3}$ such that $c^+$ or $c^-$ (resp. $r^+$ or $r^-$), say $c^+$ (resp. $r^+$), is incident with none of $E(L_2)$, and $L_{3}+\{(g,c),(g,r)\}$ is a linear forest. Again by Lemma \ref{le-9}, there are two neighbors $b$ and $h$ of $z^+$ in $B^2$ such that $b^+$ or $b^-$ (resp. $h^+$ or $h^-$), say $b^+$ (resp. $h^+$), is incident with none of $E(L_3)$ and $L_{2}+\{(z^+,b),(z^+,h)\}$ is a linear forest. For any $d\in \{b^+,h^+\}$,
$\{y^+,d\}$ is compatible to $L_{3}+\{(g,c),(g,r)\}$ and $|E(L_{3}+\{(g,c),(g,r)\})\cup F_3|\leq 2n-4$. By the induction hypothesis, $B^{3}-F_{3}$ has a H-path $P[y^{+},d]$ passing through $L_{3}+\{(g,c),(g,r)\}$. Let $q$ be the neighbor of $g$ on the segment of $P[y^+,d]$ between $y^+$ and $g$, if $y$ lies on the segment of $P[u,a]$ between $x$ and $a$; and let $q$ be the neighbor of $g$ on the segment of $P[y^+,d]$ between $y^+$ and $d$, otherwise.
Note that $\{t^+,q^+\}$ is compatible to $L_{2}+\{(z^+,b),(z^+,h)\}$ and $|E(L_{2}+\{(z^+,b),(z^+,h))\}\cup F_2|\leq 2n-4$. By the induction hypothesis, $B^{2}-F_{2}$ has a H-path $P[t^{+},q^+]$ passing through $L_{2}+\{(z^+,b),(z^+,h)\}$. Exactly one of $b$ and $h$, say $b$, lies on the segment of $P[t^+,q^+]$ between $z^+$ and $t^+$. Thus, $P[u,a]\cup P[v,s]\cup P[t^{+},q^+]\cup P[y^+,d]+\{(x^+,s),(a,g),(b,d),(q,q^+),(t,t^+),(x,x^+),(y,y^+),(z,z^+)\}-\{(x,y),(x^+,z),(t,$ $x^+),(z^+,b),$ $(g,q)\}$ is a H-path of $BH_{n}-F$ passing through $L$.

Suppose now that $w=v$. In this case, $u\neq x$.
By Lemma \ref{le-9}, there are two neighbors $y$ and $g$ of $x$ in $B^0$ such that $y^+$ or $y^-$ (resp. $g^+$ or $g^-$), say $y^+$ (resp. $g^+$), is incident with none of $E(L_3)$ and $L_{0}+\{(x,y),(x,g)\}$ is a linear forest. We claim that there is an $a\in V_0\cap Y\setminus \{y,g\}$ such that $a$ and $a^{\pm}$ are incident with none of $E(L_0)$ and $E(L_3)\cup F_3$, respectively. The reason is follows. There are |$V_0\cap Y\setminus \{y,g\}|-|E(L_0)|\geq 4^{n-1}/2-(2n-6)$ candidates of $a$. Since $E(L_3)\cup F_3$ has at most $|E(L_3)\cup F_3|$ even end vertices, each of which fails at most two candidates of such $a$. Since |$V_0\cap Y\setminus \{y,g\}|-|E(L_0)|-2|E(L_3)\cup F_3|\geq 4^{n-1}/2-(2n-6)-2(2n-6)>0$, the claim holds. Note that $\{u,a\}$ is compatible to $L_{0}+\{(x,y),(x,g)\}$, and $|E(L_{0}+\{(x,y),(x,g)\})\cup F_0|\leq 2n-4$. By the induction hypothesis, $B^{0}-F_{0}$ has a H-path $P[u,a]$ passing through $L_{0}+\{(x,y),(x,g)\}$. Exactly one of $y$ and $g$, say $y$, lies on the segment of $P[u,a]$ between $u$ and $x$.
By Lemma \ref{le-2}, there is an $s\in V_1\cap X$ such that $s$ and $s^{\pm}$ are incident with none of $E(L_1)$ and $E(L_2)\cup F_2$, respectively. By Lemma \ref{le-8}, there is a $t\in N_B^1(x^+)$ such that $(x^+,t)\notin E(L_1)$, $L_{1}+(x^+,t)$ is a linear forest, $\{v,s\}$ is compatible to $L_{1}+(x^+,t)$ and $t^+$ or $t^-$, say $t^+$, is incident with none of $E(L_2)$. By the induction hypothesis, $B^1-F_1$ has a H-path $P[v,s]$ passing through $L_{1}+(x^+,t)$. In this case, $t\neq s$.
By Lemma \ref{le-9}, $a^+$ has two neighbors $c$ and $r$ in $B^{3}$ such that $c^+$ or $c^-$ (resp. $r^+$ or $r^-$), say $c^+$ (resp. $r^+$), is incident with none of $E(L_2)$, and $L_{3}+\{(a^+,c),(a^+,r)\}$ is a linear forest. Again by Lemma \ref{le-9}, there are two neighbors $b$ and $h$ of $s^+$ in $B^2$ such that $b^+$ or $b^-$ (resp. $h^+$ or $h^-$), say $b^+$ (resp. $h^+$), is incident with none of $E(L_3)$ and $L_{2}+\{(s^+,b),(s^+,h)\}$ is a linear forest. For any $d\in \{b^+,h^+\}$,
$\{y^+,d\}$ is compatible to $L_{3}+\{(a^+,c),(a^+,r)\}$ and $|E(L_{3}+\{(a^+,c),(a^+,r)\})\cup F_3|\leq 2n-4$. By the induction hypothesis, $B^{3}-F_{3}$ has a H-path $P[y^{+},d]$ passing through $L_{3}+\{(a^+,c),(a^+,r)\}$. Exactly one of $c$ and $r$, say $c$, lies on the segment of $P[y^+,d]$ between $y^+$ and $a^+$.
Note that $\{t^+,c^+\}$ is compatible to $L_{2}+\{(s^+,b),(s^+,h)\}$ and $|E(L_{2}+\{(s^+,b),(s^+,h))\}\cup F_2|\leq 2n-4$. By the induction hypothesis, $B^{2}-F_{2}$ has a H-path $P[t^{+},c^+]$ passing through $L_{2}+\{(s^+,b),(s^+,h)\}$. Exactly one of $b$ and $h$, say $b$, lies on the segment of $P[t^+,c^+]$ between $s^+$ and $t^+$. Thus, $P[u,a]\cup P[v,s]\cup P[t^{+},c^+]\cup P[y^+,d]+\{(a,a^+),(b,d),(c,c^+),(t,t^+),(x,x^+),(y,y^+),(z,z^+)\}-\{(x,y),(x^+,t),(s^+,b),(a^+,c)\}$ is a H-path of $BH_{n}-F$ passing through $L$.

{\it Case 1.1.3.}  $L_0$ has a maximal path $P[x,r]$ with $r\neq x$ and $x^+$ is incident with none of $L_1$.

Suppose first that $u=r$. In this case, $v\neq x^+$. By Lemma \ref{le-2}, there is an $a\in V_0\cap Y$ such that $a$ and $a^{\pm}$ are incident with none of $E(L_0)$ and $E(L_3)\cup F_3$, respectively. By Lemma \ref{le-8}, there is a $y\in N_B^0(x)$ such that $(x,y)\notin E(L_0)$, $L_{0}+(x,y)$ is a linear forest, $\{u,a\}$ is compatible to $L_0+(x,y)$, and $y^+$ or $y^-$, say $y^+$, is incident with none of $E(L_3)$. By the induction hypothesis, $B^0-F_0$ has a H-path $P[u,a]$ passing through $L_{0}+(x,y)$. Then $a\neq y$.
By Lemma \ref{le-9}, there are two neighbors $z$ and $s$ of $x^{+}$ in $B^{1}$ such that $z^+$ or $z^-$ (resp. $s^+$ or $s^-$), say $z^+$ (resp. $s^+$), is incident with none of $E(L_2)$ and $L_{1}+\{(x^{+},z),(x^{+},s)\}$ is a linear forest. We claim that there is a $d\in V_1\cap X\setminus \{z,s\}$ such that $d$ and $d^{\pm}$ are incident with none of $E(L_1)$ and $E(L_2)\cup F_2$, respectively. The reason is follows. There are |$V_1\cap X\setminus \{z,s\}|-|E(L_1)|\geq 4^{n-1}/2-(2n-6)$ candidates of $d$. Since $E(L_2)\cup F_2$ has at most $|E(L_2)\cup F_2|$ odd end vertices, each of which fails at most two candidates of such $d$. Since |$V_1\cap X\setminus \{z,s\}|-|E(L_1)|-2|E(L_2)\cup F_2|\geq 4^{n-1}/2-(2n-6)-2(2n-6)>0$, the claim holds.
Note that $\{v,d\}$ is compatible to $L_{1}+\{(x^{+},z),(x^{+},s)\}$, and $|E(L_{1}+\{(x^{+},z),(x^{+},s)\})\cup F_1|\leq 2n-4$. By the induction hypothesis, $B^{1}-F_{1}$ has a H-path $P[v,d]$ passing through $L_{1}+\{(x^{+},z),(x^{+},s)\}$. Exactly one of $z$ and $s$, say $z$, lies on the segment of $P[v,d]$ between $v$ and $x^+$.
By Lemma \ref{le-9}, $a^+$ has two neighbors $c$ and $t$ in $B^{3}$ such that $c^+$ or $c^-$ (resp. $t^+$ or $t^-$), say $c^+$ (resp. $t^+$), is incident with none of $E(L_2)$, and $L_{3}+\{(a^+,c),(a^+,t)\}$ is a linear forest. Again by Lemma \ref{le-9}, there are two neighbors $b$ and $h$ of $d^+$ in $B^2$ such that $b^+$ or $b^-$ (resp. $h^+$ or $h^-$), say $b^+$ (resp. $h^+$), is incident with none of $E(L_3)$ and $L_{2}+\{(d^+,b),(d^+,h)\}$ is a linear forest. For any $g\in \{b^+,h^+\}$,
$\{y^+,g\}$ is compatible to $L_{3}+\{(a^+,c),(a^+,t)\}$ and $|E(L_{3}+\{(a^+,c),(a^+,t)\})\cup F_3|\leq 2n-4$. By the induction hypothesis, $B^{3}-F_{3}$ has a H-path $P[y^{+},g]$ passing through $L_{3}+\{(a^+,c),(a^+,t)\}$. Exactly one of $c$ and $t$, say $c$, lies on the segment of $P[y^+,g]$ between $y^+$ and $a^+$. Note that $\{z^+,c^+\}$ is compatible to $L_{2}+\{(d^+,b),(d^+,h)\}$ and $|E(L_{2}+\{(d^+,b),(d^+,h)\})\cup F_2|\leq 2n-4$. By the induction hypothesis, $B^{2}-F_{2}$ has a H-path $P[z^{+},c^+]$ passing through $L_{2}+\{(d^+,b),(d^+,h)\}$. Exactly one of $b$ and $h$, say $b$, lies on the segment of $P[z^+,c^+]$ between $z^+$ and $d^+$. Thus, $P[u,a]\cup P[v,d]\cup P[z^{+},c^+]\cup P[y^+,g]+\{(a,a^+),(b,g),(c,c^+),(d,d^+),(x,x^+),(y,y^+),(z,z^+)\}-\{(x,y),(x^+,z),(d^+,b),$ $(a^+,c)\}$ is a H-path of $BH_{n}-F$ passing through $L$.

Suppose now that $u\neq r$.
Let $(x,w)\in E(P[x,r])$. By Lemma \ref{le-10}, there are two distinct vertices $y,a\in N_B^0(x)\setminus \{w\}$ such that $L_{0}+\{(x,y),(x,a)\}-(x,w)$ is a linear forest, $y$ is not the shadow vertex of $a$, $a^+$ or $a^-$, say $a^+$, is incident with none of $E(L_3)\cup F_3$ and $y^+$ or $y^-$, say $y^+$, is not an internal vertex of $L_3$.
Note that $\{u,w\}$ is compatible to $L_{0}+\{(x,y),(x,a)\}-(x,w)$ and $|E(L_{0}+\{(x,y),(x,a)\}-(x,w))\cup F_0|\leq 2n-5$. By the induction hypothesis, $B^0-F_0$ has a H-path $P[u,w]$ passing through $L_{0}+\{(x,y),(x,a)\}-(x,w)$.
By Lemma \ref{le-2}, there is a $d\in V_1\cap X$ such that $d$ and $d^{\pm}$ are incident with none of $E(L_1)$ and $E(L_2)\cup F_2$, respectively.
By Lemma \ref{le-8}, there is a $z\in N_B^1(x^+)$ such that $(x^+,z)\notin E(L_1)$, $L_{1}+(x^+,z)$ is a linear forest, $\{v,d\}$ is compatible to $L_{1}+(x^+,z)$, and $z^+$ or $z^-$, say $z^+$, is incident with none of $E(L_2)$. Note that $|E(L_{1}+\{(x^+,z)\})\cup F_{1}|\leq 2n-5$. By the induction hypothesis, $B^{1}-F_{1}$ has a H-path $P[v,d]$ passing through $L_{1}+(x^+,z)$. Let $g=d^-$ if $z=d$ and let $g=d^+$ otherwise. Then $g\neq z^+$.
By Lemma \ref{le-9}, $a^+$ has two neighbors $c$ and $t$ in $B^{3}$ such that $c^+$ or $c^-$ (resp. $t^+$ or $t^-$), say $c^+$ (resp. $t^+$), is incident with none of $E(L_2)$, and $L_{3}+\{(a^+,c),(a^+,t)\}$ is a linear forest. Again by Lemma \ref{le-9}, there are two neighbors $b$ and $h$ of $g$ in $B^2$ such that $b^+$ or $b^-$ (resp. $h^+$ or $h^-$), say $b^+$ (resp. $h^+$), is incident with none of $E(L_3)$ and $L_{2}+\{(g,b),(g,h)\}$ is a linear forest. For any $s\in \{b^+,h^+\}$,
$\{y^+,s\}$ is compatible to $L_{3}+\{(a^+,c),(a^+,t)\}$ and $|E(L_{3}+\{(a^+,c),(a^+,t)\})\cup F_3|\leq 2n-4$. By the induction hypothesis, $B^{3}-F_{3}$ has a H-path $P[y^{+},s]$ passing through $L_{3}+\{(a^+,c),(a^+,t)\}$. Let $q$ be the neighbor of $a^+$ on the segment of $P[y^+,s]$ between $a^+$ and $s$, if $z$ lies on the segment of $P[v,d]$ between $x^+$ and $v$; and let $q$ be the neighbor of $a^+$ on the segment of $P[y^+,s]$ between $a^+$ and $y^+$, otherwise.
Note that $\{z^+,q^+\}$ is compatible to $L_{2}+\{(g,b),(g,h)\}$ and $|E(L_{2}+\{(g,b),(g,h)\})\cup F_2|\leq 2n-4$. By the induction hypothesis, $B^{2}-F_{2}$ has a H-path $P[z^{+},q^+]$ passing through $L_{2}+\{(g,b),(g,h)\}$. Exactly one of $b$ and $h$, say $b$, lies on the segment of $P[z^+,q^+]$ between $z^+$ and $g$. Thus, $P[u,w]\cup P[v,d]\cup P[z^{+},q^+]\cup P[y^+,s]+\{(x,w),(a,a^+),(b,s),(d,g),(q,q^+),(x,x^+),(y,y^+),(z,z^+)\}-\{(x,y),(x,a),(x^+,$ $z),(g,b),$ $(a^+,q)\}$ is a H-path of $BH_{n}-F$ passing through $L$.

{\it Case 1.1.4.}  $L_0$ has a maximal path $P[x,r]$ with $r\neq x$ and $L_1$ has a maximal path $P[x^+,w]$ with $w\neq x^+$.

In this case, $u\neq x$, $v\neq x^{+}$. Since $\{u,v\}$ is compatible to $L$, let $P[u,a]$ is a maximal path in $L_0$ and let $P[v,b]$ is a maximal path in $L_1$, we has $\{a,b\}\cap \{r,w\}=\emptyset$. If $u\neq r$ is similarly to the Case 1.1.3 $u\neq r$. If $u=r$, then $v\neq w$ is similarly to the Case 1.1.2 $v\neq w$.

{\it Case 1.2.}  $j=2$.

By Lemma \ref{le-2}, there is an $a\in V_0\cap Y$ such that $a$ and $a^{\pm}$ are incident with none of $E(L_0)$ and $E(L_3)\cup F_3$, respectively. By Lemma \ref{le-8}, there is a $y\in N_B^0(x)$ such that $(x,y)\notin E(L_0)$, $L_0+(x,y)$ is a linear forest and $y^+$ or $y^-$, say $y^+$, is incident with none of $E(L_3)$. Note that $\{u,a\}$ is compatible to $L_0+(x,y)$. By the induction hypothesis, $B^0-F_0$ has a H-path $P[u,a]$ passing through $L_0+(x,y)$. Let $g=a^-$, if $y=a$; and $g=a^+$, otherwise. Then $g\neq y^+$. By Lemma \ref{le-9}, $g$ has two neighbors $w$ and $t$ in $B^{3}$ such that $w^{+}$ or $w^-$ (resp. $t^{+}$ or $t^-$), say $w^+$ (resp. $t^+$), is incident with none of $E(L_2)$,
and $L_{3}+\{(g,w),(g,t)\}$ is a linear forest. We claim that there is a $c\in V_3\cap Y\setminus \{w,t\}$ such that $c$ and $c^{\pm}$ are incident with none of $E(L_3)$ and $E(L_2)\cup F_2$, respectively, and $v$ is not adjacent to $c^{\pm}$. The reason is follows. There are |$V_3\cap Y\setminus \{w,t\}|-|E(L_3)|\geq 4^{n-1}/2-(2n-6)$ candidates of $c$. Since $E(L_2)\cup F_2$ has at most $|E(L_2)\cup F_2|$ even end vertices, each of which fails at most two candidates of such $c$. Since there are $|N_B^2(v)|=2n-2$ vertices adjacent to $v$. Since |$V_3\cap Y\setminus \{w,t\}|-|E(L_3)|-2|E(L_2)\cup F_2|-|N_B^2(v)|\geq 4^{n-1}/2-(2n-6)-2(2n-6)-(2n-2)>0$, the claim holds.
Note that $\{y^+,c\}$ is compatible to $L_{3}+\{(g,w),(g,t)\}$. By the induction hypothesis, $B^3-F_3$ has a H-path $P[y^+,c]$ passing through $L_{3}+\{(g,w),(g,t)\}$. Exactly one of $w$ and $t$, say $w$, lies on the segment of $P[y^+,c]$ between $g$ and $y^+$. By Lemma \ref{le-9}, $c^+$ has two neighbors $z$ and $d$ in $B^{2}$ such that $z^{+}$ or $z^-$ (resp. $t^{+}$ or $t^-$), say $z^+$ (resp. $t^+$), is incident with none of $E(L_1)$, and $L_{2}+\{(c^+,z),(c^+,d)\}$ is a linear forest. Note that $\{w^+,v\}$ is compatible to $L_{2}+\{(c^+,z),(c^+,d)\}$. By the induction hypothesis, $B^2-F_2$ has a H-path $P[w^+,v]$ passing through $L_{2}+\{(c^+,z),(c^+,d)\}$. Exactly one of $z$ and $d$, say $z$, lies on the segment of $P[w^+,v]$ between $c^+$ and $w^+$. By the induction hypothesis, $B^1-F_1$ has a H-path $P[x^+,z^+]$ passing through $L_{1}$. Thus, $P[u,a]\cup P[x^{+},z^{+}]\cup P[w^{+},v]\cup P[y^{+},c]+\{(a,g),(c,c^+),(w,w^{+}),(x,x^{+}),(y,y^{+}),(z,z^{+})\}-\{(x,y),(c^{+},z),(g,w)\}$ is a desired H-path of $BH_n-F$.

{\it Case 1.3.}  $j=3$.

By Lemma \ref{le-4}, there is an $a\in V_0\cap Y$ such that $a$ and $a^{\pm}$ are incident with none of $E(L_0)$ and $E(L_3)\cup F_3$, respectively, and $v$ is not adjacent to $a^{\pm}$. By Lemma \ref{le-8}, there is a $y\in N_B^0(x)$ such that $(x,y)\notin E(L_0)$, $L_0+(x,y)$ is a linear forest and $y^+$ or $y^-$, say $y^+$, is incident with none of $E(L_3)$. Note that $\{u,a\}$ is compatible to $L_0+(x,y)$. By the induction hypothesis, $B^0-F_0$ has a H-path $P[u,a]$ passing through $L_0+(x,y)$. Let $g=a^-$, if $y=a$; and $g=a^+$, otherwise. Then $g\neq y^+$. By Lemma \ref{le-9}, $g$ has two neighbors $w$ and $t$ in $B^{3}$ such that $w^{+}$ or $w^-$ (resp. $t^{+}$ or $t^-$), say $w^+$ (resp. $t^+$), is incident with none of $E(L_2)$, and $L_{3}+\{(g,w),(g,t)\}$ is a linear forest. Note that $\{y^+,v\}$ is compatible to $L_{3}+\{(g,w),(g,t)\}$. By the induction hypothesis, $B^3-F_3$ has a H-path $P[y^+,v]$ passing through $L_{3}+\{(g,w),(g,t)\}$. Exactly one of $w$ and $t$, say $w$, lies on the segment of $P[y^+,v]$ between $g$ and $y^+$. By Lemma \ref{le-2}, there is a $z\in V_1\cap X$ such that $z$ (resp. $z^+$) is incident with none of $E(L_1)$ (resp. $E(L_2)$). By the induction hypothesis, $B^1-F_1$, $B^2-F_2$ have H-paths $P[x^+,z]$, $P[w^+,z^+]$ passing through $L_1$ and $L_2$, respectively. Thus, $P[u,a]\cup P[x^{+},z]\cup P[w^{+},z^+]\cup P[y^{+},v]+\{(a,g),(w,w^{+}),(x,x^{+}),(y,y^{+}),(z,z^{+})\}-\{(x,y),(g,w)\}$ is a desired H-path of $BH_n-F$.

{\it Case 2.}  $i\neq 0$.

{\it Case 2.1.}  $i=1,j=2$.

Suppose first that $\{u,x^+\}$ is compatible to $L_1$. By Lemma \ref{le-2}, there are vertices $y\in V_0\cap Y$ and $z\in V_3\cap Y$ such that $y$ (res. $y^+$) is incident with none of $E(L_0)$ (resp. $E(L_3)$), and $z$ (resp. $z^+$) is incident with none of $E(L_3)$ (resp. $E(L_2)$). By the induction hypothesis, $B^0-F_0$, $B^1-F_1$, $B^2-F_2$, $B^3-F_3$ have H-paths $P[x,y]$, $P[u,x^+]$, $P[v,z^+]$ and $P[y^+,z]$ passing through $L_0$, $L_1$, $L_2$ and $L_3$, respectively. Thus, $P[x,y]\cup P[x^{+},u]\cup P[v,z^{+}]\cup P[y^{+},z]+\{(x,x^{+}),(y,y^{+}),(z,z^{+})\}$ is a H-path of $BH_{n}-F$ passing through $L$.

Suppose now that $L_1$ has a maximal path $P[u,x^+]$. By Lemma \ref{le-2}, there is an $a\in V_1\cap Y$, such that $a$ and $a^{\pm}$ are incident with none of $E(L_1)$ and $E(L_0)$, respectively. By Lemma \ref{le-8}, there is a $z\in N_B^1(x^+)$ such that $(x^+,z)\notin E(L_1)$, $L_1+(x^+,z)$ is a linear forest and $z^+$ or $z^-$, say $z^+$, is incident with none of $E(L_2)$. Note that $\{u,a\}$ is compatible to $L_1+(x^+,z)$. By the induction hypothesis, $B^1-F_1$ has a H-path $P[u,a]$ passing through $L_1+(x^+,z)$. Since $a\neq x^+$, $a^-\neq x$.  By Lemma \ref{le-2}, there is a $d\in V_2\cap X$, such that $d$ and $d^{\pm}$ are incident with none of $E(L_2)$ and $E(L_3)\cup F_3$, respectively. By Lemma \ref{le-8}, there is a $c\in N_B^2(z^+)$ such that $(z^+,c)\notin E(L_2)$, $L_2+(z^+,c)$ is a linear forest and $c^+$ or $c^-$, say $c^+$, is incident with none of $E(L_3)$. Note that $\{v,d\}$ is compatible to $L_2+(z^+,c)$. By the induction hypothesis, $B^2-F_2$ has a H-path $P[v,d]$ passing through $L_2+(z^+,c)$. Let $g=d^-$, if $c=d$; and $g=d^+$, otherwise. Then $g\neq c^+$. By Lemma \ref{le-9}, $a^-$ has two neighbors $b$ and $t$ in $B^{0}$ such that $b^+$ or $b^-$ (resp. $t^+$ or $t^-$), say $b^+$ (resp. $t^+$), is incident with none of $E(L_3)$, and $L_{0}+\{(a^-,b),(a^-,t)\}$ is a linear forest. Again by Lemma \ref{le-9}, there are two neighbors $w$ and $r$ of $g$ in $B^3$ such that $w^+$ or $w^-$ (resp. $r^+$ or $r^-$), say $w^+$ (resp. $r^+$), is incident with none of $E(L_0)$ and $L_{3}+\{(g,w),(g,r)\}$ is a linear forest. For any $h\in \{w^+,r^+\}$,
$\{x,h\}$ is compatible to $L_{0}+\{(a^-,b),(a^-,t)\}$ and $|E(L_{0}+\{(a^-,b),(a^-,t)\})\cup F_0|\leq 2n-4$. By the induction hypothesis, $B^{0}-F_{0}$ has a H-path $P[x,h]$ passing through $L_{0}+\{(a^-,b),(a^-,t)\}$. Exactly one of $b$ and $t$, say $b$, lies on the segment of $P[x,h]$ between $a^-$ and $x$.
Note that $\{c^+,b^+\}$ is compatible to $L_{3}+\{(g,w),(g,r)\}$ and $|E(L_{3}+\{(g,w),(g,r)\})\cup F_3|\leq 2n-4$. By the induction hypothesis, $B^{3}-F_{3}$ has a H-path $P[c^{+},b^+]$ passing through $L_{3}+\{(g,w),(g,r)\}$. Exactly one of $w$ and $r$, say $w$, lies on the segment of $P[c^+,b^+]$ between $g$ and $c^+$. Thus, $P[x,h]\cup P[u,a]\cup P[v,d]\cup P[c^{+},b]+\{(a,a^-),(b,b^+),(c,c^+),(d,g),(w,h),(x,x^{+}),(z,z^{+})\}-\{(a^-,b),(x^{+},z),(z^+,c),$ $(g,w)\}$ is a desired H-path of $BH_n-F$.

{\it Case 2.2.}  $i=1,j=3$.

By Lemma \ref{le-2}, there is an $a\in V_1\cap Y$ such that $a$ and $a^{\pm}$ are incident with none of $E(L_1)$ and $E(L_0)$, respectively. By Lemma \ref{le-8}, there is a $z\in N_B^1(x^+)$ such that $L_1+(x^+,z)$ is a linear forest and $z^+$ or $z^-$, say $z^+$, is incident with none of $E(L_2)$. Note that $\{u,a\}$ is compatible to $L_1+(x^+,z)$. By the induction hypothesis, $B^1-F_1$ has a H-path $P[u,a]$ passing through $L_1+(x^+,z)$. Let $g=a^-$, if $a\neq x^+$; and $g=a^+$, otherwise. Then $g\neq x$. By Lemma \ref{le-9}, $g$ has two neighbors $y$ and $t$ in $B^0$ such that $b\notin \{y,t\}$, $y^{+}$ or $y^-$ (resp. $t^{+}$ or $t^-$), say $y^+$ (resp. $t^+$), is incident with none of $E(L_3)$ and $L_0+\{(g,y),(g,t)\}$ is a linear forest.
We claim that there is a $b\in V_0\cap Y\setminus \{y,t\}$ such that $b$ and $b^{\pm}$ are incident with none of $E(L_0)$ and $E(L_3)\cup F_3$, respectively, and $v$ is not adjacent to $b^{\pm}$. The reason is follows. There are |$V_0\cap Y\setminus \{y,t\}|-|E(L_0)|\geq 4^{n-1}/2-(2n-6)$ candidates of $b$. Since $E(L_3)\cup F_3$ has at most $|E(L_3)\cup F_3|$ even end vertices, each of which fails at most two candidates of such $b$. Since there are $|N_B^3(v)|=2n-2$ vertices adjacent to $v$. Since |$V_0\cap Y\setminus \{y,t\}|-|E(L_0)|-2|E(L_3)\cup F_3|-|N_B^3(v)|\geq 4^{n-1}/2-(2n-6)-2(2n-6)-(2n-2)>0$, the claim holds. Note that $\{x,b\}$ is compatible to $L_0+\{(g,y),(g,t)\}$.
By the induction hypothesis, $B^0-F_0$ has a H-path $P[x,b]$ passing through $L_0+\{(g,y),(g,t)\}$. Exactly one of $y$ and $t$, say $y$, lies on the segment of $P[x,b]$ between $x$ and $g$. By Lemma \ref{le-9}, $b^+$ has two neighbors $w$ and $s$ in $B^3\setminus \{v\}$ such that $w^{+}$ or $w^-$ (resp. $s^{+}$ or $s^-$), say $w^+$ (resp. $s^+$), is incident with none of $E(L_2)$ and $L_3+\{(b^+,w),(b^+,s)\}$ is a linear forest. Note that $\{y^+,v\}$ is compatible to $L_3+\{(b^+,w),(b^+,s)\}$. By the induction hypothesis, $B^3-F_3$ has a H-path $P[y^+,v]$ passing through $L_3+\{(b^+,w),(b^+,s)\}$. Exactly one of $w$ and $s$, say $w$, lies on the segment of $P[y^+,v]$ between $b^+$ and $y^+$. By the induction hypothesis, $B^2-F_2$ has a H-path $P[w^+,z^+]$ passing through $L_2$. Thus, $P[x,b]\cup P[u,a]\cup P[w^+,z^{+}]\cup P[y^{+},v]+\{(a,g),(b,b^+),(w,w^+),(x,x^{+}),(y,y^{+}),(z,z^{+})\}-\{(g,y),(x^+,z),(b^+,w)\}$ is a H-path of $BH_{n}-F$ passing through $L$.

{\it Case 2.3.}  $i=2,j=3$.

By Lemma \ref{le-2}, there are vertices $y\in V_0\cap Y$ and $z\in V_1\cap X$ such that $y$ (resp. $y^+$) is incident with none of $E(L_0)$ (resp. $E(L_3)$), and $z$ (resp. $z^+$) is incident with none of $E(L_1)$ (resp. $E(L_2)$). By the induction hypothesis, $B^0-F_0$, $B^1-F_1$, $B^2-F_2$, $B^3-F_3$ have H-paths $P[x,y]$, $P[x^+,z]$, $P[z^+,u]$ and $P[y^+,v]$ passing through $L_0$, $L_1$, $L_2$ and $L_3$, respectively. Thus, $P[x,y]\cup P[x^{+},z]\cup P[z^{+},u]\cup P[y^{+},v]+\{(x,x^{+}),(y,y^{+}),(z,z^{+})\}$ is a H-path of $BH_{n}-F$ passing through $L$.

\begin{lemma}
If $|E(L_{0})\cup F_{0}|=2n-3$, then $BH_{n}-F$ contains a H-path $P[u,v]$ passing through $L$.
\end{lemma}

Proof. In this case, $E(L_{k})\cup F_{k}=\emptyset$ for $k\in N_{4}\setminus \{0\}$. By Lemma \ref{le-1} and Theorem \ref{th-li2019}, $B^{0}-F_{0}$ has a H-cycle $C_{0}$ passing through $L_{0}$. 

{\it Case 1.}   $u,v\in V_i$.

{\it Case 1.1.}   $l=0$ or $l=3$.

The proofs of the cases $l=0$ and $l=3$ are analogous. We here consider the case $l=0$.

{\it Case 1.1.1.}  $i=0$.

Since $F_0=F\neq \emptyset$, let $f\in F_{0}$. By the induction hypothesis, $B^{0}-F_{0}\setminus \{f\}$ has a H-path $P[u,v]$ passing through $L_{0}$. Let $(x,y)\in E(P[u,v])\setminus E(L_{0})$.
Let $z,c\in V_1\cap X$, $d,w\in V_2\cap X$ be pair-wires distinct.

Suppose first that $f\notin E(P[u,v])$. By Theorem \ref{th-xu2007}, $B^{1}$, $B^2$, $B^3$ have H-paths $P[x^{+},z]$, $P[z^{+},w]$ and $P[w^{+},y^{+}]$, respectively. Thus, $P[u,v]\cup P[x^{+},z]\cup P[z^{+},w]\cup P[w^{+},y^{+}]+\{(w,w^{+}),(x,x^{+}),(y,y^{+}),(z,z^{+})\}-(x,y)$ is a H-path of $BH_{n}-F$ passing through $L$.

Suppose now that $f\in E(P[u,v])$. Let $(s,t)=f$. Without loss of generality, assume that $s\in X$ and $t\in Y$. Let $g=s^-$ (resp. $h=t^-$), if $s=x$ (resp. $t=y$); and $g=s^+$ (resp. $h=t^+$), otherwise. Then $g\neq x^+$ (resp. $h\neq y^+$). By Theorem \ref{th-cheng2014}, there exist two vertex-disjoint paths $P[x^+,z]$ and $P[g,c]$ (resp. $P[c^+,d]$ and $P[z^+,w]$) in $B^{1}$ (resp. $B^2$) such that each vertex of $B^{1}$ (resp. $B^2$) lies on one of the two paths. Theorem \ref{th-cheng2014} implies that there exist two vertex-disjoint paths $P[h,d^+]$ and $P[y^+,w^+]$ in $B^{3}$ such that each vertex of $B^{3}$ lies on one of the two paths. Thus, $P[u,v]\cup P[x^{+},z]\cup P[g,c]\cup P[z^{+},w]\cup P[c^{+},d]\cup P[y^{+},w^{+}]\cup P[h,d^{+}]+\{(c,c^{+}),(d,d^{+}),(s,g),(t,h),(w,w^{+}),(x,x^{+}),(y,y^{+}),(z,z^{+})\}-\{(x,y),(s,t)\}$ is a H-path of $BH_{n}-F$ passing through $L$.

{\it Case 1.1.2.}  $i=1$.

Let $(x,y)\in E(C_{0})\setminus E(L_{0})$. By Theorem \ref{th-xu2007}, $B^{1}$ has a H-path $P[u,v]$. Let $(x^{+},z)\in E(P[u,v])$. Let $w\in V_2\cap X$. By Theorem \ref{th-xu2007}, $B^{2}$, $B^3$ have H-paths $P[z^{+},w]$ and $P[w^{+},y^{+}]$, respectively. Hence, $C_{0}\cup P[u,v]\cup P[z^{+},w]\cup P[w^{+},y^{+}]+\{(w,w^{+}),(x,x^{+}),(y,y^{+}),(z,z^{+})\}-\{(x,y),(x^{+},z)\}$ is a H-path of $BH_{n}-F$ passing through $L$.

{\it Case 1.1.3.}  $i=2$.

Let $(x,y)\in E(C_{0})\setminus E(L_{0})$ and let $z\in V_1\cap X$. By Theorem \ref{th-xu2007}, $B^{1}$, $B^2$ have H-paths $P[x^{+},z]$ and $P[u,v]$, respectively. Let $(z^+,w)\in E(P[u,v])$. By Theorem \ref{th-xu2007}, $B^{3}$ has a H-path $P[y^+,w^+]$. Thus, $C_{0}\cup P[x^{+},z]\cup P[u,v]\cup P[y^+,w^+]+\{(w,w^{+}),(x,x^{+}),(y,y^{+}),(z,z^{+})\}-\{(x,y),(z^+,w)\}$ is a H-path of $BH_{n}-F$ passing through $L$.

{\it Case 1.1.4.}  $i=3$.

Let $(x,y)\in E(C_{0})\setminus E(L_{0})$ and let $z\in V_1\cap X, w\in V_2\cap X$. By Theorem \ref{th-xu2007}, $B^{1}$, $B^2$ have H-paths $P[x^{+},z]$ and $P[z^+,w]$, respectively. There is a neighbor of $y$ in $B^3$, say $y^+$, being not $u$, and there is a neighbor of $w$ in $B^3$, say $w^+$, being not $v$. By Theorem \ref{th-cheng2014}, there exist two vertex-disjoint paths $P[y^+,v]$ and $P[u,w^+]$ in $B^{3}$ such that each vertex of $B^{3}$ lies on one of the two paths. Thus, $C_{0}\cup P[x^{+},z]\cup P[z^+,w]\cup P[y^+,v]\cup P[u,w^+]+\{(w,w^{+}),(x,x^{+}),(y,y^{+}),(z,z^{+})\}-(x,y)$ is a H-path of $BH_{n}-F$ passing through $L$.

{\it Case 1.2.}   $l=1$ or $l=2$.

The proofs of the cases $l=1$ and $l=2$ are analogous. We here consider the case $l=1$.
Let $(a,b)\in E(C_{0})\setminus E(L_{0})$. Without loss of generality, assume that $a\in X$ and $b\in Y$.

{\it Case 1.2.1.}  $i=0$.

Since $F_{0}=F\neq \emptyset$, let $f\in F_{0}$. By the induction hypothesis, $B^{0}-F_{0}\setminus \{f\}$ has a H-path $P[u,v]$ passing through $L_{0}$. Let $(s,t)=f$, if $f$ lies on $P[u,v]$; and let $(s,t)\in E(P[u,v])\setminus E(L_{0})$, otherwise. Without loss of generality, assume that $s\in X$ and $t\in Y$. Let $y\in V_2\cap X$. By Theorem \ref{th-xu2007}, $B^{1}$, $B^{2}$, $B^{3}$ have H-paths $P[s^{+},x]$, $P[x^{+},y]$ and $P[y^{+},t^{+}]$, respectively. Thus, $P[u,v]\cup P[s^{+},x]\cup P[x^+,y]\cup P[y^+,t^+]+\{(s,s^{+}),(t,t^+),(x,x^{+}),(y,y^{+})\}-(s,t)$ is a H-path of $BH_{n}-F$ passing through $L$.

{\it Case 1.2.2.}  $i=1$.

There is a neighbor of $a$ in $B^1$, say $a^+$, being not $v$. Let $c\in V_2\cap X$. By Theorem \ref{th-xu2007}, $B^{2}$, $B^{3}$ have H-paths $P[x^{+},c]$, $P[b^{+},c^+]$, respectively.

Suppose first that $u\neq x$. By Theorem \ref{th-cheng2014}, there exist two vertex-disjoint paths $P[a^+,u]$ and $P[v,x]$ in $B^{1}$ such that each vertex of $B^{1}$ lies on one of the two paths. Thus, $C_0\cup P[a^{+},u]\cup P[v,x]\cup P[x^+,c]\cup P[b^+,c^+]+\{(a,a^{+}),(b,b^+),(c,c^+),(x,x^{+})\}-(a,b)$ is a H-path of $BH_{n}-F$ passing through $L$.

Suppose now that $u=x$. By Theorem \ref{th-lv2014}, $B^1-\{u\}$ has a H-path $P[a^+,v]$. Thus, $C_0\cup P[a^{+},v]\cup P[x^+,c]\cup P[b^+,c^+]+\{(a,a^{+}),(b,b^+),(c,c^+),(u,x^{+})\}-(a,b)$ is a H-path of $BH_{n}-F$ passing through $L$.

{\it Case 1.2.3.}  $i=2$.

By Theorem \ref{th-xu2007}, $B^{1}$, $B^{2}$ have H-paths $P[a^{+},x]$ and $P[u,v]$, respectively. Let $(x^+,c)\in E(P[u,v])$. By Theorem \ref{th-xu2007}, $B^{3}$ has a H-path $P[b^+,c^+]$. Thus, $C_0\cup P[a^{+},x]\cup P[u,v]\cup P[b^+,c^+]+\{(a,a^{+}),(b,b^+),(c,c^+),(x,x^{+})\}-\{(a,b),(x^+,c)\}$ is a H-path of $BH_{n}-F$ passing through $L$.

{\it Case 1.2.4.}  $i=3$.

By Theorem \ref{th-xu2007}, $B^{1}$, $B^{3}$ have H-paths $P[a^{+},x]$ and $P[u,v]$, respectively. Let $(b^+,c)\in E(P[u,v])$. By Theorem \ref{th-xu2007}, $B^{3}$ has a H-path $P[x^+,c^+]$. Thus, $C_0\cup P[a^{+},x]\cup P[x^+,c^+]\cup P[u,v]+\{(a,a^{+}),(b,b^+),(c,c^+),(x,x^{+})\}-\{(a,b),(b^+,c)\}$ is a H-path of $BH_{n}-F$ passing through $L$.

{\it Case 2.}  $u\in V_i$, $v\in V_j$, for $i,j\in N_{4}$ and $i\neq j$.

{\it Case 2.1.}  $l=0$ or $l=3$.

The proofs of the cases $l=0$ and $l=3$ are analogous. We here consider the case $l=0$.

{\it Case 2.1.1.}  $i=0$.

Let $(u,a)\in E(C_{0})\setminus E(L_{0})$. In this case, $P[u,a]=C_0-(u,a)$ is a H-path passing through $L_0$ of $B^0-F_0$. Let $z,b\in V_1\cap X$, $c,w\in V_2\cap X$ be pair-wires distinct.

Suppose first that $j=1$.

If $x^+\neq v$,
let $(x,y)\in E(P[u,a])\setminus E(L_0)$. Let $g=a^-$, if $y=a$; and $g=a^+$, otherwise. Then $g\neq y^+$. By Theorem \ref{th-cheng2014}, there exist two vertex-disjoint paths $P[x^+,z]$ and $P[v,b]$ (resp. $P[z^+,w]$ and $P[b^+,c]$) in $B^{1}$ (resp. $B^2$) such that each vertex of $B^{1}$ (resp. $B^2$) lies on one of the two paths. Theorem \ref{th-cheng2014} implies that there exist two vertex-disjoint paths $P[g,c^+]$ and $P[y^+,w^+]$ in $B^{3}$ such that each vertex of $B^{3}$ lies on one of the two paths. Thus, $P[u,a]\cup P[x^{+},z]\cup P[v,b]\cup P[z^{+},w]\cup P[b^{+},c]\cup P[y^{+},w^{+}]\cup P[g,c^{+}]+\{(a,g),(b,b^{+}),(c,c^{+}),(w,w^{+}),(x,x^{+}),(y,y^{+}),(z,z^{+})\}-(x,y)$ is a H-path of $BH_{n}-F$ passing through $L$.

If $x^+=v$ and $x$ is incident with none of $E(L_0)$, then $u\neq x$. Let $y$ be the neighbor of $x$ on the segment of $P[u,a]$ between $x$ and $u$. By Theorem \ref{th-lv2014}, $B^1-\{v\}$ has a H-path $P[z,b]$. By Theorem \ref{th-cheng2014}, there exist two vertex-disjoint paths $P[z^+,w]$ and $P[b^+,c]$ (resp. $P[a^+,c^+]$ and $P[y^+,w^+]$) in $B^{2}$ (resp. $B^3$) such that each vertex of $B^{2}$ (resp. $B^3$) lies on one of the two paths. Thus, $P[u,a]\cup P[z,b]\cup P[z^{+},w]\cup P[b^{+},c]\cup P[y^{+},w^{+}]\cup P[a^{+},c^{+}]+\{(a,a^{+}),(b,b^{+}),(c,c^{+}),(w,w^{+}),(x,v),(y,y^{+}),(z,z^{+})\}-(x,y)$ is a H-path of $BH_{n}-F$ passing through $L$.

If $x^+=v$ and $L_0$ has a maximal path $P[x,r]$ with $r\neq x$, in this case, $r\neq u$. Let $(x,s)\in E(P[x,r])$. Note that $\{u,s\}$ is compatible to $L_0-(x,s)$. By the induction hypothesis, $B^0-F_0$ has a H-path $P[u,s]$ passing through $L_0-(x,s)$. Let $y,t$ be the two distinct neighbors of $x$ on $P[u,s]$. By Theorem \ref{th-lv2014}, $B^1-\{v\}$ has a H-path $P[z,b]$. By Theorem \ref{th-cheng2014}, there exist two vertex-disjoint paths $P[z^+,w]$ and $P[b^+,c]$ (resp. $P[t^+,c^+]$ and $P[y^+,w^+]$) in $B^{2}$ (resp. $B^3$) such that each vertex of $B^{2}$ (resp. $B^3$) lies on one of the two paths. Thus, $P[u,s]\cup P[z,b]\cup P[z^{+},w]\cup P[b^{+},c]\cup P[y^{+},w^{+}]\cup P[t^{+},c^{+}]+\{(x,s),(b,b^{+}),(c,c^{+}),(t,t^+),(w,w^{+}),(x,v),(y,y^{+}),$ $(z,z^{+})\}-\{(x,y),(x,t)\}$ is a H-path of $BH_{n}-F$ passing through $L$.

Suppose second that $j=2$. Let $(x,y)\in E(P[u,a])\setminus E(L_0)$. Let $g=a^-$, if $y=a$; and $g=a^+$, otherwise. Then $g\neq y^+$. By Theorem \ref{th-xu2007}, $B^1$ has a H-path $P[x^+,z]$. There is a neighbor of $z$ in $B^2$, say $z^+$, being not $v$. By Theorem \ref{th-cheng2014}, there exist two vertex-disjoint paths $P[z^+,w]$ and $P[v,c]$ (resp. $P[g,c^+]$ and $P[y^+,w^+]$) in $B^{2}$ (resp. $B^3$) such that each vertex of $B^{2}$ (resp. $B^3$) lies on one of the two paths. Thus, $P[u,a]\cup P[x^+,z]\cup P[z^{+},w]\cup P[v,c]\cup P[y^{+},w^{+}]\cup P[g,c^{+}]+\{(a,g),(c,c^{+}),(w,w^{+}),(x,x^+),(y,$ $y^{+}),(z,z^{+})\}-(x,y)$ is a H-path of $BH_{n}-F$ passing through $L$.


Suppose now that $j=3$. Let $(x,y)\in E(P[u,a])\setminus E(L_0)$. Let $g=a^-$, if $y=a$; and $g=a^+$, otherwise. Then $g\neq y^+$. By Theorem \ref{th-xu2007}, $B^1$, $B^2$ have H-paths $P[x^+,z]$ and $P[z^+,w]$, respectively. There is a neighbor of $w$ in $B^3$, say $w^+$, being not $v$. By Theorem \ref{th-cheng2014}, there exist two vertex-disjoint paths $P[g,v]$ and $P[y^+,w^+]$ $B^3$ such that each vertex of $B^3$ lies on one of the two paths. Thus, $P[u,a]\cup P[x^+,z]\cup P[z^{+},w]\cup P[y^{+},w^{+}]\cup P[g,v]+\{(a,g),(w,w^{+}),(x,x^+),(y,y^{+}),(z,z^{+})\}-(x,y)$ is a H-path of $BH_{n}-F$ passing through $L$.

{\it Case 2.1.2.}  $i=1,j=2$.

Let $(x,y)\in E(C_0)\setminus E(L_0)$ and let $z\in V_2\cap X$. By Theorem \ref{th-xu2007}, $B^{1}$, $B^{2}$, $B^{3}$ have H-paths $P[u,x^{+}]$, $P[v,z]$ and $P[z^{+},y^{+}]$, respectively. Hence, $C_{0}\cup P[u,x^{+}]\cup P[v,z]\cup P[z^{+},y^{+}]+\{(x,x^{+}),(y,y^{+}),(z,z^{+})\}-(x,y)$ is a H-path of $BH_{n}-F$ passing through $L$.

{\it Case 2.1.3.}  $i=1,j=3$.

Let $(x,y)\in E(C_0)\setminus E(L_0)$. Then $P[x,y]=C_0-(x,y)$ is a H-path passing through $L_0$ of $B^0-F_0$. By Lemma \ref{le-3}, there is an edge $(a,b)\in E(P[x,y])\setminus E(L_0)$ for some $a\in X$ and $b\in Y$ such that $\{a,b\}\cap \{x,y\}=\emptyset$. By Theorem \ref{th-xu2007}, $B^{1}$, $B^{3}$ have H-paths $P[u,x^{+}]$, and $P[y^{+},v]$, respectively. Let $z$ be the neighbor of $a^+$ on the segment of $P[u,x^+]$ between $a^+$ and $x^+$, and let $c$ be the neighbor of $b^+$ on the segment of $P[y^+,v]$ between $b^+$ and $y^+$. By Theorem \ref{th-xu2007}, $B^{2}$ has a H-path $P[z^+,c^+]$. Thus, $P[x,y]\cup P[u,x^{+}]\cup P[z^{+},c^{+}]\cup P[y^+,v]+\{(a,a^+),(b,b^+),(c,c^+),(x,x^{+}),(y,y^{+}),(z,z^{+})\}-\{(a,b),(a^+,z),(b^+,c)\}$ is a H-path of $BH_{n}-F$ passing through $L$.

{\it Case 2.1.4.}  $i=2,j=3$.

Let $(x,y)\in E(C_0)\setminus E(L_0)$ and let $z\in V_1\cap X$. By Theorem \ref{th-xu2007}, $B^{1}$, $B^{2}$, $B^{3}$ have H-paths $P[x^{+},z]$, $P[u,z^+]$ and $P[y^{+},v]$, respectively. Hence, $C_{0}\cup P[x^{+},z]\cup P[u,z^+]\cup P[y^{+},v]+\{(x,x^{+}),(y,y^{+}),(z,z^{+})\}-(x,y)$ is a H-path of $BH_{n}-F$ passing through $L$.

{\it Case 2.2.}  $l=1$ or $l=2$.

The proofs of the cases $l=1$ and $l=2$ are analogous. We here consider the case $l=1$.

{\it Case 2.2.1.}  $i=0,j=1$.

Let $(u,a)\in E(C_{0})\setminus E(L_{0})$ and let $b\in V_2\cap X$. By Theorem \ref{th-xu2007}, $B^{1}$, $B^{2}$, $B^{3}$ have H-paths $P[v,x]$, $P[x^+,b]$ and $P[a^{+},b^+]$, respectively. Hence, $C_{0}\cup P[v,x]\cup P[x^+,b]\cup P[a^{+},b^+]+\{(a,a^{+}),(b,b^{+}),(x,x^{+})\}-(u,a)$ is a H-path of $BH_{n}-F$ passing through $L$.

{\it Case 2.2.2.}  $i=0,j=2$.

Let $(u,a)\in E(C_{0})\setminus E(L_{0})$. Then $P[u,a]=C_0-(u,a)$ is a H-path passing through $L_0$ of $B^0-F_0$. There are $\lfloor|E(P[u,a])|/2\rfloor=\lfloor4^{n-1}-1/2\rfloor$ edges each of which has the form $(s,t)$ with $s\in X$ and $t\in Y$ such that $t$ lies on the segment of $P[u,a]$ between $u$ and $s$. Since $\lfloor|E(P[u,a])|/2\rfloor-|E(L_0)|\geq \lfloor4^{n-1}-1/2\rfloor-(2n-4)>0$, there is at least such one edge $(s,t)$ on $P[u,a]$ that meats above requirements and furthermore $(s,t)\notin E(L_0)$. Let $b\in V_2\cap X$. By Theorem \ref{th-xu2007}, $B^{1}$, $B^2$ have H-paths $P[s^+,x]$ and $P[v,b]$, respectively. Let $(x^+,y)\in E(P[v,b])$. By Theorem \ref{th-cheng2014}, there are two vertex-disjoint paths $P[a^+,y^+]$ and $P[t^+,b^+]$ in $B^{3}$ each vertex of $B^3$ lies on one of the two paths. Thus, $P[u,a]\cup P[s^+,x]\cup P[v,b]\cup P[a^{+},y^{+}]\cup P[t^{+},b^{+}]+\{(a,a^+),(b,b^{+}),(s,s^{+}),(t,t^+),(x,x^+),(y,y^{+})\}-\{(s,t),(x^+,y)\}$ is a H-path of $BH_{n}-F$ passing through $L$.

{\it Case 2.2.3.}  $i=0,j=3$.

According to Case 2.2.2. There is an edge $(s,t)\in E(P[u,a])\setminus E(L_0)$ for some $s\in X$ and $t\in Y$ such that $t$ lies on the segment of $P[u,a]$ between $u$ and $s$. Let $b\in V_2\cap X$. By Theorem \ref{th-xu2007}, $B^{1}$, $B^2$ have H-paths $P[s^+,x]$ and $P[x^+,b]$, respectively. There is a neighbor of $b$ in $B^3$, say $b^+$, being not $v$. By Theorem \ref{th-cheng2014}, there are two vertex-disjoint paths $P[a^+,v]$ and $P[t^+,b^+]$ in $B^{3}$ each vertex of $B^3$ lies on one of the two paths. Thus, $P[u,a]\cup P[s^+,x]\cup P[x^+,b]\cup P[a^{+},v]\cup P[t^{+},b^{+}]+\{(a,a^+),(b,b^{+}),(s,s^{+}),(t,t^+),(x,x^+)\}-(s,t)$ is a H-path of $BH_{n}-F$ passing through $L$.

{\it Case 2.2.4.}  $i=1,j=2$.

Let $(a,b)\in E(C_{0})\setminus E(L_{0})$. Then $P[a,b]=C_0-(a,b)$ is a H-path passing through $L_0$ of $B^0-F_0$. There are $\lfloor|E(P[a,b])|/2\rfloor=\lfloor4^{n-1}-1/2\rfloor$ edges each of which has the form $(s,t)$ with $s\in X$ and $t\in Y$ such that $t$ lies on the segment of $P[a,b]$ between $a$ and $s$. Since $\lfloor|E(P[a,b])|/2\rfloor-|E(L_0)|\geq \lfloor4^{n-1}-1/2\rfloor-(2n-4)>0$, there is at least such one edge $(s,t)$ on $P[a,b]$ that meats above requirements and furthermore $(s,t)\notin E(L_0)$. Let $c\in V_2\cap X$. By Theorem \ref{th-xu2007}, $B^2$ has a H-path $P[v,c]$.

Suppose first that $u\neq x$. Let $(x^+,y)\in E(P[v,c])$. By Theorem \ref{th-cheng2014}, there are two vertex-disjoint paths $P[a^+,u]$ and $P[s^+,x]$ (resp. $P[b^+,y^+]$ and $P[t^+,c^+]$) in $B^{1}$ (resp. $B^3$) each vertex of $B^{1}$ (resp. $B^3$) lies on one of the two paths. Thus, $P[a,b]\cup P[s^+,x]\cup P[a^+,u]\cup P[v,c]\cup P[t^{+},c^{+}]\cup P[b^+,y^+]+\{(a,a^+),(b,b^{+}),(c,c^+),(s,s^{+}),(t,t^+),(x,x^+),(y,y^{+})\}-\{(s,t),(x^+,y)\}$ is a H-path of $BH_{n}-F$ passing through $L$.

Suppose now that $u=x$. In this case, $v\neq x^+$. Let $y$ be the neighbor of $x^+$ on the segment of $P[v,c]$ between $x^+$ and $v$. By Theorem \ref{th-cheng2014}, there are two vertex-disjoint paths $P[b^+,y^+]$ and $P[t^+,c^+]$ $B^3$ each vertex of $B^3$ lies on one of the two paths. By Theorem \ref{th-lv2014}, $B^1-\{u\}$ has a H-path $P[a^+,s^+]$. Thus, $P[a,b]\cup P[a^+,s^+]\cup P[v,c]\cup P[t^{+},c^{+}]\cup P[b^+,y^+]+\{(a,a^+),(b,b^{+}),(c,c^+),(s,s^{+}),(t,t^+),(u,x^+),(y,y^{+})\}-\{(s,t),(x^+,y)\}$ is a H-path of $BH_{n}-F$ passing through $L$.

{\it Case 2.2.5.}  $i=1,j=3$.

According to Case 2.2.4, there is an edge $(s,t)\in E(P[a,b])\setminus E(L_0)$ for some $s\in X$ and $t\in Y$ such that $t$ lies on the segment of $P[a,b]$ between $a$ and $s$. Let $c\in V_2\cap X$. By Theorem \ref{th-cheng2014}, there are two vertex-disjoint paths $P[a^+,u]$ and $P[s^+,x]$ (resp. $P[b^+,v]$ and $P[t^+,c^+]$) in $B^{1}$ (resp. $B^3$) each vertex of $B^{1}$ (resp. $B^3$) lies on one of the two paths. By Theorem \ref{th-xu2007}, $B^2$ has a H-path $P[x^+,c]$. Thus, $P[a,b]\cup P[s^+,x]\cup P[a^+,u]\cup P[x^+,c]\cup P[t^{+},c^{+}]\cup P[b^+,v]+\{(a,a^+),(b,b^{+}),(c,c^+),(s,s^{+}),(t,t^+),(x,x^+)\}-(s,t)$ is a H-path of $BH_{n}-F$ passing through $L$.

{\it Case 2.2.6.}  $i=2,j=3$.

Let $(a,b)\in E(C_{0})\setminus E(L_{0})$. By Theorem \ref{th-xu2007}, $B^{1}$, $B^2$, $B^3$ have H-paths $P[a^{+},x]$, $P[x^{+},u]$ and $P[b^{+},v]$, respectively. Thus, $C_{0}\cup P[a^{+},x]\cup P[x^{+},u]\cup P[b^{+},v]+\{(x,x^{+}),(a,a^{+}),(b,b^{+})\}-(a,b)$ is a H-path of $BH_{n}-F$ passing through $L$.
\end{proof}

\section{$|F^c|=1$, $L^c=\emptyset$}
\label{section5}

In this section, let $(s,s^+)$ be the edge of $F^c$ for some $s\in X$ and $s^+\in Y$.

\begin{lemma}
If $|E(L_{0})\cup F_{0}|\leq 2n-4$, then $BH_{n}-F$ contains a H-path $P[u,v]$ passing through $L$.
\end{lemma}

\begin{proof}
In this  scenario, $|E(L_{k})\cup F_{k}|\leq 2n-5$ for $k\in N_{4}\setminus \{0\}$. 

{\it Case 1.}  $u,v\in V_i$.


{\it Case 1.1.}  $i=0$.

By the induction hypothesis, $B^{0}-F_{0}$ has a H-path $P[u,v]$ passing through $L_{0}$. By Lemma \ref{le-3}, there is an edge $(a,b)\in E(P[u,v])\setminus E(L_0)$ for some $a\in X$ and $b\in Y$ such that $\{a,b\}\cap \{u,v\}=\emptyset$, $\{a,b\}\cap \{s,s^+\}=\emptyset$, $a^+$ or $a^+$ (resp. $b^+$ or $b^-$), say $a^+$ (resp. $b^+$), is incident with none of $E(L_1)$ (resp. $E(L_3)$). By Lemma \ref{le-2}, there are vertices $c\in V_1\cap X$ and $d\in V_2\cap X$ such that $c$ (resp. $c^+$) is incident with none of $E(L_1)$ (resp. $E(L_2)$), $d$ (resp. $d^+$) is incident with none of $E(L_2)$ (resp. $E(L_3)$) and $s\notin \{c,d\}$. By the induction hypothesis, $B^1-F_1$, $B^2-F_2$, $B^3-F_3$ have H-paths $P[a^+,c]$, $P[c^+,d]$, $P[b^+,d^+]$ passing through $L_1$, $L_2$ and $L_3$, respectively. Thus, $P[u,v]\cup P[a^+,c]\cup P[c^+,d]\cup P[b^+,d^+]+\{(a,a^+),(b,b^+),(c,c^+),(d,d^+)\}-(a,b)$ is a H-path of $BH_{n}-F$ passing through $L$.


{\it Case 1.2.}  $i=1$.

By the induction hypothesis, $B^{1}-F_{1}$ has a H-path $P[u,v]$ passing through $L_{1}$. By Lemma \ref{le-3}, there is an edge $(a,b)\in E(P[u,v])\setminus E(L_1)$ for some $a\in X$ and $b\in Y$ such that $\{a,b\}\cap \{u,v\}=\emptyset$, $\{a,b\}\cap \{s,s^+\}=\emptyset$, $a^+$ or $a^+$ (resp. $b^+$ or $b^-$), say $a^+$ (resp. $b^+$), is incident with none of $E(L_2)$ (resp. $E(L_0)$). By Lemma \ref{le-2}, there are vertices $c\in V_2\cap X$ and $d\in V_0\cap Y$ such that $c$ (resp. $c^+$) is incident with none of $E(L_2)$ (resp. $E(L_3)$), $d$ (resp. $d^+$) is incident with none of $E(L_0)$ (resp. $E(L_3)$) and $\{c,d\}\cap \{s,s^+\}=\emptyset$. By the induction hypothesis, $B^0-F_0$, $B^2-F_2$, $B^3-F_3$ have H-paths $P[b^+,d]$, $P[a^+,c]$, $P[c^+,d^+]$ passing through $L_0$, $L_2$ and $L_3$, respectively. Thus, $P[b^+,d]\cup P[u,v]\cup P[a^+,c]\cup P[c^+,d^+]+\{(a,a^+),(b,b^+),(c,c^+),(d,d^+)\}-(a,b)$ is a H-path of $BH_{n}-F$ passing through $L$.

{\it Case 1.3.}  $i=2$.

By the induction hypothesis, $B^{2}-F_{2}$ has a H-path $P[u,v]$ passing through $L_{2}$. By Lemma \ref{le-3}, there is an edge $(a,b)\in E(P[u,v])\setminus E(L_2)$ for some $a\in X$ and $b\in Y$ such that $\{a,b\}\cap \{u,v\}=\emptyset$, $\{a,b\}\cap \{s,s^+\}=\emptyset$, $a^+$ or $a^+$ (resp. $b^+$ or $b^-$), say $a^+$ (resp. $b^+$), is incident with none of $E(L_3)$ (resp. $E(L_1)$). By Lemma \ref{le-2}, there are vertices $c\in V_0\cap Y$ and $d\in V_0\cap X$ such that $c$ (resp. $c^+$) is incident with none of $E(L_0)$ (resp. $E(L_3)$), $d$ (resp. $d^+$) is incident with none of $E(L_0)$ (resp. $E(L_1)$) and $\{c,d\}\cap \{s,s^+\}=\emptyset$. By the induction hypothesis, $B^0-F_0$, $B^1-F_1$, $B^3-F_3$ have H-paths $P[d,c]$, $P[d^+,b^+]$, $P[a^+,c^+]$ passing through $L_0$, $L_1$ and $L_3$, respectively. Thus, $P[d,c]\cup P[d^+,b^+]\cup P[u,v]\cup P[a^+,c^+]+\{(a,a^+),(b,b^+),(c,c^+),(d,d^+)\}-(a,b)$ is a H-path of $BH_{n}-F$ passing through $L$.


{\it Case 1.4.}  $i=3$.

By the induction hypothesis, $B^{3}-F_{3}$ has a H-path $P[u,v]$ passing through $L_{3}$. By Lemma \ref{le-3}, there is an edge $(a,b)\in E(P[u,v])\setminus E(L_3)$ for some $a\in X$ and $b\in Y$ such that $\{a,b\}\cap \{u,v\}=\emptyset$, $\{a,b\}\cap \{s,s^+\}=\emptyset$, $a^+$ or $a^+$ (resp. $b^+$ or $b^-$), say $a^+$ (resp. $b^+$), is incident with none of $E(L_0)$ (resp. $E(L_2)$). By Lemma \ref{le-2}, there are vertices $c\in V_0\cap X$ and $d\in V_1\cap X$ such that $c$ (resp. $c^+$) is incident with none of $E(L_0)$ (resp. $E(L_1)$), $d$ (resp. $d^+$) is incident with none of $E(L_1)$ (resp. $E(L_2)$) and $s\notin \{c,d\}$. By the induction hypothesis, $B^0-F_0$, $B^1-F_1$, $B^2-F_2$ have H-paths $P[a^+,c]$, $P[c^+,d]$, $P[b^+,d^+]$ passing through $L_0$, $L_1$ and $L_2$, respectively. Thus, $P[a^+,c]\cup P[c^+,d]\cup P[b^+,d^+]\cup P[u,v]+\{(a,a^+),(b,b^+),(c,c^+),(d,d^+)\}-(a,b)$ is a H-path of $BH_{n}-F$ passing through $L$.




{\it Case 2.}  $u\in V_i$ and $v\in V_j$ for $i,j\in N_{4}$ and $i\neq j$.


{\it Case 2.1.}  $i=0$.

For $n=3$, $|E(L_0)\cup F_0|\leq 2n-6\leq 0$. By Lemma \ref{le-4}, there is an $a\in V_0\cap Y$ such that $a\neq s^+$, $a$ and $a^{\pm}$ are incident with none of $E(L_0)$ and $E(L_3)\cup F_3$, respectively, and $s^+$ is not adjacent to $a^{\pm}$.
By the induction hypothesis, $B^{0}-F_{0}$ has a H-path $P[u,a]$ passing through $L_{0}$.

{\it Case 2.1.1.}  $j=1$.

By Lemma \ref{le-2}, there are vertices $b\in V_1\cap X$ and $c\in V_2\cap X$ such that $b$ (resp. $b^+$) is incident with none of $E(L_1)$ (resp. $E(L_2)$), $c$ (resp. $c^+$) is incident with none of $E(L_2)$ (resp. $E(L_3)$) and $s\notin \{b,c\}$. By the induction hypothesis, $B^1-F_1$, $B^2-F_2$, $B^3-F_3$ have H-paths $P[v,b]$, $P[b^+,c]$, $P[a^+,c^+]$ passing through $L_1$, $L_2$ and $L_3$, respectively. Thus, $P[u,a]\cup P[v,b]\cup P[b^+,c]\cup P[a^+,c^+]+\{(a,a^+),(b,b^+),(c,c^+)\}$ is a H-path of $BH_{n}-F$ passing through $L$.


{\it Case 2.1.2.}  $j=2$.

In this scenario, $|E(L_k)\cup F_k|\leq 2n-5$ for $k\in N_4\setminus \{0\}$. By Lemma \ref{le-3}, there is an edge $(x,y)\in E(P[u,a])\setminus E(L_{0})$ for some $x\in X$ and $y\in Y$ such that $x^{+}$ or $x^{-}$ (resp. $y^{+}$ or $y^{-}$), say $x^+$ (resp. $y^+$), is incident with none of $E(L_1)$ (resp. $E(l_3)$), $\{x,y\}\cap \{u,a\}=\emptyset$ and $\{x,y\}\cap \{s,s^+\}=\emptyset$. Let $g=a^-$, if $y=a$; and $g=a^+$, otherwise. Then $g\neq y^+$.

Suppose first that $|E(L_3)\cup F_3|=2n-5$. In this case, $|E(L_1)\cup F_1|\leq \min\{\sum_{k\in N_4\setminus \{0\}}|E(L_k)\cup F_k|, |E(L_0)\cup F_0|\}\leq 1$. Then $|E(L_2)\cup F_2|\leq 1$. By Lemma \ref{le-2}, there is a $b\in V_3\cap Y$ such that $b\notin \{s,s^+\}$, $b$ and $b^{\pm}$ are incident with none of $E(L_3)$ and $E(L_2)$, respectively. By the induction hypothesis, $B^3-F_3$ has a H-path $P[g,b]$ passing through $L_3$. Let $c$ be the neighbor of $y^+$ on the segment of $P[g,b]$ between $y^+$ and $g$. Since $|E(L_2)|\leq 1$, $c^+$ or $c^-$, say $c^+$, is not incident with none of $E(L_2)$. By Theorem \ref{th-yang2019}, $B^2-F_2$ has a H-path $P[v,b^+]$ passing through $L_2$. Let $z$ be the neighbor of $c^+$ on the segment of $P[v,b^+]$ between $c^+$ and $b^+$. By Theorem \ref{th-yang2019}, $B^1-F_1$ has a H-path $P[x^+,z^+]$ passing through $L_1$. Thus, $P[u,a]\cup P[x^+,z^+]\cup P[v,b^+]\cup P[g,b]+\{(a,g),(b,b^+),(c,c^+),(x,x^+),(y,y^+),(z,z^+)\}-\{(x,y),(c^+,z),(y^+,c)\}$ is a H-path of $BH_{n}-F$ passing through $L$.

Suppose second that $|E(L_3)\cup F_3|\leq 2n-6$ and $|E(L_m)\cup F_m|\leq 2n-6$ for $m\in \{1,2\}$. By Lemma \ref{le-9}, $g$ has two neighbors $c$ and $d$ in $B^3$ such that $c^+$ or $c^-$ (resp. $d^+$ or $d^-$), say $c^+$ (resp. $d^+$), is incident with none of $E(L_2)$, and $L_3+\{(g,c),(g,d)\}$ is a linear forest. By Lemma \ref{le-4}, there is a $b\in V_3\cap Y$ such that $b\neq s^+$, $b$ and $b^{\pm}$ are incident with none of $E(L_3)$ and $E(L_2)\cup F_2$, respectively, and $s^+$ is not adjacent to $b^{\pm}$. Note that $\{y^+,b\}$ is compatible to $L_3+\{(g,c),(g,d)\}$. By the induction hypothesis, $B^3-F_3$ has a H-path $P[y^+,b]$ passing through $L_3+\{(g,c),(g,d)\}$. Exactly one of $c$ and $d$, say $c$, lies on the segment of $P[y^+,b]$ between $g$ and $y^+$. By Lemma \ref{le-9}, $b^+$ has two neighbors $z$ and $w$ in $B^2$ such that $z^+$ or $z^-$ (resp. $w^+$ or $w^-$), say $z^+$ (resp. $w^+$), is incident with none of $E(L_1)$, and $L_2+\{(b^+,z),(b^+,w)\}$ is a linear forest. Note that $\{v,c^+\}$ is compatible to $L_2+\{(b^+,z),(b^+,w)\}$. By the induction hypothesis, $B^2-F_2$ has a H-path $P[v,c^+]$ passing through $L_2+\{(b^+,z),(b^+,w)\}$. Exactly one of $z$ and $w$, say $z$, lies on the segment of $P[c^+,v]$ between $b^+$ and $c^+$. By the induction hypothesis, $B^1-F_1$ has a H-path $P[x^+,z^+]$ passing through $L_1$. Thus, $P[u,a]\cup P[x^+,z^+]\cup P[v,c^+]\cup P[y^+,b]+\{(a,g),(b,b^+),(c,c^+),(x,x^+),(y,y^+),(z,z^+)\}-\{(x,y),(b^+,z),(g,c)\}$ is a H-path of $BH_{n}-F$ passing through $L$.

Suppose now that $|E(L_3)\cup F_3|\leq 2n-6$ and $|E(L_m)\cup F_m|=2n-5$ for some $m\in \{1,2\}$. If $n=3$, then $|E(L_3)\cup F_3|\leq 2n-6\leq 0$. If $n\geq 4$, then $|E(L_3)\cup F_3|\leq |E(L)\cup F|-|F^c|-|E(L_0)\cup F_0|-|E(L_m)\cup F_m|<0$. Thus, $E(L_3)\cup F_3=\emptyset$ for $n\geq 3$. By Lemma \ref{le-4}, there is a $z\in V_1\cap X\setminus \{s\}$ such that $z$ and $z^{\pm}$ are incident with none of $E(L_1)$ and $E(L_2)$, respectively, and $s$ is not adjacent to $z^{\pm}$. By the induction hypothesis, $B^1-F_1$ has a H-path $P[x^+,z]$ passing through $L_1$. There is a neighbor of $z$ in $B^2$, say $z^+$, being not $v$. By Lemma \ref{le-2}, there is a $b\in V_2\cap X\setminus \{s\}$ such that $b$ is incident with none of $E(L_2)$. By the induction hypothesis, $B^2-F_2$ has a H-path $P[v,b]$ passing through $L_2$. Let $c$ be the neighbor of $z^+$ on the segment of $P[v,b]$ between $z^+$ and $v$. By Theorem \ref{th-cheng2014}, there are two vertex-disjoint paths $P[g,c^+]$ and $P[y^+,b^+]$ in $B^3$ such that each vertex of $B^3$ lies on one of the two paths. Thus, $P[u,a]\cup P[x^+,z]\cup P[v,b]\cup P[y^+,b^+]\cup P[g,c^+]+\{(a,g),(b,b^+),(c,c^+),(x,x^+),(y,y^+),(z,z^+)\}-\{(x,y),(z^+,z)\}$ is a H-path of $BH_{n}-F$ passing through $L$.

{\it Case 2.1.3.}  $j=3$.

In this scenario, $|E(L_k)\cup F_k|\leq 2n-5$ for $k\in N_4\setminus \{0\}$. By Lemma \ref{le-3}, there is an edge $(x,y)\in E(P[u,a])\setminus E(L_{0})$ for some $x\in X$ and $y\in Y$ such that $x^{+}$ or $x^{-}$ (resp. $y^{+}$ or $y^{-}$), say $x^+$ (resp. $y^+$), is incident with none of $E(L_1)$ (resp. $E(l_3)$), $\{x,y\}\cap \{u,a\}=\emptyset$ and $\{x,y\}\cap \{s,s^+\}=\emptyset$. Let $g=a^-$, if $y=a$; and $g=a^+$, otherwise. Then $g\neq y^+$. By Lemma \ref{le-2}, there is a $z\in V_1\cap X\setminus \{s\}$ such that $z$ (resp. $z^+$) is incident with none of $E(L_1)$ (resp. $E(L_2)$). By the induction hypothesis, $B^1-F_1$ has a H-path $P[x^+,z]$ passing through $L_1$.

Suppose first that $|E(L_3)\cup F_3|=2n-5$. In this case, $|E(L_2)\cup F_2|\leq \min\{\sum_{k\in N_4\setminus \{0\}}|E(L_k)\cup F_k|, |E(L_0)\cup F_0|\}\leq 1$. By the induction hypothesis, $B^3-F_3$ has a H-path $P[y^+,v]$ passing through $L_3$. Let $c$ be the neighbor of $g$ on the segment of $P[y^+,v]$ between $y^+$ and $g$. Since $|E(L_2)|\leq 1$, $c^+$ or $c^-$, say $c^+$, is not incident with none of $E(L_2)$. By Theorem \ref{th-yang2019}, $B^2-F_2$ has a H-path $P[z^+,c^+]$ passing through $L_2$. Thus, $P[u,a]\cup P[x^+,z]\cup P[z^+,c^+]\cup P[y^+,v]+\{(a,g),(c,c^+),(x,x^+),(y,y^+),(z,z^+)\}-\{(x,y),(g,c)\}$ is a H-path of $BH_{n}-F$ passing through $L$.

Suppose second that $|E(L_3)\cup F_3|\leq 2n-6$ and $|E(L_m)\cup F_m|\leq 2n-6$ for $m\in \{1,2\}$. By Lemma \ref{le-9}, $g$ has two neighbors $c$ and $d$ in $B^3$ such that $c^+$ or $c^-$ (resp. $d^+$ or $d^-$), say $c^+$ (resp. $d^+$), is incident with none of $E(L_2)$, and $L_3+\{(g,c),(g,d)\}$ is a linear forest. Note that $\{y^+,v\}$ is compatible to $L_3+\{(g,c),(g,d)\}$. By the induction hypothesis, $B^3-F_3$ has a H-path $P[y^+,v]$ passing through $L_3+\{(g,c),(g,d)\}$. Exactly one of $c$ and $d$, say $c$, lies on the segment of $P[y^+,v]$ between $g$ and $y^+$. By the induction hypothesis, $B^2-F_2$ has a H-path $P[z^+,c^+]$ passing through $L_2$. Thus, $P[u,a]\cup P[x^+,z]\cup P[z^+,c^+]\cup P[y^+,v]+\{(a,g),(c,c^+),(x,x^+),(y,y^+),(z,z^+)\}-\{(x,y),(g,c)\}$ is a H-path of $BH_{n}-F$ passing through $L$.

Suppose now that $|E(L_3)\cup F_3|\leq 2n-6$ and $|E(L_m)\cup F_m|=2n-5$ for some $m\in \{1,2\}$. If $n=3$, then $|E(L_3)\cup F_3|\leq 2n-6\leq 0$. If $n\geq 4$, then $|E(L_3)\cup F_3|\leq |E(L)\cup F|-|F^c|-|E(L_0)\cup F_0|-|E(L_m)\cup F_m|<0$. Thus, $E(L_3)\cup F_3=\emptyset$ for $n\geq 3$. By Lemma \ref{le-2}, there is a $c\in V_2\cap X\setminus \{s\}$ such that $c$ is incident with none of $E(L_2)$. By the induction hypothesis, $B^2-F_2$ has a H-path $P[z^+,c]$ passing through $L_2$. There is a neighbor of $c$ in $B^3$, say $c^+$, being not $v$. By Theorem \ref{th-cheng2014}, there are two vertex-disjoint paths $P[g,v]$ and $P[y^+,c^+]$ in $B^3$ such that each vertex of $B^3$ lies on one of the two paths. Thus, $P[u,a]\cup P[x^+,z]\cup P[z^+,c]\cup P[y^+,c^+]\cup P[g,v]+\{(a,g),(c,c^+),(x,x^+),(y,y^+),(z,z^+)\}-(x,y)$ is a H-path of $BH_{n}-F$ passing through $L$.

{\it Case 2.2.}  $i\neq 0$.

By Lemma \ref{le-4}, there are vertices $a\in V_0\cap X\setminus \{s\}$ and $b\in V_0\cap Y\setminus \{s^+\}$ such that $a$ (resp. $a^+$) is incident with none of $E(L_0)$ (resp. $E(L_1)$), $b$ (resp. $b^+$) is incident with none of $E(L_0)$ (resp. $E(L_3)$), and $s$ (resp. $s^+$) is not adjacent to $a^+$ (resp. $b^+$). By the induction hypothesis, $B^0-F_0$ has a H-path $P[a,b]$ passing through $L_0$.

{\it Case 2.2.1.}  $i=1,j=2$.

By the induction hypothesis, $B^1-F_1$ has a H-path $P[a^+,u]$ passing through $L_1$. By Lemma \ref{le-2}, there is a $c\in V_3\cap Y\setminus \{s^+\}$ such that $c$ (resp. $c^+$) is incident with none of $E(L_3)$ (resp. $E(L_2)$). By the induction hypothesis, $B^2-F_2$, $B^3-F_3$ have H-paths $P[v,c^+]$, $P[b^+,c]$ passing through $L_2$ and $L_3$, respectively. Hence, $P[a,b]\cup P[a^{+},u]\cup P[v,c^{+}]\cup P[b^{+},c]+\{(a,a^{+}),(b,b^{+}),(c,c^{+})\}$ is a H-path of $BH_{n}-F$ passing through $L$.


{\it Case 2.2.2.}  $i=1,j=3$.

In this scenario, $|E(L_k)\cup F_k|\leq 2n-5$ for $k\in N_4\setminus \{0\}$. By Lemma \ref{le-3}, there is an edge $(x,y)\in E(P[a,b])\setminus E(L_{0})$ for some $x\in X$ and $y\in Y$ such that $x^{+}$ or $x^{-}$ (resp. $y^{+}$ or $y^{-}$), say $x^+$ (resp. $y^+$), is incident with none of $E(L_1)$ (resp. $E(l_3)$), $\{x,y\}\cap \{a,b\}=\emptyset$ and $\{x,y\}\cap \{s,s^+\}=\emptyset$. 

Suppose first that $|E(L_3)\cup F_3|=2n-5$. In this case, $|E(L_1)\cup F_1|\leq \min\{\sum_{k\in N_4\setminus \{0\}}|E(L_k)\cup F_k|, |E(L_0)\cup F_0|\}\leq 1$. Then $|E(L_2)\cup F_2|\leq 1$. By the induction hypothesis, $B^3-F_3$ has a H-path $P[y^+,v]$ passing through $L_3$. Let $c$ be the neighbor of $b^+$ on the segment of $P[y^+,v]$ between $y^+$ and $b^+$. Since $|E(L_2)|\leq 1$, $c^+$ or $c^-$, say $c^+$, is not incident with none of $E(L_2)$. By Theorem \ref{th-yang2019}, $B^1-F_1$ has a H-path $P[x^+,u]$ passing through $L_1$. Let $z$ be the neighbor of $a^+$ on the segment of $P[x^+,u]$ between $x^+$ and $a^+$. By Theorem \ref{th-yang2019}, $B^2-F_2$ has a H-path $P[z^+,c^+]$ passing through $L_2$. Thus, $P[a,b]\cup P[x^+,u]\cup P[z^+,c^+]\cup P[y^+,v]+\{(a,a^+),(b,b^+),(c,c^+),(x,x^+),(y,y^+),(z,z^+)\}-\{(x,y),(a^+,z),(b^+,c)\}$ is a H-path of $BH_{n}-F$ passing through $L$.

Suppose second that $|E(L_3)\cup F_3|\leq 2n-6$ and $|E(L_m)\cup F_m|\leq 2n-6$ for $m\in \{1,2\}$. By Lemma \ref{le-9}, $b^+$ has two neighbors $c$ and $d$ in $B^3$ such that $c^+$ or $c^-$ (resp. $d^+$ or $d^-$), say $c^+$ (resp. $d^+$), is incident with none of $E(L_2)$, and $L_3+\{(b^+,c),(b^+,d)\}$ is a linear forest. Note that $\{y^+,v\}$ is compatible to $L_3+\{(b^+,c),(b^+,d)\}$. By the induction hypothesis, $B^3-F_3$ has a H-path $P[y^+,v]$ passing through $L_3+\{(b^+,c),(b^+,d)\}$. Exactly one of $c$ and $d$, say $c$, lies on the segment of $P[y^+,v]$ between $b^+$ and $y^+$. By Lemma \ref{le-9}, $a^+$ has two neighbors $z$ and $w$ in $B^1$ such that $z^+$ or $z^-$ (resp. $w^+$ or $w^-$), say $z^+$ (resp. $w^+$), is incident with none of $E(L_2)$, and $L_1+\{(a^+,z),(a^+,w)\}$ is a linear forest. Note that $\{x^+,u\}$ is compatible to $L_1+\{(a^+,z),(a^+,w)\}$. By the induction hypothesis, $B^1-F_1$ has a H-path $P[x^+,u]$ passing through $L_1+\{(a^+,z),(a^+,w)\}$. Exactly one of $z$ and $w$, say $z$, lies on the segment of $P[x^+,u]$ between $a^+$ and $x^+$.
By the induction hypothesis, $B^2-F_2$ has a H-path $P[z^+,c^+]$ passing through $L_2$. Thus, $P[a,b]\cup P[x^+,u]\cup P[z^+,c^+]\cup P[y^+,v]+\{(a,a^+),(b,b^+),(c,c^+),(x,x^+),(y,y^+),(z,z^+)\}-\{(x,y),(a^+,z),(b^+,c)\}$ is a H-path of $BH_{n}-F$ passing through $L$.

Suppose third that $|E(L_3)\cup F_3|\leq 2n-6$ and $|E(L_1)\cup F_1|=2n-5$. If $n=3$, then $|E(L_3)\cup F_3|\leq 2n-6\leq 0$. If $n\geq 4$, then $|E(L_3)\cup F_3|\leq |E(L)\cup F|-|F^c|-|E(L_0)\cup F_0|-|E(L_1)\cup F_1|<0$. Thus, $E(L_3)\cup F_3=\emptyset$ for $n\geq 3$. Then $E(L_2)\cup F_2=\emptyset$. By the induction hypothesis, $B^1-F_1$ has a H-path $P[x^+,u]$ passing through $L_1$. Let $z$ be the neighbor of $a^+$ on the segment of $P[x^+,u]$ between $a^+$ and $x^+$. Let $c\in V_2\cap X\setminus \{s\}$. By Theorem \ref{th-xu2007}, $B^2$ has a H-path $P[z^+,c]$.
There is a neighbor of $c$ in $B^3$, say $c^+$, being not $v$. By Theorem \ref{th-cheng2014}, there are two vertex-disjoint paths $P[b^+,v]$ and $P[y^+,c^+]$ in $B^3$ such that each vertex of $B^3$ lies on one of the two paths. Thus, $P[a,b]\cup P[x^+,u]\cup P[z^+,c]\cup P[y^+,c^+]\cup P[b^+,v]+\{(a,a^+),(b,b^+),(c,c^+),(x,x^+),(y,y^+),(z,z^+)\}-\{(x,y),(a^+,z)\}$ is a H-path of $BH_{n}-F$ passing through $L$.

Suppose now that $|E(L_3)\cup F_3|\leq 2n-6$ and $|E(L_2)\cup F_2|=2n-5$. In this scenario, $E(L_3)\cup F_3=E(L_1)\cup F_1=\emptyset$. By Lemma \ref{le-2}, there are vertices $z\in V_2\cap Y\setminus \{s^+\}$ and $c\in V_2\cap X\setminus \{s\}$ such that $z$ and $c$ are incident with none of $E(L_2)$. By the induction hypothesis, $B^2-F_2$ has a H-path $P[z,c]$ passing through $L_2$. There is a neighbor of $z$ (resp. $c$) in $B^1$ (resp. $B^3$), say $z^+$ (resp. $c^+$), being not $u$ (resp. $v$). By Theorem \ref{th-cheng2014}, there are two vertex-disjoint paths $P[a^+,u]$ and $P[x^+,z^+]$ (resp. $P[b^+,v]$ and $P[y^+,c^+]$) in $B^1$ (resp. $B^3$) such that each vertex of $B^1$ (resp. $B^3$) lies on one of the two paths. Thus, $P[a,b]\cup P[a^+,u]\cup P[x^+,z^+]\cup P[z,c]\cup P[y^+,c^+]\cup P[b^+,v]+\{(a,a^+),(b,b^+),(c,c^+),(x,x^+),(y,y^+),(z,z^+)\}-(x,y)$ is a H-path of $BH_{n}-F$ passing through $L$.

{\it Case 2.2.3.}  $i=2,j=3$.

By Lemma \ref{le-2}, there is a $c\in V_1\cap X\setminus \{s\}$, such that $c$ (resp. $c^{+}$) is incident with none of $E(L_1)$ (resp. $E(L_{2})$). By the induction hypothesis, $B^{1}-F_{1}$, $B^{2}-F_{2}$, $B^{3}-F_{3}$ have H-paths $P[a^{+},c]$, $P[c^{+},u]$, $P[b^{+},v]$ passing through $L_1$, $L_2$ and $L_{3}$, respectively. Hence, $P[a,b]\cup P[a^{+},c]\cup P[c^{+},u]\cup P[b^{+},v]+\{(a,a^{+}),(b,b^{+}),(c,c^{+})\}$ is a H-path of $BH_{n}-F$ passing through $L$.
\end{proof}

\begin{lemma}
If $|E(L_{0})\cup F_{0}|=2n-3$, then $BH_{n}-F$ contains a H-path $P[u,v]$ passing through $L$.
\end{lemma}

\begin{proof}
In this scenario, $E(L_{k})\cup F_{k}=\emptyset$ for $k\in N_{4}\setminus \{0\}$.

{\it Case 1.}  $u,v\in V_i$.

{\it Case 1.1.}  $i=0$.

Since $\{u,v\}$ is compatible to $L$ and $E(L_0)\neq \emptyset$, there is a path in $L_0$ such that at least one of the two end vertices, say $x$, is not in $\{u,v\}$. Without loss of generality, assume that $x\in X$. Let $(x,y)\in E(L_0)$. By the induction hypothesis, $B^{0}-F_{0}$ has a H-path $P[u,v]$ passing through $L_{0}-(x,y)$. Let $c\in V_1\cap X\setminus \{s\}$, $d\in V_2\cap X\setminus \{s\}$.

Suppose first that $(x,y)\in E(P[u,v])$. Let $(a,b)$ be an arbitrary edge in $P[u,v]\setminus E(L_0)$ for some $a\in X$ and $b\in Y$. Since $|F^c|=1$, $(a,a^+)$ or $(a,a^-)$ (resp. $(b,b^+)$ or $(b,b^-)$), say $(a,a^+)$ (resp. $(b,b^+)$), is not in $F^c$. By Theorem \ref{th-xu2007}, $B^1$, $B^2$, $B^3$ have H-paths $P[a^+,c]$, $P[c^+,d]$ and $P[b^+,d^+]$, respectively. Thus, $P[u,v]\cup P[a^+,c]\cup P[c^+,d]\cup P[b^+,d^+\{(a,a^+),(b,b^+),(c,c^+),(d,d^+)\}-(a,b)$ is a H-path of $BH_{n}-F$ passing through $L$.

Suppose now that $(x,y)\notin E(P[u,v])$. No matter $b$ is $v$ or not, there is a neighbor $x$ of $b$ on $P[u,v]$ such that $(b,x)\notin E(L_0)$. Let $(a,y)\in E(P[u,v])$ such that exactly one of $\{x,y\}$ lies on the segment of $P[u,v]$ between $a$ and $b$. By Theorem \ref{th-xu2007}, $B^1$, $B^2$, $B^3$ have H-paths $P[x^+,c]$, $P[c^+,d]$ and $P[y^+,d^+]$, respectively. Thus, $P[u,v]\cup P[x^+,c]\cup P[c^+,d]\cup P[y^+,d^+]+\{(a,b),(c,c^+),(d,d^+),(x,x^+),(y,y^+)\}-\{(a,y),(b,x)\}$ is a H-path of $BH_{n}-F$ passing through $L$.

{\it Case 1.2.}  $i\neq 0$.

By Theorem \ref{th-li2019} and Lemma \ref{le-1}, $B^0-F_0$ has a H-cycle $C_0$ passing through $L_0$. Let $(a,b)\in E(C_0)\setminus E(L_0)$ for some $a\in X$ and $b\in Y$ such that $\{a,b\}\cap \{s,s^+\}=\emptyset$. Thus, $P[a,b]=C_0-(a,b)$ a H-path passing through $L_0$ of $BH_n-F$.

{\it Case 1.2.1.}  $i=1$.

By Theorem \ref{th-xu2007}, $B^1-F_1$ has a H-path $P[u,v]$. Let $(a^+,c)\in E(P[u,v])$. Since $|F^c|=1$, $(c,c^+)$ or $(c,c^-)$, say $(c,c^+)$, is not in $F^c$. Let $d\in V_2\cap X\setminus \{s\}$. By Theorem \ref{th-xu2007}, $B^2-F_2$, $B^3-F_3$ have H-paths $P[c^+,d]$ and $P[b^+,d^+]$, respectively. Thus, $P[a,b]\cup P[u,v]\cup P[c^+,d]\cup P[b^+,d^+]+\{(a,a^{+}),(b,b^{+}),(c,c^{+}),(d,d^{+})\}-(a^{+},c)$ is a H-path of $BH_{n}-F$ passing through $L$.


{\it Case 1.2.2.}  $i=2$.

Let $c\in V_1\cap X\setminus \{s\}$. By Theorem \ref{th-xu2007}, $B^1-F_1$, $B^2-F_2$ have H-paths $P[a^+,c]$ and $P[u,v]$, respectively. Let $(c^+,d)\in E(P[u,v])$. Since $|F^c|=1$, $(d,d^+)$ or $(d,d^-)$, say $(d,d^+)$, is not in $F^c$. By Theorem \ref{th-xu2007}, $B^3-F_3$ has a H-path $P[b^+,d^+]$. Thus, $P[a,b]\cup P[a^+,c]\cup P[u,v]\cup P[b^+,d^+]+\{(a,a^{+}),(b,b^{+}),(c,c^{+}),(d,d^{+})\}-(c^{+},d)$ is a H-path of $BH_{n}-F$ passing through $L$.


{\it Case 1.2.3.}  $i=3$.

By Theorem \ref{th-xu2007}, $B^3-F_3$ has a H-path $P[u,v]$. Let $(b^+,d)\in E(P[u,v])$. Since $|F^c|=1$, $(d,d^+)$ or $(d,d^-)$, say $(d,d^+)$, is not in $F^c$. Let $c\in V_1\cap X\setminus \{s\}$. By Theorem \ref{th-xu2007}, $B^1-F_1$, $B^2-F_2$ have H-paths $P[a^+,c]$ and $P[c^+,d^+]$, respectively. Thus, $P[a,b]\cup P[a^+,c]\cup P[c^+,d^+]\cup P[u,v]+\{(a,a^{+}),(b,b^{+}),(c,c^{+}),(d,d^{+})\}-(b^{+},d)$ is a H-path of $BH_{n}-F$ passing through $L$.


{\it Case 2.}  $u\in V_i$ and $v\in V_j$ for $i,j\in N_{4}$ and $i\neq j$.

{\it Case 2.1.}  $i=0$.

By Theorem \ref{th-li2019} and Lemma \ref{le-1}, $B^0-F_0$ has a H-cycle $C_0$ passing through $L_0$. Let $(u,a)\in E(C_0)\setminus E(L_0)$. Thus, $P[u,a]=C_0-(u,a)$ a H-path passing through $L_0$ of $BH_n-F$. Since $|F^c|=1$, $(a,a^+)$ or $(a,a^-)$, say $(a,a^+)$, is not in $F^c$.


{\it Case 2.1.1.}  $j=1$.

Let $b\in V_1\cap X\setminus \{s\}$ and $c\in V_2\cap X\setminus \{s\}$. By Theorem \ref{th-xu2007}, $B^1-F_1$, $B^2-F_2$, $B^3-F_3$ have H-paths $P[v,b]$, $P[b^+,c]$ and $P[a^+,c^+]$, respectively. Thus, $P[u,a]\cup P[v,b]\cup P[b^+,c]\cup P[a^+,c^+]+\{(a,a^{+}),(b,b^{+}),(c,c^{+})\}$ is a H-path of $BH_{n}-F$ passing through $L$.


{\it Case 2.1.2.}  $j=2$.

By Lemma \ref{le-3}, there is an edge $(x,y)\in E(P[u,a])\setminus E(L_0)$ such that $\{x,y\}\cap \{u,a\}=\emptyset$ and $\{x,y\}\cap \{s,s^+\}=\emptyset$. Let $z\in V_1\cap X\setminus \{s\}$ and $w\in V_2\cap X\setminus \{s\}$. By Theorem \ref{th-xu2007}, $B^1-F_1$, $B^2-F_2$ have H-paths $P[x^+,z]$ and $P[v,w]$, respectively. There is a neighbor of $z$ in $B^2$, say $z^+$, being not $v$. Let $c$ be the neighbor of $z^+$ on the segment of $P[v,w]$ between $z^+$ and $v$. By Theorem \ref{th-cheng2014}, there are two vertex-disjoint paths $P[a^+,c^+]$ and $P[y^+,w^+]$ in $B^3$ such that each vertex of $B^3$ lies on one of the two paths. Thus, $P[u,a]\cup P[x^+,z]\cup P[v,w]\cup P[a^+,c^+]\cup P[y^+,w^+]+\{(a,a^{+}),(c,c^{+}),(w,w^{+}),(x,x^{+}),(y,y^{+}),(z,z^{+})\}-\{(x,y),(z^{+},c)\}$ is a H-path of $BH_{n}-F$ passing through $L$.


{\it Case 2.1.3.}  $i=3$.

By Lemma \ref{le-3}, there is an edge $(x,y)\in E(P[u,a])\setminus E(L_0)$ such that $\{x,y\}\cap \{u,a\}=\emptyset$ and $\{x,y\}\cap \{s,s^+\}=\emptyset$. Let $z\in V_1\cap X\setminus \{s\}$ and $w\in V_2\cap X\setminus \{s\}$. By Theorem \ref{th-xu2007}, $B^1-F_1$, $B^2-F_2$ have H-paths $P[x^+,z]$ and $P[z^+,w]$, respectively. There is a neighbor of $w$ in $B^3$, say $w^+$, being not $v$. By Theorem \ref{th-cheng2014}, there are two vertex-disjoint paths $P[a^+,v]$ and $P[y^+,w^+]$ in $B^3$ such that each vertex of $B^3$ lies on one of the two paths. Thus, $P[u,a]\cup P[x^+,z]\cup P[z^+,w]\cup P[a^+,v]\cup P[y^+,w^+]+\{(a,a^{+}),(w,w^{+}),(x,x^{+}),(y,y^{+}),(z,z^{+})\}-(x,y)$ is a H-path of $BH_{n}-F$ passing through $L$.


{\it Case 2.2.} $i\neq 0$.

Let $(a,b)\in E(C_0)\setminus E(L_0)$ such that $\{a,b\}\cap \{s,s^+\}=\emptyset$. Then $P[a,b]=C_0-(a,b)$ is a H-path passing through $L_0$ of $B^0-F_0$.

{\it Case 2.2.1.}  $i=1,j=2$.

Let $c\in V_2\cap X\setminus \{s\}$. By Theorem \ref{th-xu2007}, $B^1-F_1$, $B^2-F_2$, $B^3-F_3$ have H-paths $P[a^+,u]$, $P[v,c]$ and $P[b^+,c^+]$, respectively. Thus, $P[a,b]\cup P[a^+,u]\cup P[v,c]\cup P[b^+,c^+]+\{(a,a^{+}),(b,b^{+}),(c,c^{+})\}$ is a H-path of $BH_{n}-F$ passing through $L$.

{\it Case 2.2.2.}  $i=1,j=3$.

By Lemma \ref{le-3}, there is an edge $(x,y)\in E(P[a,b])\setminus E(L_0)$ such that $\{x,y\}\cap \{a,b\}=\emptyset$ and $\{x,y\}\cap \{s,s^+\}=\emptyset$. Let $z\in V_2\cap Y\setminus \{s^+\}$ and $w\in V_2\cap X\setminus \{s\}$. By Theorem \ref{th-xu2007}, $B^2$ has a H-path $P[z,w]$. There is a neighbor of $z$ (resp. $w$) in $B^1$ (resp. $B^3$), say $z^+$ (resp. $w^+$), being not $u$ (resp. $v$). By Theorem \ref{th-cheng2014}, there are two vertex-disjoint paths $P[a^+,u]$ and $P[x^+,z^+]$ (resp. $P[b^+,v]$ and $P[y^+,w^+]$) in $B^1$ (resp. $B^3$) such that each vertex of $B^1$ (resp. $B^3$) lies on one of the two paths. Thus, $P[a,b]\cup P[x^+,z^+]\cup P[a^+,u]\cup P[z,w]\cup P[b^+,v]\cup P[y^+,w^+]+\{(a,a^{+}),(b,b^+),(w,w^{+}),(x,x^{+}),(y,y^{+}),(z,z^{+})\}-(x,y)$ is a H-path of $BH_{n}-F$ passing through $L$.

{\it Case 2.2.3.}  $i=2,j=3$.

Let $c\in V_1\cap X\setminus \{s\}$. By Theorem \ref{th-xu2007}, $B^1-F_1$, $B^2-F_2$, $B^3-F_3$ have H-paths $P[a^+,c]$, $P[c^+,u]$ and $P[b^+,v]$, respectively. Thus, $P[a,b]\cup P[a^+,c]\cup P[c^+,u]\cup P[b^+,v]+\{(a,a^{+}),(b,b^{+}),(c,c^{+})\}$ is a H-path of $BH_{n}-F$ passing through $L$.
\end{proof}

\section{Conclusions}

Let $F\subset BH_n$ and let $L$ be a linear forest of $BH_n-F$ such that $|F|+|E(L)|\leq 2n-2$. For any two vertices $u$ and $v$ of opposite parts in $BH_n$ that are compatible to $L$, we bent to show that there is a hamiltonian path of $BH_n-F$ between $u$ and $v$ passing through $L$. The proof was carried out by induction on $n$. Some known results indicates the assertation holds for the base case $n=2$. Assume the assertation holds for $n-1$ and prove it also holds for $n$ with $n\geq 3$. If $|F|=2n-3$ and the lins of $F$ are incident with a common node, then we choose some dimension such that $F$ has at least two links in this dimension and $L$ has no link in this dimension; Otherwise, we choose some dimension such that the total number of $F$ and $L$ in this dimension does not exceed $1$. 
No matter which case, without loss of generality, assume that the chosen dimension is dimension $n-1$. Partition $BH_n$ into $4$ disjoint copies of $BH_{n-1}$ along dimension $n-1$. On the basis of the above partition of $BH_n$, we complete the proof for the case that there is at most $1$ faulty link in dimension $n-1$. According to the classification method, the case that $F$ has exactly two links in dimension $n-1$ was solved in Section 4 of \cite{Yang-Song}.

An interesting related problem is to investigate the fault-tolerant-prescribed hamiltonian laceability of balanced hypercubes in the hybrid faluts model.

\end{document}